\documentclass{gtart_h}  

\def\ifplaintex{\expandafter\ifx\csname documentclass\endcsname\relax}

\def\ifplaintex{\expandafter\ifx\csname documentclass\endcsname\relax}

\ifplaintex 
\else
\fi

\expandafter\ifx\csname epsfbox\endcsname\relax\input epsf\fi

\def\gt{{\mathsurround=0pt\it $\cal G\mskip-2mu$eometry \&\ 
$\cal T\!\!$opology}}        %  journal title in recommended style

\def\gtp{{\mathsurround=0pt\it $\cal G\mskip-2mu$eometry \&\ 
$\cal T\!\!$opology $\cal P\!$ublications}}  % GT publications

\def\lognumber#1{\def\thelognumber{#1}}
\def\volumenumber#1{\def\thevolumenumber{#1}}
\def\papernumber#1{\def\thepapernumber{#1}}
\def\volumeyear#1{\def\thevolumeyear{#1}}

\def\pagenumbers#1#2{\def\startpage{#1}\def\finishpage{#2}}
\def\published#1{\def\publishdate{#1}}
\def\proposed#1{\def\theproposer{#1}}
\def\seconded#1{\def\theseconders{#1}}
\def\received#1{\def\receiveddate{#1}}

\def\accepted#1{\def\accepteddate{#1}}

\long\def\asciiabstract#1{\long\def\theasciiabstract{#1}}
\def\asciikeywords#1{\def\theasciikeywords{#1}}

\let\\\par\let\thelognumber\relax
\let\thevolumenumber\relax\let\thepapernumber\relax
\let\thevolumeyear\relax\let\thesamplenumber\relax\let\startpage\relax
\let\finishpage\relax\let\publishdate\relax\let\receiveddate\relax
\let\reviseddate\relax\let\accepteddate\relax\let\theasciititle\relax
\let\theasciiauthors\relax
\let\theasciiabstract\relax\let\theasciikeywords\relax
\let\theasciiemail\relax\let\theshortauthors\relax\let\theshorttitle\relax

\long\def\maketitlep{   % start of definition of \maketitlep

\count0=\startpage

\gt\hfill      %   Journal title (top left) 
\hbox to 77pt{\vbox to 0pt{\vglue -15pt\epsfbox{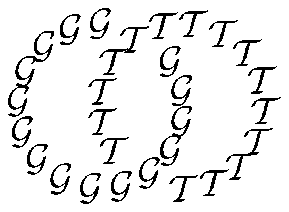}\vss}\hss}
\break
{\small\ifx\thesamplenumber\relax % sample?  
Volume \else Sample
\fi\thevolumenumber\ (\thevolumeyear)
\startpage--\finishpage\nl
Published: \publishdate}
\vglue 0.5truein plus 0.4fil minus 0.1truein

{\parskip=0pt\leftskip 0pt plus 1fil\def\\{\par\smallskip}{\ifplaintex\large
\else\Large\fi\bf\thetitle}\par\medskip}   

\vglue 0pt plus 0.1fil 

{\parskip=0pt\leftskip 0pt plus 1fil\def\\{\par}{\sc\theauthors}
\par\medskip}

\vglue 0pt plus 0.1fil 

{\small\parskip=0pt\let\newline\\
{\leftskip 0pt plus 1fil\def\\{\par}{\sl\theaddress}\par}
\expandafter\ifx\theemail\relax    % email address?
\relax\else\vglue 5pt plus 0.02fil minus 2pt\def\\{\stdspace{\rm 
and}\stdspace} 
\cl{Email:\stdspace\tt\theemail}\fi
\ifx\theurl\relax                  % URL given?
\relax\else\vglue 5pt plus 0.02fil minus 2pt\def\\{\stdspace{\rm 
and}\stdspace}
\cl{URL:\stdspace\tt\theurl}\fi\par}

\vglue 7pt plus 0.3fil minus 3pt

{\bf Abstract}
\vglue 5pt plus 0.1fil minus 2pt

\theabstract

\vglue 7pt plus 0.3fil minus 3pt

{\bf AMS Classification numbers}\quad Primary:\quad \theprimaryclass

Secondary:\quad \thesecondaryclass

\vglue 5pt plus 0.3fil minus 2pt

{\bf Keywords:}\quad \thekeywords

\vglue 10pt plus 0.5fil minus 5pt

{\small  Proposed: \theproposer\hfill Received: \receiveddate\nl
Seconded: \theseconders\hfill 
\ifx\reviseddate\relax                         % paper revised?
Accepted: \accepteddate                        % no
\else
Revised: \reviseddate                          % yes
\fi}
\eject
}       %  end of definition of \maketitlep

\font\phead=cmsl9 scaled 950
\font=cmsl9 scaled 1050
\font\pnum=cmbx10 scaled 913
\font\lnum=cmbx10 
\font\pfoot=cmsl9 scaled 950
\font=cmsl9 scaled 1050
\ifplaintex
\headline{\vbox to 0pt{\vskip -4.5mm\line{\small\phead\ifnum
\count0=\startpage ISSN 1364-0380 (on line)
1465-3060 (printed) \hfill {\pnum\folio}\else\ifodd\count0\def\\{ }% 
\ifx\theshorttitle\relax\thetitle\else\theshorttitle\fi\hfill{\pnum\folio}
\else\def\\{ and }{\pnum\folio}\hfill\ifx\theshortauthors\relax\theauthors
\else\theshortauthors\fi\fi\fi}\vss}}
\footline{\vbox to 0pt{\vglue 0mm\line{\small\pfoot\ifnum\count0=\startpage
\copyright\ \gtp\hfill\else
\gt, Volume \thevolumenumber\ (\thevolumeyear)\hfill\fi}\vss
}}
\else
\makeatletter
\def\@oddhead{{\small\lhead\ifnum\count0=\startpage ISSN 1364-0380 (on line)
1465-3060 (printed) \hfill {\lnum\number\count0}\else\ifodd\count0
\def\\{ }\ifx\theshorttitle\relax \thetitle \else\theshorttitle\fi\hfill
{\lnum\number\count0}\else\def\\{ and }{\lnum\number\count0}
\hfill\ifx\theshortauthors\relax 
\theauthors\else\theshortauthors\fi\fi\fi}}\def\@evenhead{\@oddhead}
\def\@oddfoot{\small\lfoot\ifnum\count0=\startpage\copyright\ \gtp\hfill\else
\gt, Volume \thevolumenumber\ (\thevolumeyear)\hfill\fi}
\def\@evenfoot{\@oddfoot}
\makeatother
\fi

\newwrite\gtoutfile
\long\gdef\makeheadfile{  %%% start of definition of \makeheadfile
{\def\\{, }\def\s{ }
\immediate\openout\gtoutfile head.xxx
\immediate\write\gtoutfile{Proxy-for: \ifx\theasciiauthors\relax
\theauthors\else\theasciiauthors\fi\s<\ifx\theasciiemail\relax\theemail\else\theasciiemail\fi>}
\immediate\write\gtoutfile{\noexpand\\}
\immediate\write\gtoutfile{Authors: \ifx\theasciiauthors\relax
\theauthors\else\theasciiauthors\fi}
{\def\\{ }\immediate\write\gtoutfile{Title: \ifx\theasciititle\relax
\thetitle\else\theasciititle\fi}}
\immediate\write\gtoutfile{Subj-class: GT or SG or MG etc}
\immediate\write\gtoutfile{MSC-class: \theprimaryclass\ifx\thesecondaryclass\relax\else, \thesecondaryclass\fi}
\immediate\write\gtoutfile{Journal-ref: Geom. Topol. \thevolumenumber
(\thevolumeyear) \startpage-\finishpage}
\immediate\write\gtoutfile{Comments: Published by Geometry and Topology at}
\immediate\write\gtoutfile{\s\s http://www.maths.warwick.ac.uk/gt/GTVol\thevolumenumber/paper\thepapernumber.abs.html}
\immediate\write\gtoutfile{\noexpand\\}
\immediate\write\gtoutfile{}
\ifx\theasciiabstract\relax
\immediate\write\gtoutfile{\theabstract}\else
\immediate\write\gtoutfile{\theasciiabstract}\fi
\immediate\write\gtoutfile{}
\immediate\write\gtoutfile{\noexpand\\}
\immediate\write\gtoutfile{}
\immediate\closeout\gtoutfile}}  %%% end of definition of \makeheadfile

\def\maketitlepage{\maketitlep\makeheadfile}
\let\maketitle\maketitlepage

\lognumber{512}
\received{3 November 2004}
\volumenumber{9}\papernumber{33}\volumeyear{2005}
\pagenumbers{1443}{1499}   
\published{8 August 2005}
\accepted{04 July 2005}
\proposed{Vaughan Jones}
\seconded{Robion Kirby, Cameron Gordon}

\usepackage{
  amsmath,graphicx,dbnsymb,amssymb,color,picins,hyphenat,
  multicol,slashbox,rotating,floatflt,
}
\usepackage[all]{xy}

\def\figref#1{\hyperlink{#1anchor}{Figure~\ref*{#1}}}
\def\fref#1{\hyperlink{#1anchor}{\ref*{#1}}}
\def\anchor#1{\noindent\hypertarget{#1anchor}{\smash{$\phantom{99}$}}}

\theoremstyle{plain}

\newtheorem{theorem}{Theorem}
\newtheorem{proposition}{Proposition}[section]

\newtheorem{lemma}[proposition]{Lemma}

\newtheorem{conjecture}{Conjecture}

\theoremstyle{definition}
\newtheorem{definition}[proposition]{Definition}

\newtheorem{problem}[proposition]{Problem}

\theoremstyle{remark}
\newtheorem{example}[proposition]{Example}
\newtheorem{exercise}[proposition]{Exercise}
\newtheorem{hint}[proposition]{Hint}
\newtheorem{remark}[proposition]{Remark}

\newlength{\standardunitlength}
\newcommand{\tr}{\operatorname{tr}}

\newenvironment{myitemize}{\bgroup
        \begin{list}{$\bullet$}{\setlength{\leftmargin}{16pt}
        \setlength{\labelwidth}{12pt}
        \setlength{\labelsep}{4pt}}
}{
        \end{list}\egroup
}

\def\llbracket{\left[\!\!\left[}
\def\rrbracket{\right]\!\!\right]}

\def\bbF{{\mathbb F}}
\def\bbQ{{\mathbb Q}}
\def\bbR{{\mathbb R}}
\def\bbZ{{\mathbb Z}}
\def\hatJ{{\hat J}}
\def\calA{{\mathcal A}}
\def\calC{{\mathcal C}}

\def\calF{{\mathcal F}}
\def\calG{{\mathcal G}}
\def\calO{{\mathcal O}}
\def\calP{{\mathcal P}}
\def\calS{{\mathcal S}}
\def\calT{{\mathcal T}}

\newcommand{\Alg}{{\mathcal Alg}}
\newcommand{\FourTu}{{\text{\it 4Tu}}}
\newcommand{\Cob}{{\mathcal Cob}}
\newcommand{\Cobd}{{\mathcal Cob}_\bullet}
\newcommand{\Cobdl}{{\mathcal Cob}_{\bullet/l}}
\newcommand{\Cobi}{{\mathcal Cob}_{/i}}
\newcommand{\Cobl}{{\mathcal Cob}_{/l}}
\newcommand{\Cobz}{{\mathcal Cob}_{0/l}}
\newcommand{\Kh}{{\text{\it Kh}}}
\newcommand{\qdim}{\operatorname{{\it q}dim}}
\newcommand{\Kob}{\operatorname{Kob}}
\newcommand{\Kobh}{{\operatorname{Kob}_{/h}}}
\newcommand{\Kom}{\operatorname{Kom}}
\newcommand{\Komh}{\operatorname{Kom}_{/h}}
\newcommand{\Mat}{\operatorname{Mat}}
\newcommand{\MM}{\operatorname{MM}}
\newcommand{\Mor}{\operatorname{Mor}}
\newcommand{\Obj}{\operatorname{Obj}}
\newcommand{\PKob}{\operatorname{Kob_{/\pm}}}
\newcommand{\PKobh}{\operatorname{Kob}_{/\pm h}}
\newcommand{\Top}{{\mathcal T\!op}}
\newcommand{\TX}{{T\!X}}
\newcommand{\Vect}{\operatorname{Vect}}
\newcommand{\ZMod}{{\bbZ\!\operatorname{Mod}}}
\newcommand{\Ztwo}{{\bbZ_{(2)}}}

\newcommand{\eps}[2]{{\hspace{-3pt}\begin{array}{c}%
  \raisebox{-2.5pt}{\includegraphics[width=#1]{figs/#2.eps}}%
\end{array}\hspace{-3pt}}}
\newcommand{\epsg}[2]{{\hspace{-3pt}\begin{array}{c}%
  \raisebox{0pt}{\includegraphics[width=#1]{figs/#2.eps}}%
\end{array}\hspace{-3pt}}}

\begin{document}

\title{Khovanov's homology for tangles and cobordisms}

\author{Dror Bar-Natan}
\address{Department of Mathematics, University of 
Toronto\\Toronto, Ontario M5S 3G3, Canada}
\email{drorbn@math.toronto.edu}
\urladdr{http://www.math.toronto.edu/~drorbn}


\asciikeywords{2-knots, canopoly, categorification, cobordism, Euler
  characteristic, Jones polynomial, Kauffman bracket, Khovanov, knot
  invariants, movie moves, planar algebra, skein modules, tangles,
  trace groups}

\begin{abstract}
  We give a fresh introduction to the Khovanov Homology theory for knots
  and links, with special emphasis on its extension to tangles, cobordisms and
  2--knots. By staying within a world of topological pictures a little
  longer than in other articles on the subject, the required extension
  becomes essentially tautological. And then a simple application of an
  appropriate functor (a ``TQFT'') to our pictures takes them to the 
  familiar realm of complexes of (graded) vector spaces and ordinary
  homological invariants. 
\end{abstract}

\asciiabstract{%
  We give a fresh introduction to the Khovanov Homology theory for knots
  and links, with special emphasis on its extension to tangles, cobordisms and
  2-knots. By staying within a world of topological pictures a little
  longer than in other articles on the subject, the required extension
  becomes essentially tautological. And then a simple application of an
  appropriate functor (a `TQFT') to our pictures takes them to the 
  familiar realm of complexes of (graded) vector spaces and ordinary
  homological invariants.}

\maketitle

\section{Introduction} \label{sec:intro}

The Euler characteristic of a space is a fine invariant; it is a key
ingredient in any discussion of the topology of surfaces (indeed, it
separates closed orientable surfaces), and it has further uses in higher
dimensions as well. But homology is way better. The alternating sum of the
ranks of the homology groups of a space is its Euler characteristic, so
homology is at least as strong an invariant, and it is easy to find
examples showing that homology is a strictly stronger invariant than the
Euler characteristic.

And then there's more. Unlike the Euler characteristic, homology is a
{\em functor} --- continuous maps between spaces induce maps between
their homologies, and it is this property of homology that makes it one
of the cornerstones of algebraic topology; not merely the fact that it
is a little better then the Euler characteristic at telling spaces
apart.

In his seminal paper~\cite{Khovanov:Categorification} Khovanov
explained that the Jones polynomial of a link $L$ is the (graded) Euler
characteristic of a certain ``link homology'' theory. In my follow up
paper~\cite{Bar-Natan:Categorification} I have computed the Khovanov
homology of many links and found that indeed it is a stronger invariant
than the Jones polynomial, as perhaps could be expected in the light of
the classical example of the Euler characteristic and the homology of
spaces.

In further analogy with the classical picture,
Jacobsson~\cite{Jacobsson:Cobordisms} and
Khovanov \cite{Khovanov:Cobordisms} found what seems to be the
appropriate functoriality property of the Khovanov homology. What
replaces continuous maps between spaces is 4D cobordisms between
links:  Given such a cobordism $C$ between links $L_1$ and $L_2$
Jacobsson and Khovanov show how to construct a map $\Kh(C)\co
\Kh(L_1)\to \Kh(L_2)$ (defined up to a $\pm$ sign) between the
corresponding Khovanov homology groups $\Kh(L_1)$ and
$\Kh(L_2)$.\footnote{Check~\cite{Rasmussen:SliceGenus} for a
topological application and~\cite{CarterSaitoSatoh:KhovanovJacobsson}
for some computations.} Note that if $L_1$ and $L_2$ are both the empty
link, then a cobordism
between $L_1$ and $L_2$ is a 2--knot in $\bbR^4$ (see
eg,~\cite{CarterSaito:KnottedSurfaces}) and $\Kh(C)\co \Kh(L_1)\to
\Kh(L_2)$ becomes a group homomorphism $\bbZ\to\bbZ$, hence a single
scalar (defined up to a sign). That is, this special case of $\Kh(C)$
yields a numerical invariant of
2--knots.\footnote{Added July 2005: The latter is trivial;
see~\cite{Rasmussen:ClosedSurfaces}.
}

Given a ``movie presentation'' of a cobordism $C$ (as
in~\cite{CarterSaito:KnottedSurfaces}), the construction of $\Kh(C)$ is
quite simple to describe. But the proofs that $\Kh(C)$ is independent of
the specific movie presentation chosen for $C$ are quite involved.
Jacobsson's proof involves a large number of complicated case by case
computations.  Khovanov's proof is more conceptual, but it relies
on his rather complicated ``functor-valued invariant of
tangles''~\cite{Khovanov:Functor} and even then there remains some
case-checking to do.

A major purpose of this article is to rewrite Khovanov's proof in a 
simpler language. Thus a major part of our work is to simplify (and at
the same time extend!) Khovanov's treatment of tangles. As side benefits we
find a new homology theory for knots/links (Section~\ref{subsec:4TuZ2}) and
what we believe is the ``right'' way to see that the Euler characteristic
of Khovanov homology is the Jones polynomial (Section~\ref{sec:Euler}).

\subsection{Executive summary}
This quick section is for the experts. If you aren't one, skip straight to
Section~\ref{subsec:Plan}.

\subsubsection{Why are tangles relevant to cobordisms?}
\label{subsubsec:WhyTangles}
Tangles are knot pieces, cobordisms are movies starring
knots and links~\cite{CarterSaito:KnottedSurfaces}. Why is the former
relevant to the study of the latter (in the context of Khovanov
homology)?

The main difficulty in showing that cobordisms induce maps of homology
groups is to show that trivial movies induce trivial maps on homology.
A typical example of such a trivial movie is the circular clip
\begin{equation} \label{eq:MMDemo}
  \begin{array}{c}
    \includegraphics[height=0.7in]{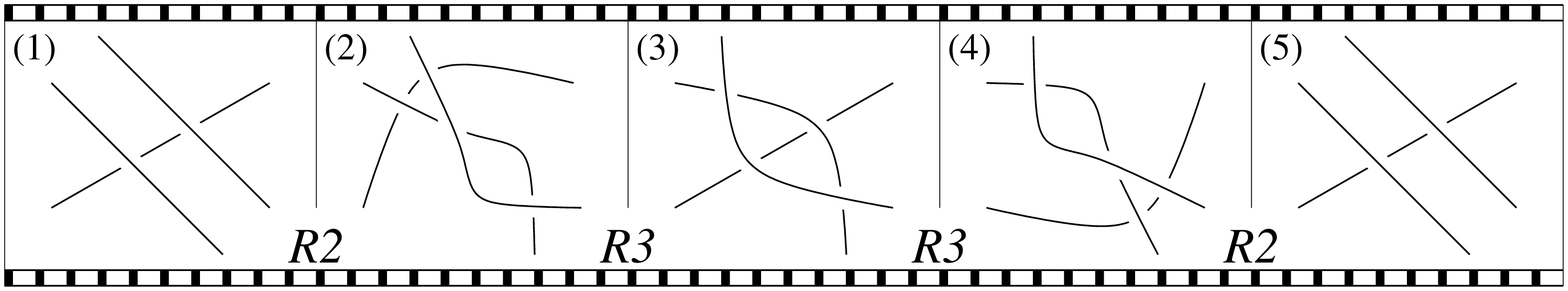}
  \end{array}
\end{equation}
(this is $\MM_6$ of \figref{fig:MM6-10}). Using a nice theory for
tangles which we will develop later, we will be able to replace the
composition of morphisms corresponding to the above clip by the
following composition, whose ``core'' (circled below) remains the same
as in clip~\eqref{eq:MMDemo}:
\begin{equation} \label{eq:SimpleClip}
  \ \hspace{-2mm}\begin{array}{c}
    \includegraphics[width=5.15in]{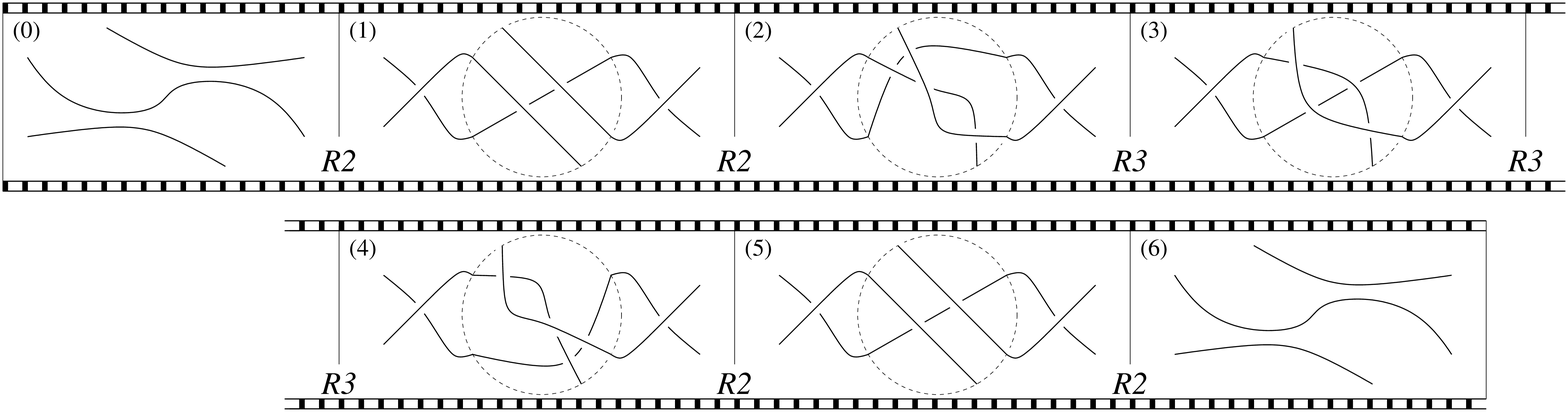}
  \end{array}
\end{equation}
But this composition is an automorphism of the complex $K$ of the
crossingless tangle $T=\begin{array}{c}
\includegraphics[width=0.5in]{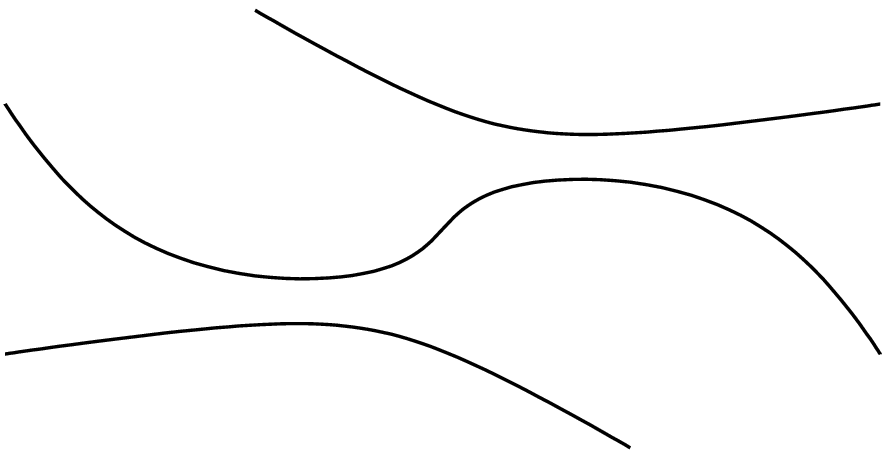} \end{array}$
and that complex is very simple --- as $T$ has no crossings, $K$ ought
to\label{OughtTo}\footnote{We say ``ought to'' only because at this point
our theory of tangles is not yet
defined} consist only of one chain group and no differentials.  With
little luck, this would mean that $K$ is also simple in a technical
sense --- that it has no automorphisms other than multiples of
the identity. Thus indeed the circular clip of
Equation~\eqref{eq:MMDemo} induces a trivial (at least up to a scalar)
map on homology.

In the discussion of the previous paragraph it was crucial that the
complex for the tangle $T=\begin{array}{c}
\includegraphics[width=0.5in]{figs/CrossinglessTangle.eps} \end{array}$ be
simple, and that it would be possible to manipulate tangles as in the
transition from~\eqref{eq:MMDemo} to~\eqref{eq:SimpleClip}. Thus a
``good'' theory of tangles is useful for the study of cobordisms.

\subsubsection{How we deal with tangles} \label{subsubsec:HowTangles}
As defined in~\cite{Khovanov:Categorification} (or
in~\cite{Bar-Natan:Categorification} or
in~\cite{Viro:KhovanovHomology}), the Khovanov homology theory does not
lend itself naturally to an extension to tangles. In order to define
the chain spaces one needs to count the cycles in each smoothing, and
this number is not known unless all `ends' are closed, ie, unless the
tangle is really a link. In~\cite{Khovanov:Functor} Khovanov solves the
problem by taking the chain space of a tangle to be the direct sum of
all chain spaces of all possible closures of that tangle. Apart from
being quite cumbersome (when all the details are in place;
see~\cite{Khovanov:Functor}), as written, Khovanov's solution only
allows for `vertical' compositions of tangles, whereas one would wish
to compose tangles in arbitrary planar ways, in the spirit of V~Jones'
planar algebras~\cite{Jones:PlanarAlgebrasI}.

We deal with the extension to tangles in a different manner. Recall
that the Khovanov picture of~\cite{Khovanov:Categorification} can be
drawn in two steps. First one draws a `topological picture' $\Top$
made of smoothings of a link diagram and of cobordisms between these
smoothings.  Then one applies a certain functor $\calF$ (a
$(1+1)$--dimensional TQFT) to this topological picture, resulting in an
algebraic picture $\Alg$, a complex involving modules and module morphisms
whose homology is shown to be a link invariant. Our trick is to postpone
the application of $\calF$ to a later stage, and prove the invariance
already at the level of the topological picture. To allow for
that, we first need to `mod out' the topological picture by the `kernel' of
the `topology to algebra' functor $\calF$. Fortunately it turns out that
that `kernel' can be described in completely local terms and hence our
construction is completely local and allows for arbitrary compositions. For
the details, read on.

\subsection{The plan} \label{subsec:Plan} A traditional math paper sets
out many formal definitions, states theorems and moves on to proving
them, hoping that a ``picture'' will emerge in the reader's mind as
(s)he struggles to interpret the formal definitions. In our case the
``picture'' can be summarized by a rather fine picture that can be
uploaded to one's mind even without the formalities, and, in fact, the
formalities won't necessarily make the upload any smoother.  Hence we
start our article with the picture, \figref{fig:Main}, and follow
it in Section~\ref{sec:Narrative} by a narrative description thereof,
without yet assigning any meaning to it and without describing the
``frame'' in which it lives --- the category in which it is an object.
We fix that in Sections~\ref{sec:Frame} and~\ref{sec:Invariance}: in
the former we describe a certain category of complexes where our
picture resides, and in the latter we show that within that category
our picture is a homotopy invariant. The nearly tautological
Section~\ref{sec:PlanarAlgebra} discusses the good behaviour of our
invariant under arbitrary tangle compositions.  In
Section~\ref{sec:grading} we refine the picture a bit by introducing
gradings, and in Section~\ref{sec:Homology} we explain that by
applying an appropriate functor $\calF$ (a $1+1$--dimensional TQFT) we
can get a computable homology theory which yields honest knot/link
invariants.

While not the technical heart of this paper,
Sections~\ref{sec:EmbeddedCobordisms}--\ref{sec:Euler} are its raison
d'\^etre. In Section~\ref{sec:EmbeddedCobordisms} we explain how our
machinery allows for a simple and conceptual explanation of the
functoriality of the Khovanov homology under tangle cobordisms. In
Section~\ref{sec:MoreOnCobl} we further discuss the ``frame'' of
Section~\ref{sec:Frame} finding that in the case of closed tangles
(ie, knots and links) and over rings that contain $\frac12$ it frames
very little beyond the original Khovanov homology while if $2$ is not
invertible our frame appears richer than the original. In
Section~\ref{sec:Euler} we introduce a generalized notion of Euler
characteristic which allows us to ``localize'' the assertion ``The
Euler characteristic of Khovanov Homology is the Jones polynomial''.

The final Section~\ref{sec:OddsAndEnds} contains some further ``odds
and ends''.

\subsection{Acknowledgement} I wish to thank B~Chorny, L~Kauffman,
M~Khovanov, A~Kricker, G~Naot, J~Przytycki, A~Referee, J~Roberts,
S\,D~Schack, K~Shan, A~Sikora, A~Shumakovich, J~Stasheff,
D~Thurston and Y~Yokota for their comments and suggestions and
S~Carter and M~Saito for allowing me to use some figures
from~\cite{CarterSaito:KnottedSurfaces}. This work was partially supported
by NSERC grant RGPIN 262178.

\section{A picture's worth a thousand words} \label{sec:Narrative}

\begin{figure}[p]\anchor{fig:Main}
\begin{center}\begin{sideways}
  \input figs/Main.pstex_t
  %\[ \includegraphics[width=6.4in]{figs/Main.eps} \]
\end{sideways}\end{center}
\caption{
  The main picture. See the narrative in Section~\ref{sec:Narrative}.
} \label{fig:Main}
\end{figure}

As promised in the introduction --- we'd like to start with
\figref{fig:Main} at a completely descriptive level. All
interpretations will be postponed to later sections.

\parpic[r]{$\begin{array}{c}
  \includegraphics[width=1in]{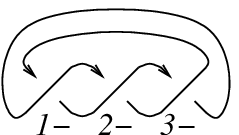}
\end{array}$}
\medskip{\bf2.1\qua Knot}\qua
On the upper left of the figure we see the left-handed trefoil knot $K$
with its $n=3$ crossings labeled $1$, $2$ and $3$. It is enclosed in
double brackets ($\llbracket\cdot\rrbracket$) to indicate that the
rest of the figure shows the {\em formal Khovanov Bracket} of the
left-handed trefoil. As we describe the rest of the figure we will also
indicate how it changes if the left-handed trefoil is replaced by an
arbitrary other knot or link.

\parpic[r]{$\begin{array}{c}
  \includegraphics[width=0.6in]{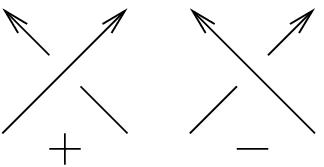}
\end{array}$}
\medskip{\bf2.2\qua Crossings}\qua
On the figure of $K$ we have also marked the signs of its crossings --- $(+)$
for overcrossings ($\overcrossing$) and $(-)$ for undercrossings
($\undercrossing$). Let $n_+$ and $n_-$ be the numbers of $(+)$ crossings
and $(-)$ crossings in $K$, respectively. Thus for the left-handed trefoil
knot, $(n_+, n_-)=(0,3)$.

\parpic[r]{$\begin{array}{c}
  \includegraphics[width=1.5in]{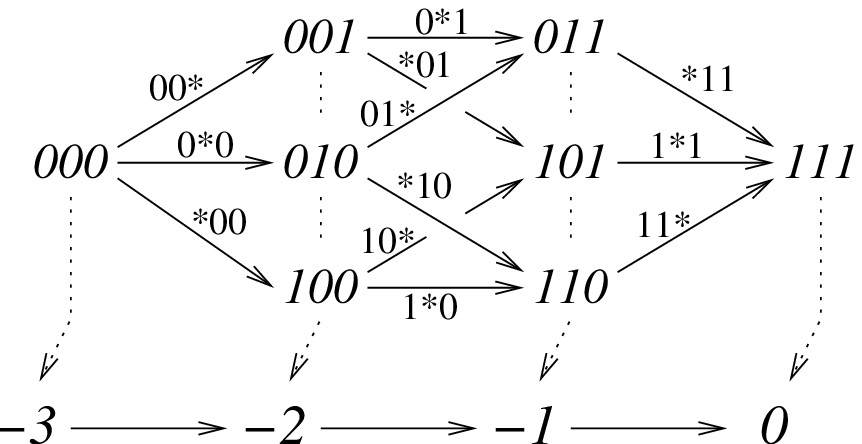}
\end{array}$}
\medskip{\bf2.3\qua Cube}\qua
The main part of the figure is the 3--dimensional cube whose vertices
are all the 3--letter strings of $0$'s and $1$'s. The edges of the
cube are marked in the natural manner by 3--letter strings of $0$'s,
$1$'s and precisely one $\star$ (the $\star$ denotes the coordinate which
changes from $0$ to $1$ along a given edge). The cube is skewered along its
main diagonal, from $000$ to $111$. More precisely, each vertex of the
cube has a ``height'', the sum of its coordinates, a number between $0$
and $3$. The cube is displayed in such a way so that vertices of height
$k$ project down to the point $k-n_-$ on a line marked below the cube. We've
indicated these projections with dashed arrows and tilted them a bit to
remind us of the $-n_-$ shift.

\parpic[r]{\raisebox{-24mm}{$\begin{array}{c}
  \includegraphics[width=2in]{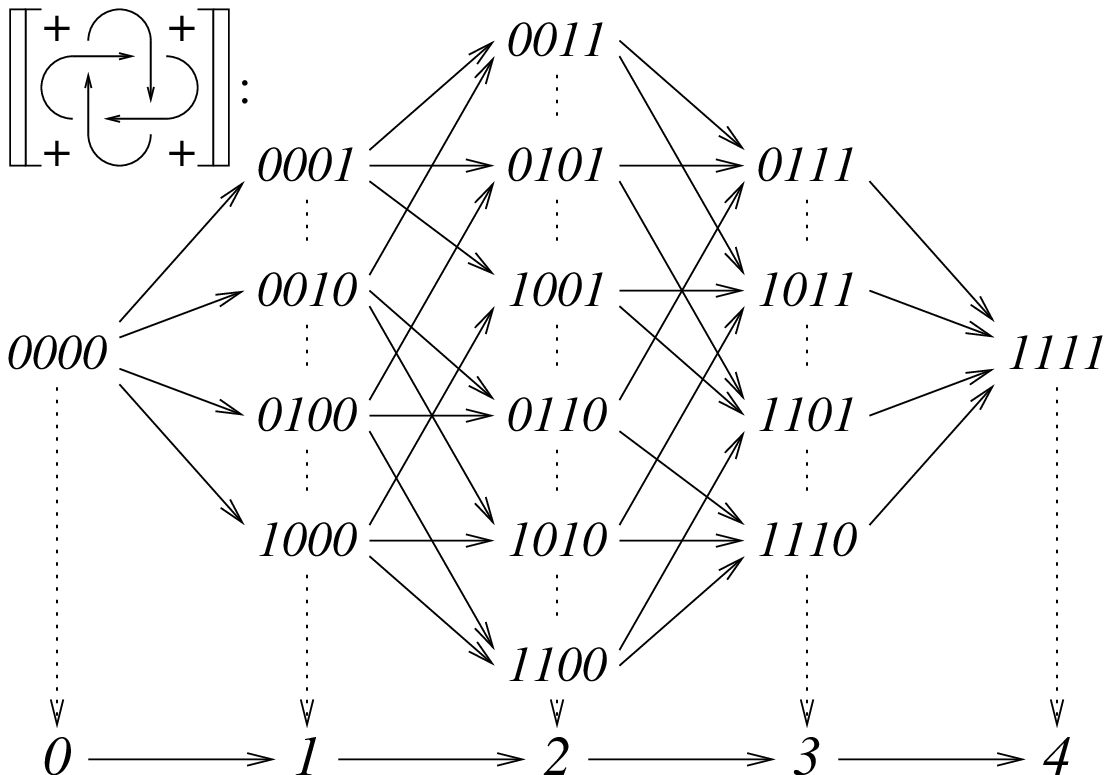}
\end{array}$}}
\vskip 4mm
{\bf2.4\qua Aside: More crossings}\qua
Had we been talking about some $n$ crossing knot rather than the
$3$--crossing left-handed trefoil, the core of our picture would have
been the $n$\hyp dimensional cube with vertices $\{0,1\}^n$, projected by
the ``shifted height'' to the integer points on the interval $[-n_-,n_+]$.
\vskip 4mm

\begin{figure}\anchor{fig:Smoothings}
\[ \eps{4in}{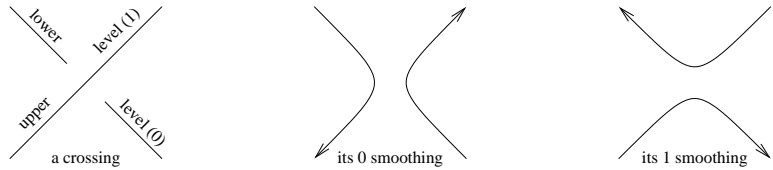} \]
\caption{
  A crossing is an interchange involving two highways. The
  $0$--smoothing is when you enter on the lower level (level $0$) and
  turn right at the crossing. The $1$--smoothing is when you enter on
  the upper level (level $1$) and turn right at the crossing.
} \label{fig:Smoothings}
\end{figure}

\parpic[r]{$\begin{array}{c}
  \includegraphics[width=1in]{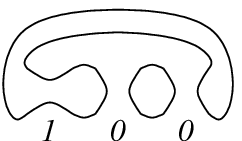}
\end{array}$}
{\bf2.5\qua Vertices}\qua
Each vertex of the cube carries a {\em smoothing} of $K$ --- a planar
diagram obtained by replacing every crossing $\slashoverback$ in the given
diagram of $K$ with either a ``$0$--smoothing'' ($\smoothing$) or with a
``$1$--smoothing'' ($\hsmoothing$) (see \figref{fig:Smoothings} for the
distinction). As our $K$ has $3$ crossings, it has
$2^3=8$ smoothings. Given the ordering on the crossings of $K$ these $8$
smoothings naturally correspond to the vertices of the 3--dimensional cube
$\{0,1\}^3$.

\parpic[r]{$\begin{array}{c}
  \includegraphics[width=1.9in]{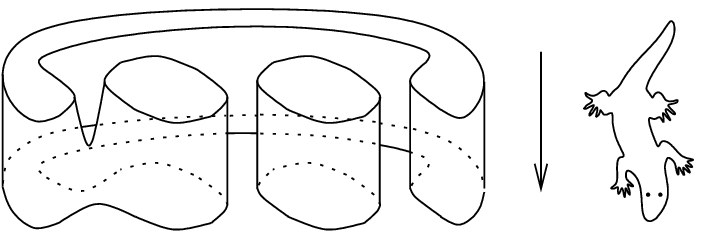}
\end{array}$}
{\bf2.6\qua Edges}\qua
Each edge of the cube is labeled by a {\em cobordism} between the
smoothing on the tail of that edge and the smoothing on its head ---
an oriented two dimensional surfaces embedded in $\bbR^2\times[0,1]$
whose boundary lies entirely in $\bbR^2\times\{0,1\}$ and whose ``top''
boundary is the ``tail'' smoothing and whose ``bottom'' boundary
is the ``head'' smoothing. Specifically, to get the cobordism for an edge
$(\xi_i)\in\{0,1,\star\}^3$ for which $\xi_j=\star$ we remove a disk
neighborhood of the crossing $j$ from the smoothing
$\xi(0):=\xi|_{\star\to 0}$ of $K$, cross with $[0,1]$, and fill the
empty cylindrical slot around the missing crossing with a saddle cobordism
$\begin{array}{c}
  \includegraphics[height=24pt]{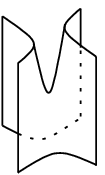}
\end{array}$.
Only one such cobordism is displayed in full in \figref{fig:Main}
--- the one corresponding to the edge $\star 00$. The other $11$
cobordisms are only shown in a diagrammatic form, where the diagram-piece
$\HSaddleSymbol$ stands for the saddle cobordism with top $\smoothing$ and
bottom $\hsmoothing$.

\parpic[r]{$\begin{array}{c}
  \includegraphics[width=2.5in]{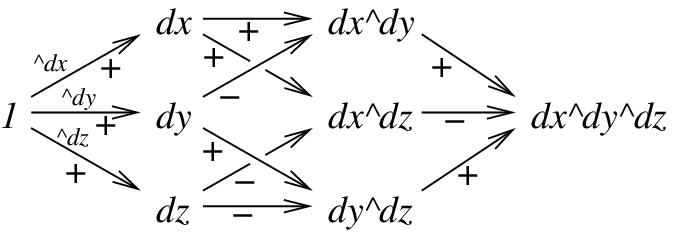}
\end{array}$}
{\bf2.7\qua Signs}\qua
While easy to miss at first glance, the final ingredient in
\figref{fig:Main} is nevertheless significant. Some of the edge
cobordisms (namely, the ones on edges $01\star$, $10\star$, $1\!\star\!0$
and $1\!\star\!1$) also carry little `minus' ($-$) signs. The picture on
the right explains how these signs are determined from a
basis of the exterior algebra in $3$ (or in general, $n$) generators
and from the exterior multiplication operation. Alternatively,
if an edge $\xi$ is labeled by a sequence $(\xi_i)$ in the alphabet
$\{0,1,\star\}$ and if $\xi_j=\star$, then the sign on the edge $\xi$
is $(-1)^\xi:=(-1)^{\sum_{i<j}\xi_i}$.

\parpic[r]{$\begin{array}{c}
  \begin{picture}(0,0)%
\includegraphics{figs/Main2.pstex}%
\end{picture}%
%
%  pstex_opts: -m 0.65 
%
\setlength{\unitlength}{2565sp}%
\begingroup\makeatletter\ifx\SetFigFont\undefined%
\gdef\SetFigFont#1#2#3#4#5{%
  \reset@font\fontsize{#1}{#2pt}%
  \fontfamily{#3}\fontseries{#4}\fontshape{#5}%
  \selectfont}%
\fi\endgroup%
\begin{picture}(5894,4144)(1189,-3968)
\put(1801,-1186){\makebox(0,0)[b]{\smash{\SetFigFont{8}{9.6}{\rmdefault}{\mddefault}{\updefault}{\color[rgb]{0,0,0}$(n_+,n_-)=(2,0)$}%
}}}
\end{picture}

\end{array}$}
{\bf2.8\qua Tangles}\qua
It should be clear to the reader how construct a picture similar to the
one in \figref{fig:Main} for an arbitrary link diagram with
possibly more (or less) crossings. In fact, it should also be clear how
to construct such a picture for any {\em tangle} (a part of a link
diagram bounded inside a disk); the main difference is that now all
cobordisms are bounded within a cylinder, and the part of their
boundary on the sides of the cylinder is a union of vertical straight
lines. An example is on the right.

\section{A frame for our picture} \label{sec:Frame}

Following the previous section, we know how to associate an intricate
but so-far-meaningless picture of a certain $n$--dimensional cube of
smoothings and cobordisms to every link or tangle diagram $T$. We plan
to interpret such a cube as a complex (in the sense of homological
algebra), denoted $\llbracket T\rrbracket$, by thinking of all
smoothings as spaces and of all cobordisms as maps. We plan to set the
$r$'th chain space $\llbracket T\rrbracket^{r-n_-}$ of the complex
$\llbracket T\rrbracket$ to be the ``direct sum'' of the
$\binom{n}{r}$ ``spaces'' (ie, smoothings) at height $r$ in the cube
and to sum the given ``maps'' (ie, cobordisms) to get a
``differential'' for $\llbracket T\rrbracket$.

The problem, of course, is that smoothings aren't spaces and cobordisms
aren't maps. They are, though, objects and morphisms respectively in
a certain category $\Cob^3(\partial T)$ defined below.

\parpic[r]{$\begin{array}{c}\raisebox{-72pt}{
  \includegraphics[width=2in]{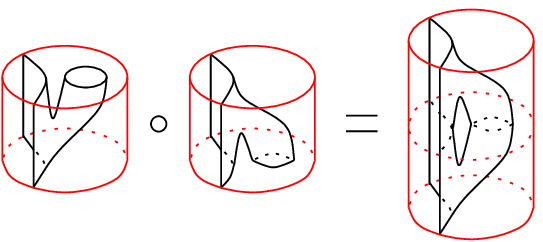}
}\end{array}$}
\begin{definition}
$\Cob^3(\emptyset)$ is the category whose objects are smoothings (ie,
simple curves in the plane) and whose morphisms are cobordisms between such
smoothings as in Section~\ref{sec:Narrative}.6, regarded up to
boundary-preserving isotopies\footnote{
  A slightly different alternative for the choice of morphisms is mentioned
  in Section~\ref{subsec:AbstractCobordism}.
}. Likewise if $B$ is a finite set of points on
the circle (such as the boundary $\partial T$ of a tangle $T$), then
$\Cob^3(B)$ is the category whose objects are smoothings with boundary
$B$ and whose morphisms are cobordisms between such
smoothings as in Section~\ref{sec:Narrative}.8, regarded up to 
boundary-preserving isotopies. In either case the composition of morphisms
is given by placing one cobordism atop the other.  We will use the notation 
$\Cob^3$ as a generic reference either to $\Cob^3(\emptyset)$ or
to $\Cob^3(B)$ for some $B$.
\end{definition}

Next, let us see how in certain parts of homological algebra general
``objects'' and ``morphisms'' can replace spaces and maps; ie, how
arbitrary categories can replace the Abelian categories of vector spaces
and/or $\bbZ$--modules which are more often used in homological algebra.

An {\em pre-additive category} is a category in which the sets of morphisms
(between any two given objects) are Abelian groups and the composition
maps are bilinear in the obvious sense. Let $\calC$ be some arbitrary
category. If $\calC$ is pre-additive, we leave it untouched. If it isn't
pre-additive to start with, we first make it pre-additive by extending every
set of morphisms $\Mor(\calO,\calO')$ to also allow formal
$\bbZ$--linear combinations of ``original'' morphisms
and by extending the composition maps in the
natural bilinear manner. In either case $\calC$ is now pre-additive.

\begin{definition} Given a pre-additive category $\calC$ as above,
the pre-additive category $\Mat(\calC)$ is defined as follows:
\begin{itemize}
\item The objects of $\Mat(\calC)$ are formal direct sums (possibly empty)
  $\oplus_{i=1}^n\calO_i$ of objects $\calO_i$ of $\calC$.
\item If $\calO=\oplus_{i=1}^m\calO_i$ and $\calO'=\oplus_{j=1}^n\calO'_j$,
  then a morphism $F\co \calO'\to\calO$ in $\Mat(\calC)$ will be an $m\times
  n$ matrix $F=(F_{ij})$ of morphisms $F_{ij}\co \calO'_j\to\calO_i$ in
  $\calC$.
\item Morphisms in $\Mat(\calC)$ are added using matrix addition.
\item Compositions of morphisms in $\Mat(\calC)$ are defined by a rule
  modeled on matrix multiplication, but with compositions in $\calC$
  replacing the multiplication of scalars,
  \[
    \left((F_{ij})\circ(G_{jk})\right)_{ik} := \sum_jF_{ij}\circ G_{jk}.
  \]
\end{itemize}
$\Mat(\calC)$ is often called ``the additive closure of $\calC$''.
\end{definition}

\begin{figure}[ht!]
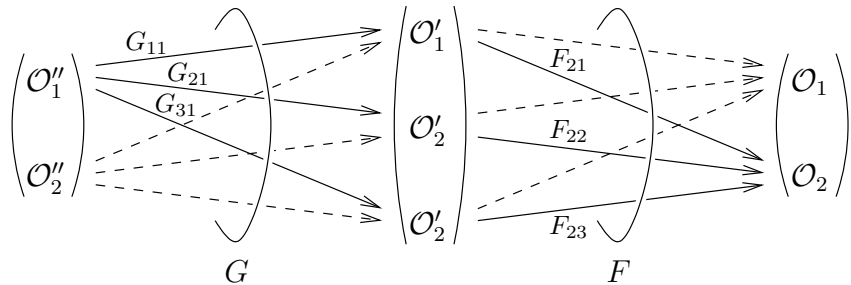
\anchor{fig:Mat}
\begin{center}
  \input figs/Mat.pstex_t
\end{center}
\caption{
  Matrices as bundles of morphisms, the composition $F\circ G$ and the
  matrix element $(F\circ G)_{21}=F_{21}\circ G_{11}+F_{22}\circ
  G_{21}+F_{23}\circ G_{31}$ (in solid lines).
} \label{fig:Mat}
\end{figure}

It is often convenient to represent objects of $\Mat(\calC)$ by column
vectors and morphisms by bundles\footnote{``Bundle'' in the
non-technical sense. Merriam--Webster: bundle: ``a group of things
fastened together for convenient handling''.} of arrows pointing from
one column to another. With this image, the composition $(F\circ
G)_{ik}$ becomes a sum over all routes from $\calO''_k$ to $\calO_i$
formed by connecting arrows.  See \figref{fig:Mat}.

A quick glance at \figref{fig:Main} should convince the reader that
it can be interpreted as a chain of morphisms $\llbracket K\rrbracket
=\left(\xymatrix{\llbracket K\rrbracket^{-3}\ar[r]& \llbracket
K\rrbracket^{-2}\ar[r]& \llbracket K\rrbracket^{-1}\ar[r]&
\llbracket K\rrbracket^0}\right)$ in $\Mat(\Cob^3)$, where
$K=\begin{array}{c}\includegraphics[width=0.7in]{figs/Trefoil.eps}\end{array}$
(if an arrow is missing, such as between the vertices $001$ and $110$
of the cube, simply regard it as $0$).  Likewise, if $T$ is an
$n$--crossing tangle, Section~\ref{sec:Narrative} tells us how it can be
interpreted as a length $n$ chain $\llbracket
T\rrbracket=\left(\xymatrix{\llbracket
T\rrbracket^{-n_-}\ar[r]&\llbracket
T\rrbracket^{-n_-+1}\ar[r]&\dots\ar[r]&\llbracket
T\rrbracket^{n_+}}\right)$. (Strictly speaking, we didn't specify how
to order the equal-height layers of an $n$--dimensional cube as ``column
vector''. Pick an arbitrary such ordering.)  Let us make a room for
such chains by mimicking the standard definition of complexes:

\parpic(78mm,18mm)[r]{\raisebox{-10mm}{$\xymatrix{
  \dots \ar[r] &
  \Omega_a^{r-1} \ar[r]^{d_a^{r-1}} \ar[d]_{F^{r-1}} &
  \Omega_a^r \ar[r]^{d_a^r} \ar[d]_{F^r} &
  \Omega_a^{r+1} \ar[r] \ar[d]_{F^{r+1}} &
  \dots  \\
  \dots \ar[r] &
  \Omega_b^{r-1} \ar[r]^{d_b^{r-1}} &
  \Omega_b^r \ar[r]^{d_b^r} &
  \Omega_b^{r+1} \ar[r] &
  \dots
}$}}
\vskip6mm
\begin{definition} Given a pre-additive category $\calC$, let
$\Kom(\calC)$ be the category of complexes over $\calC$, whose objects
are chains of finite length $\xymatrix{ \dots \ar[r] & \Omega^{r-1}
\ar[r]^{d^{r-1}} & \Omega^r \ar[r]^{d^r} & \Omega^{r+1} \ar[r] & \dots
}$ for which the composition $d^r\circ d^{r-1}$ is $0$ for all $r$, and
whose morphisms $F\co (\Omega_a^r,\,d_a)\to(\Omega_b^r,\,d_b)$ are
commutative diagrams as displayed on the right, in which all arrows are
morphisms in $\calC$. Like in ordinary homological algebra, the
composition $F\circ G$ in $\Kom(\calC)$ is defined via $(F\circ
G)^r:=F^r\circ G^r$.
\end{definition}

\begin{proposition} For any tangle (or knot/link) diagram $T$ the chain
$\llbracket T\rrbracket$ is a complex in $\Kom(\Mat(\Cob^3(\partial
T)))$. That is, $d^r\circ d^{r-1}$ is always $0$ for these chains.
\end{proposition}

\parpic[r]{$\begin{array}{c}
  \includegraphics[width=2in]{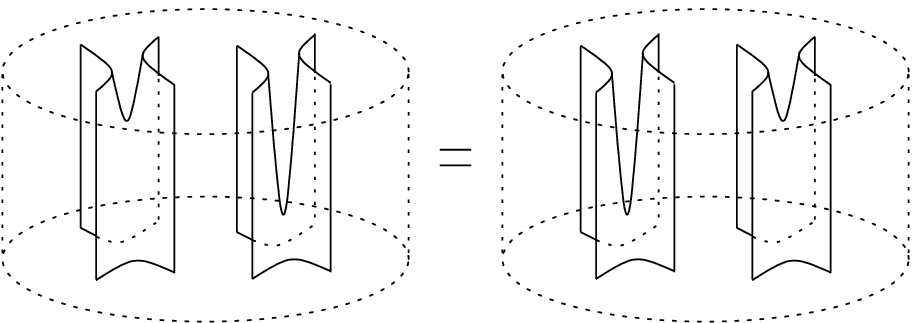}
\end{array}$}
\begin{proof} We have to show that every square face of morphisms in the
cube of $T$ anti-commutes. Every square face of the cube of signs of
Section~\ref{sec:Narrative}.6 carries an odd number of minus signs (this
follows readily from the anti-commutativity of exterior multiplication,
$dx_i\wedge dx_j=-dx_j\wedge dx_i$). Hence, all signs forgotten, we have to
show that every square face in the cube of $T$ positively commutes. This is
simply the fact that spatially separated saddles can be time-reordered
within a cobordism by an isotopy.
\end{proof}

\section{Invariance} \label{sec:Invariance}

\subsection{Preliminaries} \label{subsec:InvPrels}

The ``formal'' complex $\llbracket T\rrbracket$ is not a tangle
invariant in any sense. We will claim and prove, however, that $\llbracket
T\rrbracket$, regarded within $\Kom(\Mat(\Cobl^3))$ for some quotient
$\Cobl^3$ of $\Cob^3$, is invariant up to homotopy. But first we have
to define these terms.

\subsubsection{Homotopy in formal complexes} \label{subsubsec:Homotopy}
\parpic[r]{$\xymatrix@C=2cm{
  \Omega_a^{r-1} \ar[r]^{d_a^{r-1}}
    \ar@<-2pt>[d]_{F^{r-1}} \ar@<2pt>[d]^{G^{r-1}} &
  \Omega_a^r \ar[r]^{d_a^r} \ar[ld]_{h^r}
    \ar@<-2pt>[d]_{F^r} \ar@<2pt>[d]^{G^r} &
  \Omega_a^{r+1} \ar[ld]_{h^{r+1}}
    \ar@<-2pt>[d]_{F^{r+1}} \ar@<2pt>[d]^{G^{r+1}} \\
  \Omega_b^{r-1} \ar[r]^{d_b^{r-1}} &
  \Omega_b^r \ar[r]^{d_b^r} &
  \Omega_b^{r+1}
}$}
\vskip3mm
Let $\calC$ be a category.  Just like in ordinary homological algebra, we
say that two morphisms $F,G\co (\Omega_a^r)\to(\Omega_b^r)$ in $\Kom(\calC)$
are {\em homotopic} (and we write $F\sim G$) if there exists ``backwards
diagonal'' morphisms $h^r\co \Omega_a^r\to \Omega_b^{r-1}$ so that
$F^r-G^r=h^{r+1}d^r+d^{r-1}h^r$ for all $r$.

Many of the usual properties of homotopies remain true in the formal
case, with essentially the same proofs. In particular, homotopy is an
equivalence relation and it is invariant under composition both on the
left and on the right; if $F$ and $G$ are homotopic and $H$ is some
third morphism in $\Kom(\calC)$, then $F\circ H$ is homotopic to
$G\circ H$ and $H\circ F$ is homotopic to $H\circ G$, whenever these
compositions make sense. Thus we can make the following definition:

\begin{definition} $\Komh(\calC)$ is $\Kom(\calC)$ modulo homotopies. That
is, $\Komh(\calC)$ has the same objects as $\Kom(\calC)$ (formal
complexes), but homotopic morphisms in $\Kom(\calC)$ are declared to be the
same in $\Komh(\calC)$ (the $/h$ stands for ``modulo homotopy'').
\end{definition}

As usual, we say that two complexes $(\Omega_a^r)$ and $(\Omega_b^r)$
in $\Kom(\calC)$ are {\em homotopy equivalent} (and we write
$(\Omega_a^r)\sim(\Omega_b^r)$) if they are isomorphic in
$\Komh(\calC)$. That is, if there are morphisms
$F\co (\Omega_a^r)\to(\Omega_b^r)$ and $G\co (\Omega_b^r)\to(\Omega_a^r)$ so
that the compositions $G\circ F$ and $F\circ G$ are homotopic to the
identity automorphisms of $(\Omega_a^r)$ and $(\Omega_b^r)$,
respectively. It is routine to verify that homotopy equivalence is an
equivalence relation on complexes.

\subsubsection{The quotient $\protect{\Cobl^3}$ of $\protect{\Cob^3}$}
\label{subsubsec:STFourTu}
We mod out the morphisms of the category $\Cob^3$ by the relations
$S$, $T$ and $\FourTu$ defined below and call the
resulting quotient $\Cobl^3$ (the $/l$ stands for ``modulo local
relations'').

\parpic[r]{$\begin{array}{c}
  \includegraphics[height=1cm]{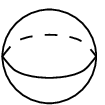}
\end{array}\hspace{-2mm}=0$}
The $S$ relation says that whenever a cobordism contains a connected
component which is a closed sphere (with no boundary), it is set equal
to zero (remember that we make all categories pre-additive, so $0$ always
makes sense).

\parpic[r]{$\begin{array}{c}
  \includegraphics[height=1cm]{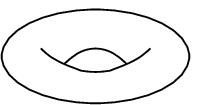}
\end{array}$\hspace{-1mm}=2}
The $T$ relation says that whenever a cobordism contains a connected
component which is a closed torus (with no boundary), that component
may be dropped and replaced by a numerical factor of $2$ (remember that
we make all categories pre-additive, so multiplying a cobordism by a
numerical factor makes sense).

\parpic[r]{$\begin{array}{c}
  \includegraphics[height=1.5cm]{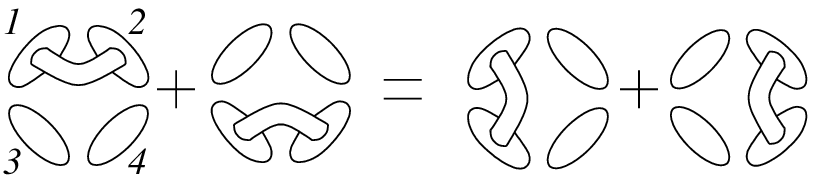}
\end{array}$}
To understand $\FourTu$, start from some given cobordism $C$ and assume
its intersection with a certain ball is the union of four disks $D_1$
through $D_4$ (these disks may well be on different connected
components of $C$). Let $C_{ij}$ denote the result of removing $D_i$
and $D_j$ from $C$ and replacing them by a tube that has the same
boundary. The ``four tube'' relation $\FourTu$ asserts that
$C_{12}+C_{34}=C_{13}+C_{24}$.

The local nature of the $S$, $T$ and $\FourTu$ relations implies that the
composition operations remain well defined in $\Cobl^3$ and hence it is also
a pre-additive category.

\subsection{Statement} \label{subsec:Statement}

Throughout this paper we will often short $\Kob(\emptyset)$, $\Kob(B)$
and $\Kob$ for $\Kom(\Mat(\Cobl^3(\emptyset)))$,
$\Kom(\Mat(\Cobl^3(B)))$ and $\Kom(\Mat(\Cobl^3))$. Likewise we will
often short $\Kobh(\emptyset)$, $\Kobh(B)$ and $\Kobh$ for
$\Komh(\Mat(\Cobl^3(\emptyset)))$, $\Komh(\Mat(\Cobl^3(B)))$ and
$\Komh(\Mat(\Cobl^3))$ respectively.

\begin{figure}[ht!]\anchor{fig:ReidMoves}
\[ \includegraphics[width=5in]{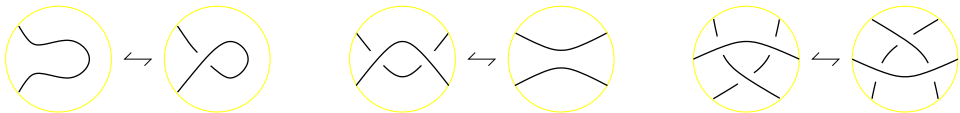} \]
\caption{
  The three Reidemeister moves $R1$, $R2$ and $R3$.
} \label{fig:ReidMoves}
\end{figure}

\begin{theorem}[The Invariance Theorem] \label{thm:invariance}
The isomorphism class of the complex $\llbracket T\rrbracket$
regarded in $\Kobh$ is an invariant of the tangle $T$. That
is, it does not depend on the ordering of the layers of a cube as column
vectors and on the ordering of the crossings and it is invariant under
the three Reidemeister moves (reproduced in \figref{fig:ReidMoves}).
\end{theorem}

\subsection{Proof} \label{subsec:Proof}

The independence on the ordering of the layers of a cube as column
vectors is left as an exercise to the reader with no hints supplied.
The independence on the ordering of the crossings is also left as an
exercise to the reader (hint: when reordering, take signs from the signs
that appear in the usual action of a symmetric group on the basis of an
exterior algebra). Invariance under the Reidemeister moves is shown,
in this section,
{\em just for the ``local'' tangles representing these moves}.
Invariance under the Reidemeister moves applied within
larger tangles or knots or links follows from the nearly tautological
good behaviour of $\llbracket T\rrbracket$ with respect to tangle
compositions discussed in Section~\ref{sec:PlanarAlgebra} (see also
Exercise~\ref{ex:loc2glob}).

\begin{figure}[ht!]\anchor{fig:R1Invariance}
\[ \xymatrix@C=25mm@R=26mm{
  \eps{8mm}{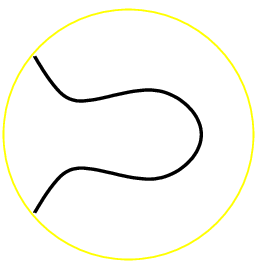}
    \ar[r]^0 \ar@<2pt>[d]^{F^0=\eps{8mm}{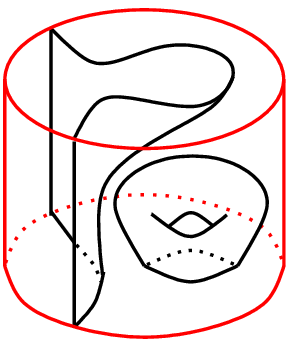}-\eps{8mm}{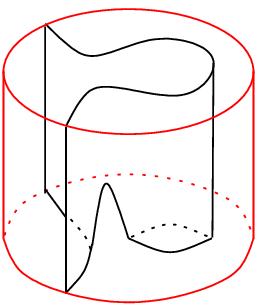}} &
  0 \ar@<2pt>[d]^0 \\
  \eps{8mm}{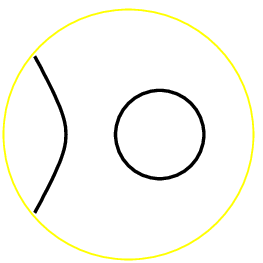}
    \ar@<2pt>[u]^{G^0=\eps{8mm}{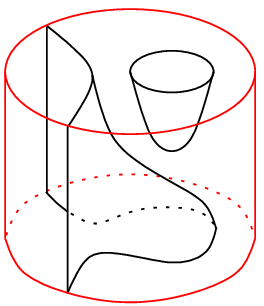}} \ar@<2pt>[r]^{d=\eps{8mm}{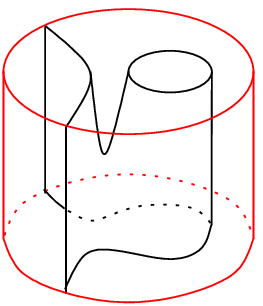}} &
  \eps{8mm}{R1-1}
    \ar@<2pt>[u]^0 \ar@<2pt>[l]^{h=\eps{8mm}{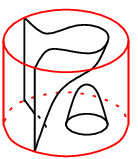}}
} \]
\caption{Invariance under $R1$} \label{fig:R1Invariance}
\end{figure}

{\bf Invariance under the Reidemeister move $R1$}\qua (See
\figref{fig:R1Invariance}) We have to
show that the formal complex
$
  \llbracket\eps{5mm}{R1-1}\rrbracket
  = \left(\xymatrix{0\ar[r]&\underline{\eps{5mm}{R1-1}}\ar[r]&0}\right)
$
is homotopy equivalent to the formal complex
$
  \llbracket\eps{5mm}{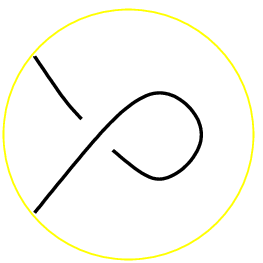}\rrbracket
  = \left(\!\!\xymatrix{
    0\ar[r] &
    \underline{\eps{5mm}{R1-3}}\ar[r]^d &
    \eps{5mm}{R1-1}\ar[r] &
    0
  }\!\!\right)
$,
in which $d=\eps{7mm}{R1-7}$ (in both complexes we have underlined the
$0$th term). To do this we construct (homotopically inverse) morphisms
$F\co \llbracket\eps{5mm}{R1-1}\rrbracket \to
\llbracket\eps{5mm}{R1-2}\rrbracket$ and
$G\co \llbracket\eps{5mm}{R1-2}\rrbracket \to
\llbracket\eps{5mm}{R1-1}\rrbracket$. The morphism $F$ is defined by
$F^0=\eps{7mm}{R1-4}-\eps{7mm}{R1-5}$ (in words: a vertical curtain
union a torus with a downward-facing disk removed, minus a simple
saddle) and $F^{\neq 0}=0$. The morphism $G$ is defined by
$G^0=\eps{7mm}{R1-6}$ (a vertical curtain union a cup) and $G^{\neq
0}=0$. The only non-trivial commutativity to verify is $dF^0=0$,
which follows from
$\eps{7mm}{R1-7}\circ\eps{7mm}{R1-4}=\eps{7mm}{R1-7}\circ\eps{7mm}{R1-5}$,
and where the latter identity holds because both of its sides are the
same --- vertical curtains with an extra handle attached. If follows
from the $T$ relation that $GF=I$.

\parpic[r]{$\eps{15mm}{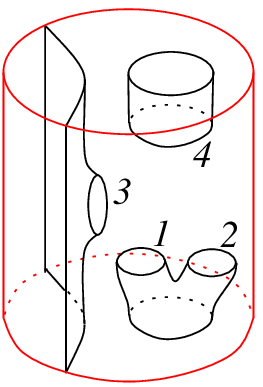}$}
Finally, consider the (homotopy) morphism
$h=\eps{7mm}{R1-8}\co \llbracket\eps{5mm}{R1-2}\rrbracket^1=\eps{5mm}{R1-1}
\to \eps{5mm}{R1-2} = \llbracket\eps{5mm}{R1-2}\rrbracket^0$. Clearly,
$F^1G^1-I+dh=-I+dh=0$. We claim that it follows from the \FourTu{}
relation that $F^0G^0-I+hd=0$ and hence $FG\sim I$ and we have proven
that $\llbracket\eps{5mm}{R1-1}\rrbracket \sim
\llbracket\eps{5mm}{R1-2}\rrbracket$. Indeed, let $C$ be the cobordism
(with four punctures labeled $1$--$4$) shown on the right, and consider
the cobordisms $C_{ij}$ constructed from it as in
Section~\ref{subsubsec:STFourTu}. Then $C_{12}$ and $C_{13}$ are the
first and second summands in $F^0 G^0$, $C_{24}$ is the identity
morphisms $I$ and $C_{34}$ is $hd$. Hence the $\FourTu$ relation
$C_{12}-C_{13}-C_{24}+C_{34}=0$ is precisely our assertion,
$F^0G^0-I+hd=0$. \qed

\begin{figure}[ht!]\anchor{fig:Reid2Proof}
\[ \eps{5.1in}{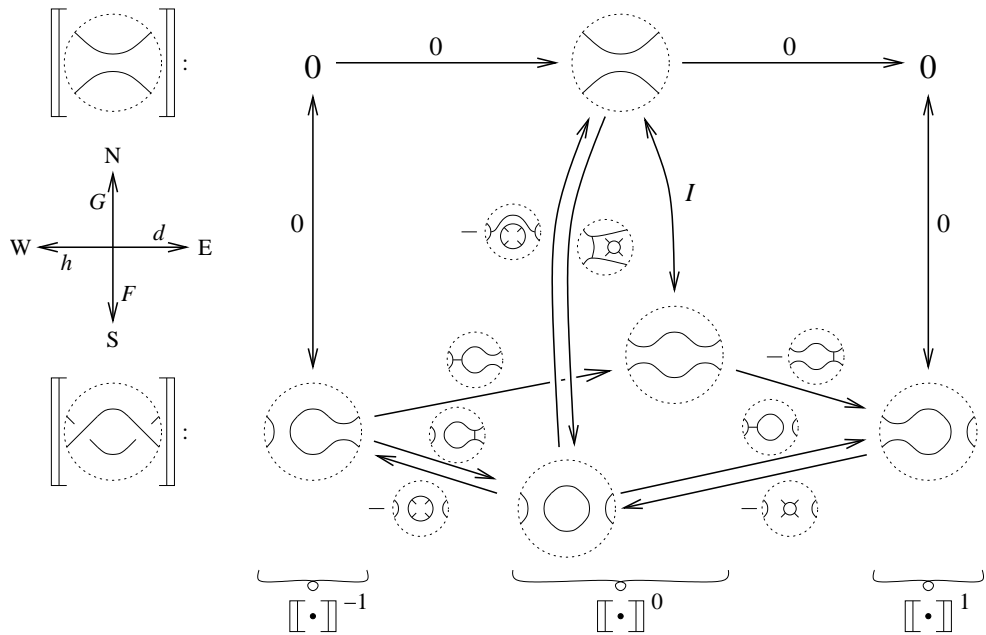} \]
\caption{Invariance under the Reidemeister move $R2$.}
\label{fig:Reid2Proof}
\end{figure}

\vskip 2mm
{\bf Invariance under the Reidemeister move $R2$}\qua This
invariance proof is very similar in spirit to the previous one, and
hence we will allow ourselves to be brief. The proof appears in whole
in \figref{fig:Reid2Proof}; let us just add some explanatory
words.  In that figure the top row is the formal complex
$\llbracket\eps{4mm}{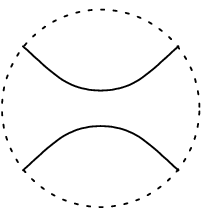}\rrbracket$ and the bottom row is the
formal complex $\llbracket\eps{4mm}{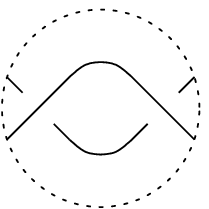}\rrbracket$. Also, all
eastward arrows are (components of) differentials, the southward arrows
are the components of a morphism
$F\co \llbracket\eps{4mm}{R2-1}\rrbracket \to
\llbracket\eps{4mm}{R2-2}\rrbracket$, the northward arrows are the
components of a morphism $G\co \llbracket\eps{4mm}{R2-2}\rrbracket \to
\llbracket\eps{4mm}{R2-1}\rrbracket$, and the westward arrows are the
non-zero components of a homotopy $h$ proving that $FG\sim I$.
Finally, in \figref{fig:Reid2Proof} the symbol $\HSaddleSymbol$ (or
its variants, $\ISaddleSymbol$, etc.) stands for the saddle cobordism
$\HSaddleSymbol\co \smoothing\to\hsmoothing$
(or $\ISaddleSymbol\co \hsmoothing\to\smoothing$, etc.) as in
\figref{fig:Main}, and the symbols $\fourwheel$ and $\fourinwheel$
stand for the cap morphism $\fourwheel\co \emptyset\to\bigcirc$ and the
cup morphism $\fourinwheel\co \bigcirc\to\emptyset$.

We leave it for the reader to verify the following facts which together
constitute a proof of invariance under the Reidemeister move $R2$:

\begin{itemize}
\item $dF=0$ (only uses isotopies).
\item $Gd=0$ (only uses isotopies).
\item $GF=I$ (uses the relation $S$).
\item The hardest --- $FG-I=hd+dh$ (uses the \FourTu{}
  relation). \qed
\end{itemize}

\begin{remark} \label{rem:R2isSDR}
The morphism $G\co \llbracket\eps{4mm}{R2-2}\rrbracket \to
\llbracket\eps{4mm}{R2-1}\rrbracket$ at the heart of the above proof
is a little more than a homotopy equivalence. A routine check shows
that it is in fact a strong deformation retract in the sense of the following
definition.
\end{remark}

\parpic[r]{\parbox{1.6in}{\vbox to -4mm{\raisebox{-64pt}{
  $\eps{1.5in}{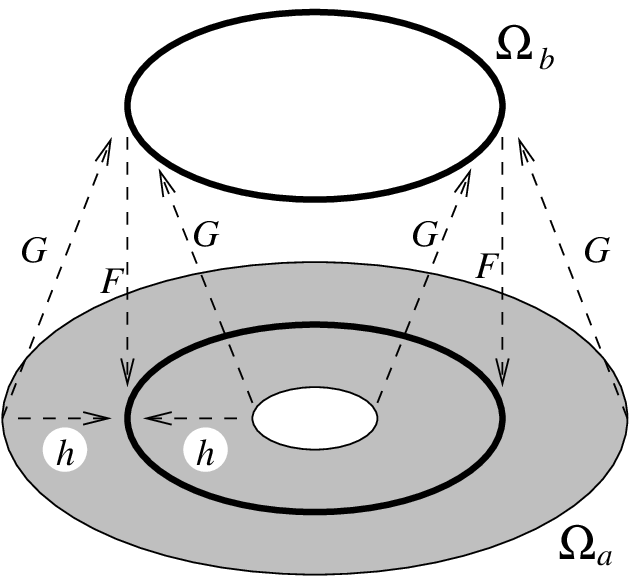}$
}}}}
\begin{definition} A morphism of complexes $G\co \Omega_a\to\Omega_b$ is
said to be a {\em strong deformation retract} if there is a morphism
$F\co \Omega_b\to\Omega_a$ and homotopy maps $h$ from $\Omega_a$ to itself
so that $GF=I$, $I-FG=dh+hd$ and $hF=0$. In this case we say that
$F$ is the {\em inclusion in a strong deformation retract}. Note that a
strong deformation retract is in particular a homotopy equivalence.
The geometric origin of this notion is the standard notion of a strong
deformation retract in homotopy theory as sketched on the right.
\end{definition}

{\bf Invariance under the Reidemeister move $R3$}\qua This is
the easiest and hardest move. Easiest because it doesn't require any
further use of the $S$, $T$ and \FourTu{} relations --- it just follows from
the $R2$ move and some `soft' algebra (just as in the case of the Kauffman
bracket, whose invariance under $R3$ follows `for free' from its invariance
under $R2$; see eg,~\cite[Lemma~2.4]{Kauffman:Bracket}). Hardest because
it involves the most crossings and hence the most complicated complexes. We
will attempt to bypass that complexity by appealing to some standard
constructions and results from homological algebra.

\parpic[r]{$\xymatrix@C=15mm{
  \Omega_0^r \ar[r]^{-d_0^r} \ar[d]^{\Psi^r} &
    \Omega_0^{r+1} \ar[r]^{-d_0^{r+1}} \ar[d]^{\Psi^{r+1}}
      \ar@{.}[dl]|{\displaystyle \oplus} &
    \Omega_0^{r+2} \ar[d]^{\Psi^{r+2}} \ar@{.}[dl]|{\displaystyle \oplus} \\
  \Omega_1^r \ar[r]_{d_1^r} &
    \Omega_1^{r+1} \ar[r]_{d_1^{r+1}} &
    \Omega_1^{r+2}
}$}
Let $\Psi\co (\Omega_0^r,\,d_0)\to(\Omega_1^r,\,d_1)$ be a morphism of
complexes. The {\em cone} $\Gamma(\Psi)$ of $\Psi$ is the complex with
chain spaces $\Gamma^r(\Psi)=\Omega_0^{r+1}\oplus\Omega_1^r$ and with
differentials $\tilde{d}^r=\left(\begin{matrix} -d_0^{r+1} & 0
\\ \Psi^{r+1} & d_1^r \end{matrix}\right)$. The following two lemmas
explain the relevance of cones to the task at hand and are
easy\footnote{Hard, if one is punctual about signs\ldots.} to verify:

\begin{lemma} \label{lem:cones}
$\llbracket\overcrossing\rrbracket =
\Gamma(\llbracket\HSaddleSymbol\rrbracket)[-1]$ and
$\llbracket\undercrossing\rrbracket =
\Gamma(\llbracket\ISaddleSymbol\rrbracket)$, where
$\llbracket\HSaddleSymbol\rrbracket$ and
$\llbracket\ISaddleSymbol\rrbracket$ are the saddle morphisms
$\llbracket\HSaddleSymbol\rrbracket\co  \llbracket\smoothing\rrbracket
\to \llbracket\hsmoothing\rrbracket$ and
$\llbracket\ISaddleSymbol\rrbracket\co  \llbracket\hsmoothing\rrbracket
\to \llbracket\smoothing\rrbracket$ and where $\cdot[s]$ is the
operator that shifts complexes $s$ units to the left:
$\Omega[s]^r:=\Omega^{r+s}$. \qed
\end{lemma}

\begin{figure}[ht!]\anchor{fig:Reid3Proof}
\[
  \Psi_L:\ \eps{50mm}{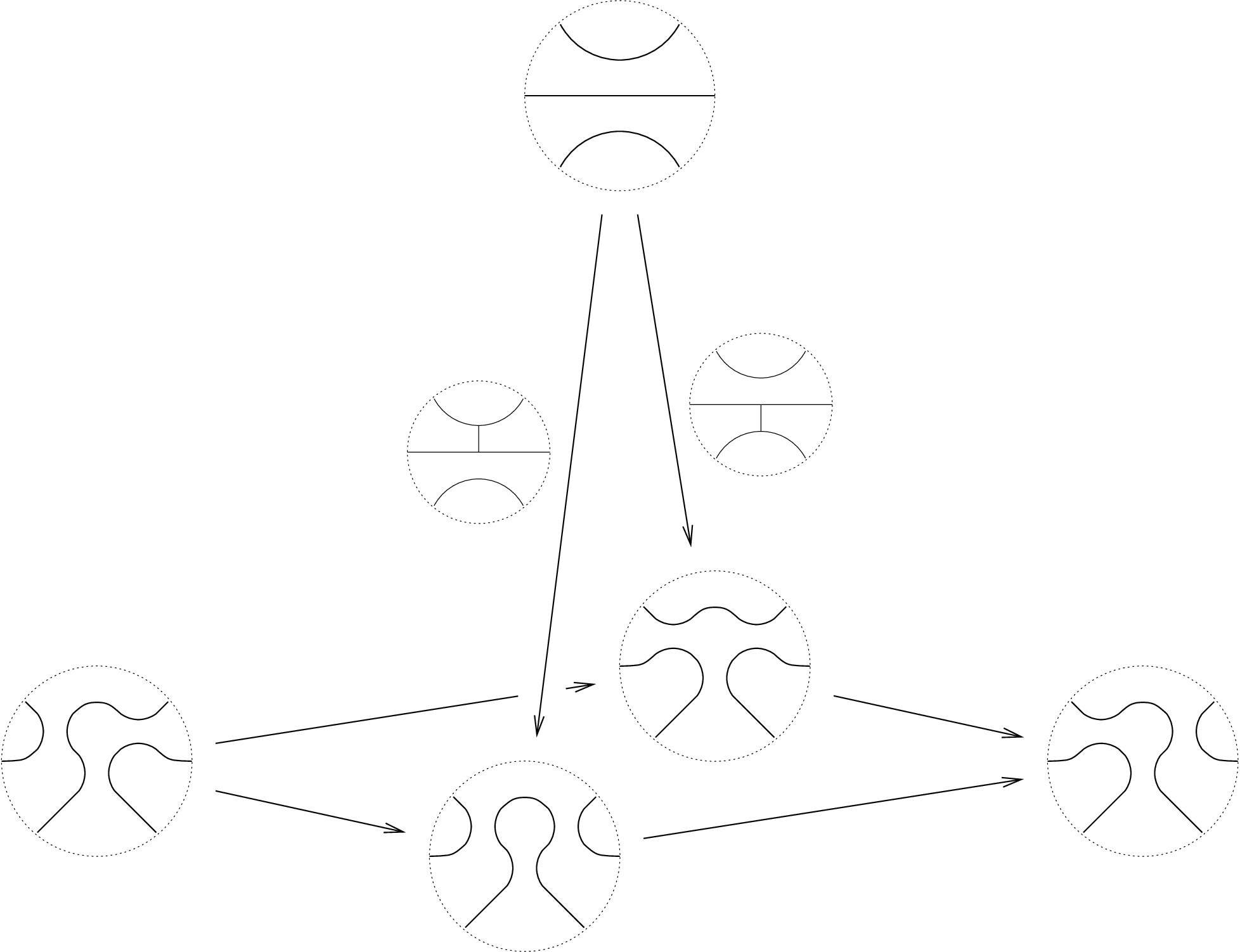}
  \qquad
  \Psi_R:\ \eps{50mm}{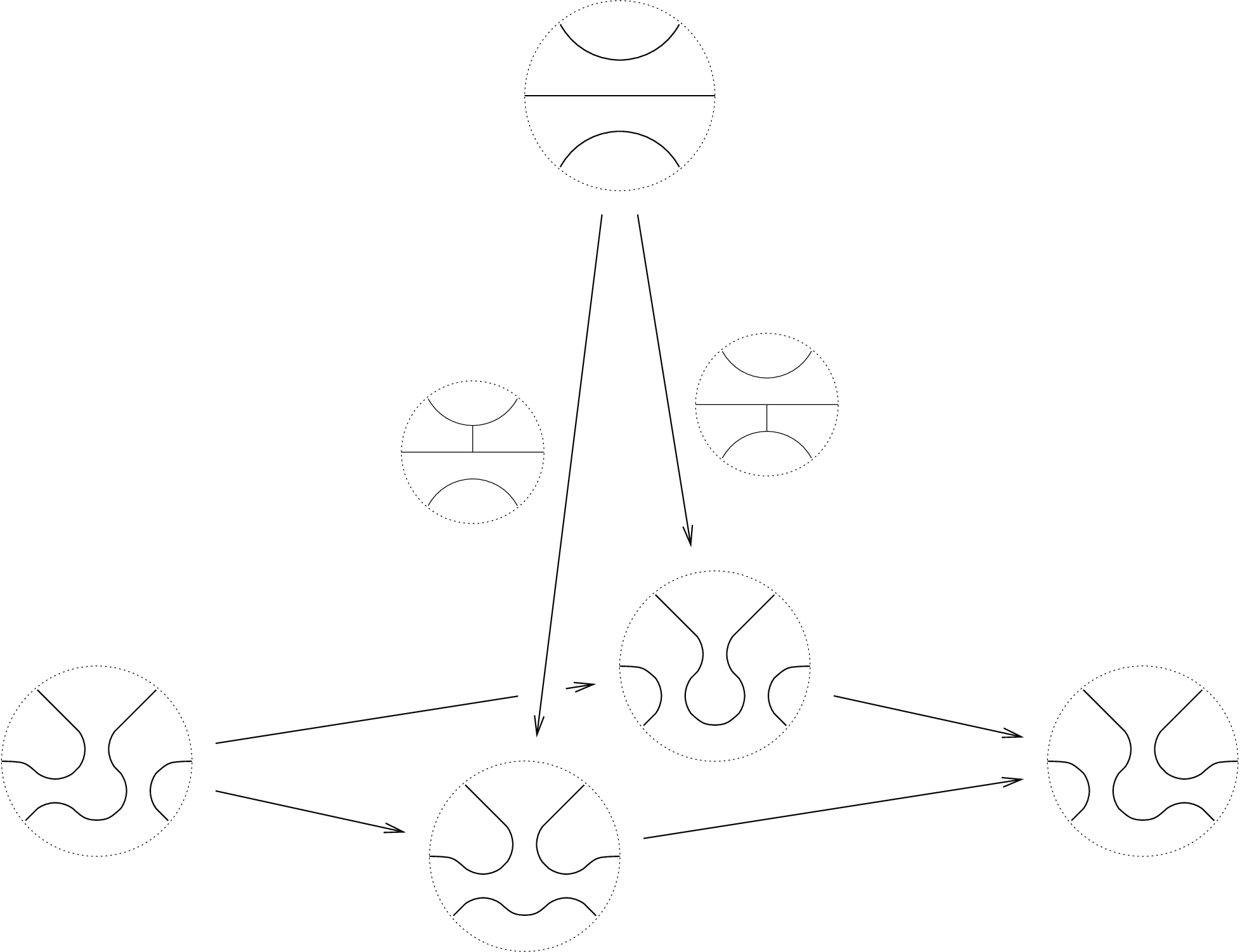}
\]
\caption{The two sides of the Reidemeister move $R3$.}
\label{fig:Reid3Proof}
\end{figure}
\parpic[r]{\parbox{1in}{\vbox to -4mm{\raisebox{-30pt}{
  $\xymatrix{
    \Omega_{0a} \ar[d]^\Psi \ar@<2pt>[r]^{G_0} &
      \Omega_{0b} \ar@<2pt>[l]^{F_0} \\
    \Omega_{1a} \ar@<2pt>[r]^{F_1} &
      \Omega_{1b} \ar@<2pt>[l]^{G_1}
  }$
}}}}
\begin{lemma} \label{lem:ConeHomotopy}The cone construction is invariant 
up to homotopy under compositions with the inclusions in strong
deformation retracts. That is, consider the diagram of complexes and
morphisms that appears on the right. If in that diagram $G_0$ is a
strong deformation retract with inclusion $F_0$, then\break
\end{lemma}

{\sl
the cones
$\Gamma(\Psi)$ and $\Gamma(\Psi F_0)$ are homotopy equivalent, and if
$G_1$ is a strong deformation retract with inclusion $F_1$, then
the cones $\Gamma(\Psi)$ and $\Gamma(F_1\Psi)$ are homotopy
equivalent. (This lemma remains true if $F_{0,1}$ are strong
deformation retracts and $G_{0,1}$ are the corresponding inclusions,
but we don't need that here).}

\parpic[r]{$\eps{45mm}{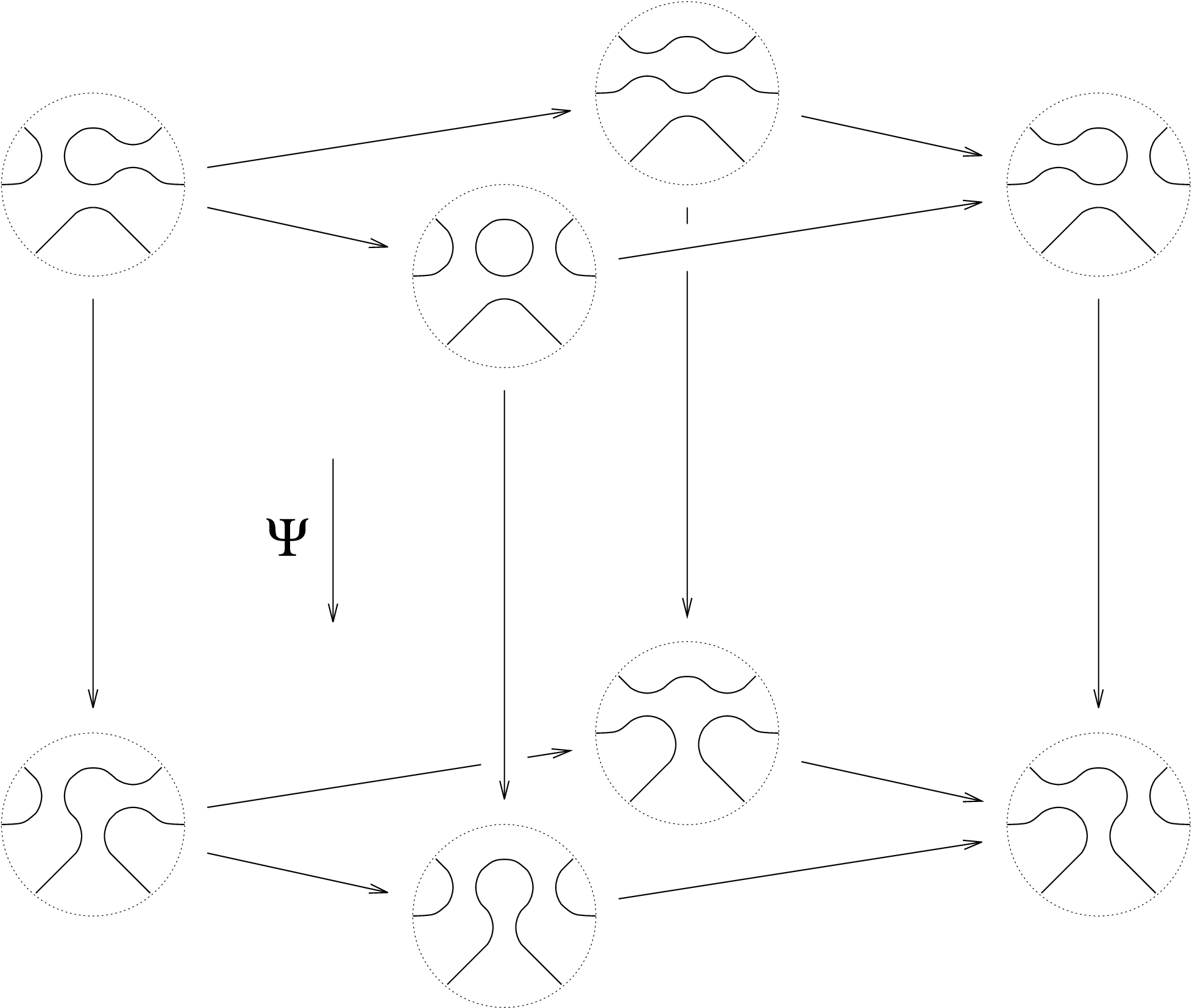}$}
We note that Lemma~\ref{lem:cones} can also be interpreted
(and remains true) in a ``skein theoretic'' sense, where each
of $\overcrossing$ and $\HSaddleSymbol$ (or $\undercrossing$ and
$\ISaddleSymbol$) represents just a small disk neighborhood inside an
otherwise-equal bigger tangle.  Thus, applying Lemma~\ref{lem:cones}
to the bottom crossing in the tangle $\eps{5mm}{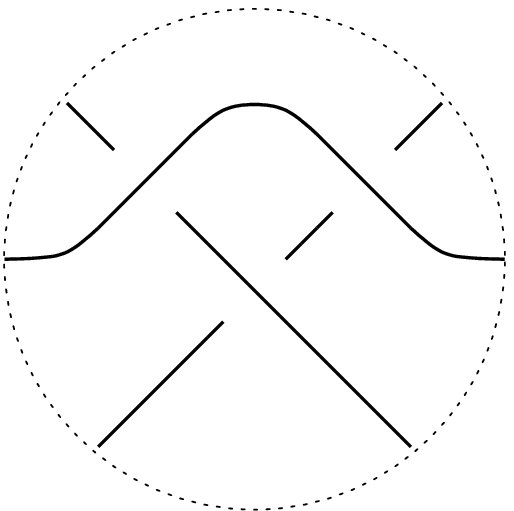}$ we find that
the complex $\llbracket\eps{5mm}{R3-1}\rrbracket$ is the cone
of the morphism $\Psi=\llbracket\eps{5mm}{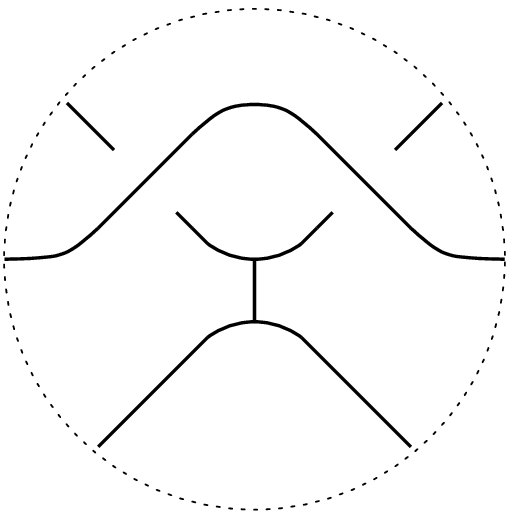}\rrbracket:
\llbracket\eps{5mm}{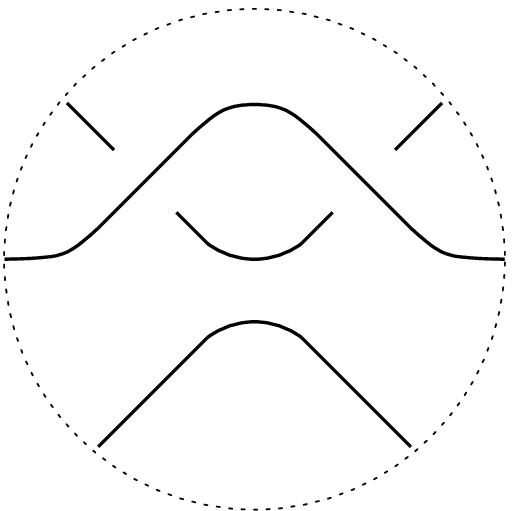}\rrbracket \to
\llbracket\eps{5mm}{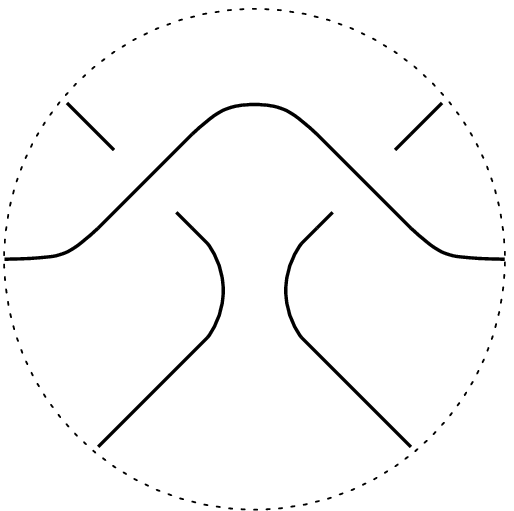}\rrbracket$, which in itself is the bundle
of four morphisms corresponding to the four smoothings of the two
remaining crossings of $\eps{5mm}{R3-2}$ (see the diagram on the
right). We can now use Lemma~\ref{lem:ConeHomotopy} and the inclusion
$F$ of \figref{fig:Reid2Proof} (notice Remark~\ref{rem:R2isSDR})
to replace the top layer of this cube by a single object,
$\eps{5mm}{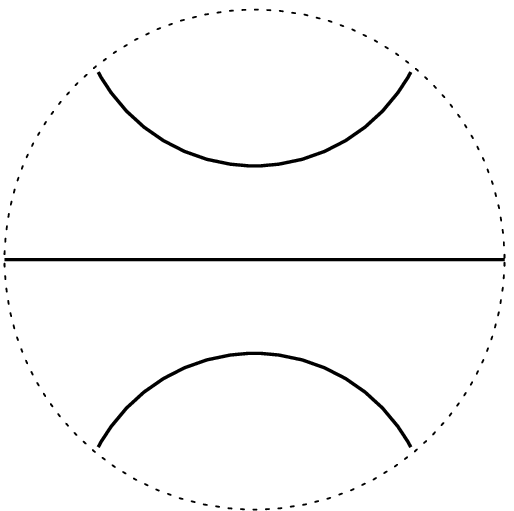}$. Thus $\llbracket\eps{5mm}{R3-1}\rrbracket$ is homotopy
equivalent to the cone of the vertical morphism $\Psi_L=\Psi F$ of
\figref{fig:Reid3Proof}. A similar treatment applied to the complex
$\llbracket\eps{5mm}{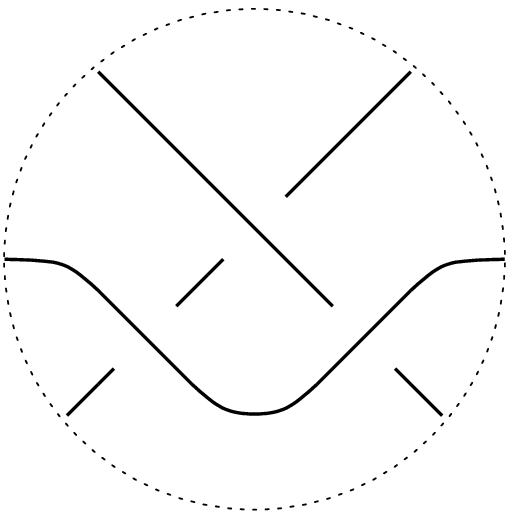}\rrbracket$ yields the cone of the morphism
$\Psi_R$ of \figref{fig:Reid3Proof}. But up to isotopies $\Psi_L$
and $\Psi_R$ are the same.

Given that we distinguish left from right, there
is another variant of the third Reidemeister move to check ---
$\eps{15mm}{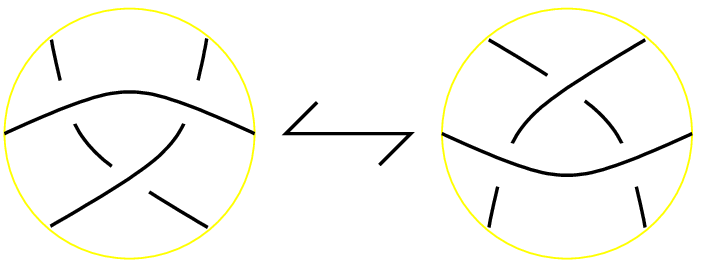}$. We leave it to the reader to verify that invariance
here can be proven in a similar way, except using the second half of
Lemma~\ref{lem:ConeHomotopy}. \qed

\proof[Proof of Lemma~\ref{lem:ConeHomotopy}]
Let $h_0\co \Omega_{0a}^\star\to\Omega_{0a}^{\star-1}$ be a homotopy for which
$I-F_0G_0=dh_0+h_0d$ and $h_0F_0=0$.
Then the diagram in \figref{diag:hom} defines morphisms
$\xymatrix{\Gamma(\Psi F_0) \ar@<2pt>[r]^(0.55){\tilde{F}_0} &
\Gamma(\Psi) \ar@<2pt>[l]^(0.45){\tilde{G}_0}}$ and a homotopy
$\tilde{h}_0\co \Gamma(\Psi)^\star\to\Gamma(\Psi)^{\star-1}$.
\begin{figure}[ht!]\anchor{diag:hom}
\[ \xymatrix{
  \Gamma(\Psi F_0): & &
  \left({
    \begin{matrix} \Omega_{0b}^{r+1} \\ \Omega_{1a}^r \end{matrix}
  }\right)
    \ar[rrrrr]^{
      \tilde{d} = \left({
        \begin{matrix} -d & 0 \\ \Psi F_0 & d \end{matrix}
      }\right)
    }
    \ar@<2pt>[dd]^{
      \tilde{F}_0^r := \left({
        \begin{matrix} F_0 & 0 \\ 0 & I \end{matrix}
      }\right)
    } &&&&&
  \left({
    \begin{matrix} \Omega_{0b}^{r+2} \\ \Omega_{1a}^{r+1} \end{matrix}
  }\right)
    \ar@<2pt>[dd]^{\tilde{F}_0^{r+1}} \\ \\
  \Gamma(\Psi): & &
  \left({
    \begin{matrix} \Omega_{0a}^{r+1} \\ \Omega_{1a}^r \end{matrix}
  }\right)
    \ar@<2pt>[rrrrr]^{
      \tilde{d} =  \left({
        \begin{matrix} -d & 0 \\ \Psi & d \end{matrix}
      }\right)
    }
    \ar@<2pt>[uu]^{
      \tilde{G}_0^r := \left({
        \begin{matrix} G_0 & 0 \\ \Psi h_0 & I \end{matrix}
      }\right)
    } &&&&&
  \left({
    \begin{matrix} \Omega_{0a}^{r+2} \\ \Omega_{1a}^{r+1} \end{matrix}
  }\right)
    \ar@<2pt>[lllll]^{
      \tilde{h}_0 := \left({
        \begin{matrix} -h_0 & 0 \\ 0 & 0 \end{matrix}
      }\right)
    }
    \ar@<2pt>[uu]^{\tilde{G}_0^{r+1}}
} \]
\nocolon\caption{}\label{diag:hom}
\end{figure}
We leave it to the reader to verify that $\tilde{F}_0$ and $\tilde{G}_0$
are indeed morphisms of complexes and that $\tilde{G}_0\tilde{F}_0=I$ and
$I-\tilde{F}_0\tilde{G}_0=\tilde{d}\tilde{h}_0+\tilde{h}_0\tilde{d}$ and
hence $\Gamma(\Psi F_0)$ and $\Gamma(\Psi)$ are homotopy equivalent. A
similar argument shows that $\Gamma(\Psi)$ and $\Gamma(F_1\Psi)$
are also homotopy equivalent. \qed

\begin{figure}[p]\anchor{fig:R3Full}
\begin{center}\begin{sideways}
  $\eps{6.4in}{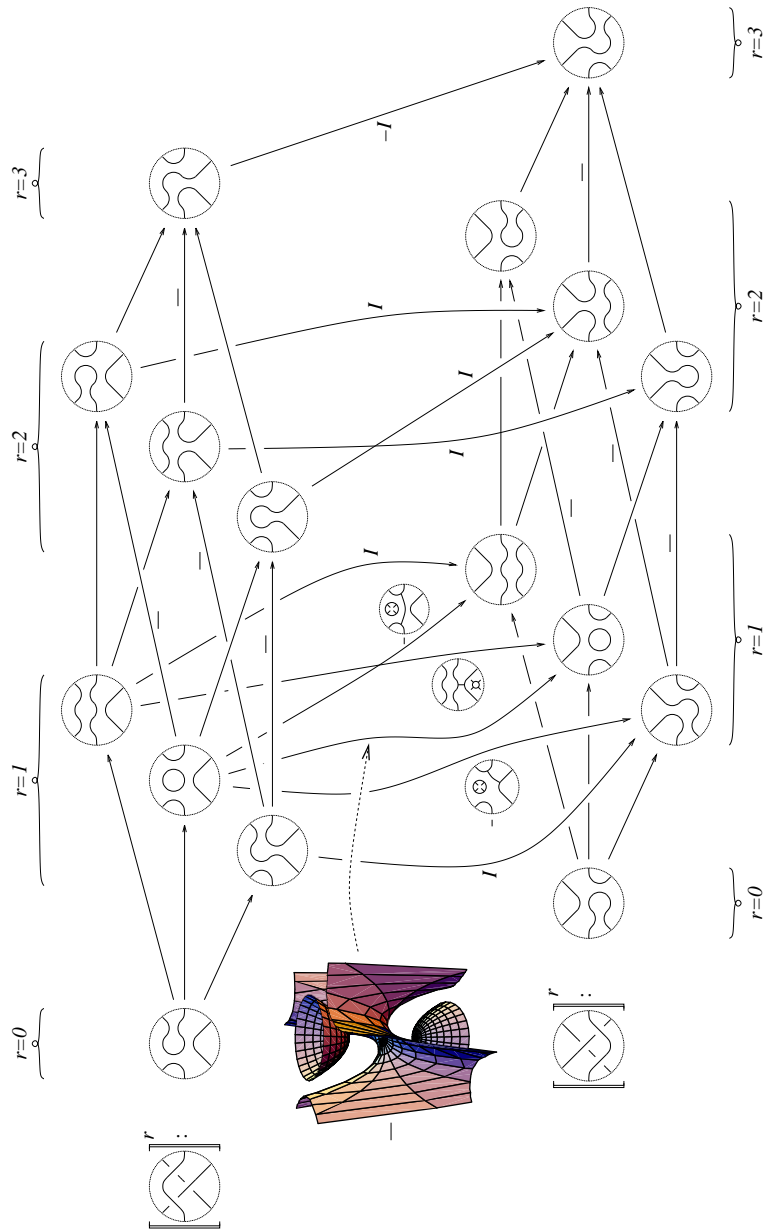}$
\end{sideways}\end{center}
\caption{
  Invariance under $R3$ in more detail than is strictly necessary.
  Notice the minus signs and consider all missing arrows between the
  top layer and the bottom layer as $0$.
} \label{fig:R3Full}
\end{figure}

\newpage
\begin{remark} The proof of invariance under the Reidemeister move $R3$
was presented in a slightly roundabout way, using cones and their
behaviour under retracts. But there is no difficulty in unraveling
everything to get concrete (homotopically invertible) morphisms between
the formal complexes at the two sides of $R3$. This is done (for just one
of the two variants of $R3$) in \figref{fig:R3Full}. The most
interesting cobordism in \figref{fig:R3Full} is displayed --- a
cubic saddle\footnote{A ``monkey saddle'', comfortably seating a monkey
with two legs and a tail.} $(z,3\operatorname{Re}(z^3))$ bound in the
cylinder $[z\leq 1]\times[-1,1]$ plus a cup and a cap.
\end{remark}

\begin{exercise} Verify that \figref{fig:R3Full} indeed defines a map
between complexes --- that the morphism defined by the downward arrows
commutes with the (right pointing) differentials.
\end{exercise}

\begin{exercise} \label{ex:loc2glob} Even though this will be done in
a formal manner in the next section, we recommend that the reader will
pause here to convince herself that the ``local'' proofs above generalize
to Reidemeister moves performed within larger tangles, and hence that
Theorem~\ref{thm:invariance} is verified. The mental picture you will thus
create in your mind will likely be a higher form of understanding than
its somewhat arbitrary serialization into a formal stream of words below.
\end{exercise}

\newpage

\section{Planar algebras and tangle compositions} \label{sec:PlanarAlgebra}

\parpic[r]{$\eps{45mm}{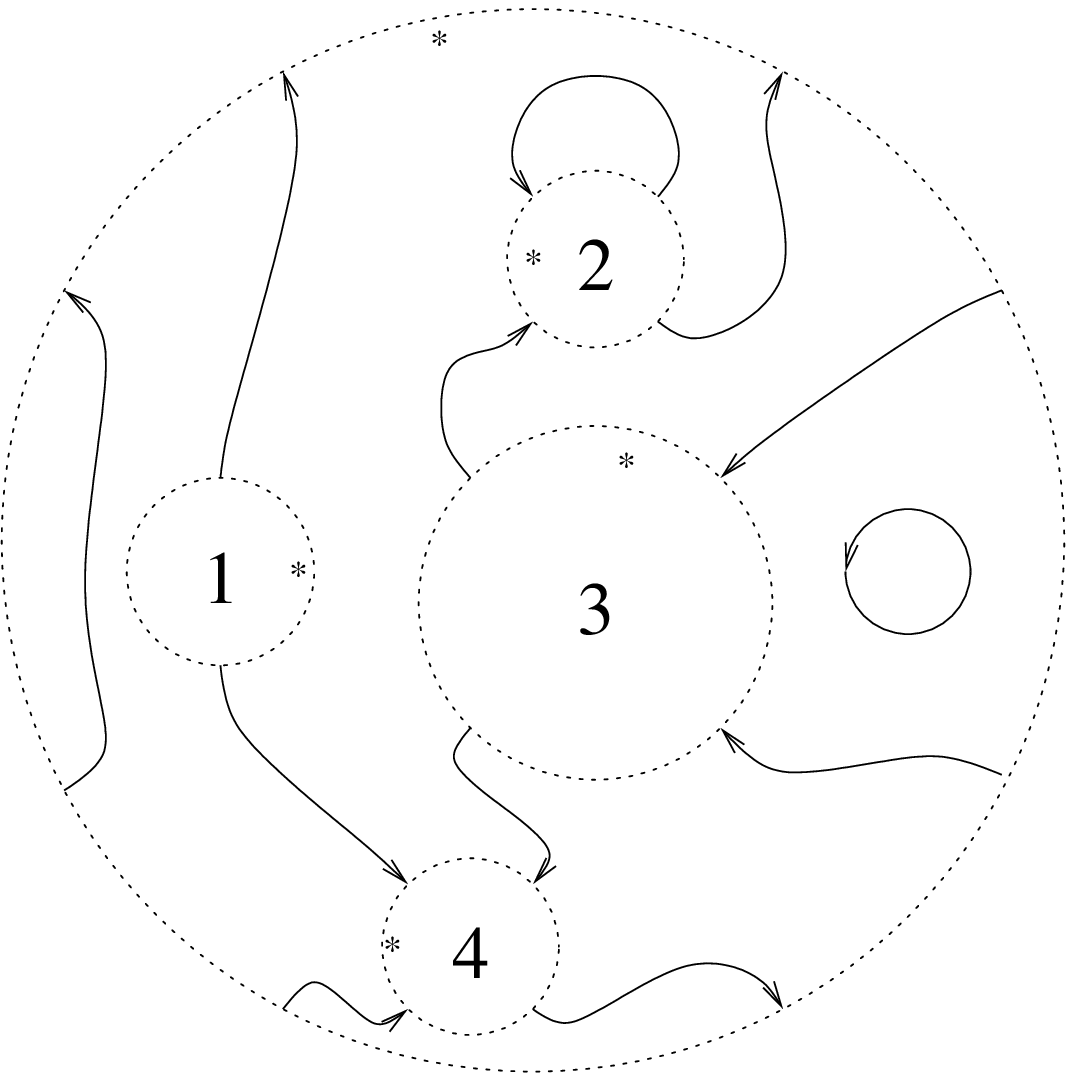}$}
\noindent{\bf Overview}\qua Tangles can be composed in a variety of ways.
Indeed, any $d$--input ``planar arc diagram'' $D$ (such as the $4$--input
example on the right) yields an operator taking $d$ tangles as inputs
and producing a single ``bigger'' tangle as an output by placing the
$d$ input tangles into the $d$ holes of $D$. The purpose of this
section is to define precisely what ``planar arc diagrams'' are and
explain how they turn the collection of tangles into a ``planar
algebra'', to explain how formal complexes in $\Kob$ also form a planar
algebra and to note that the Khovanov bracket
$\llbracket\cdot\rrbracket$ is a planar algebra morphism from the
planar algebra of tangles to the planar algebra of such complexes. Thus
Khovanov brackets ``compose well''. In particular, the invariance
proofs of the previous section, carried out at the local level, lift to
global invariance under Reidemeister moves.

\begin{definition} A $d$--input {\em planar arc diagram} $D$ is a big
``output'' disk with $d$ smaller ``input'' disks removed, along with a
collection of disjoint embedded oriented arcs that are either closed or
begin and end on the boundary. The input disks are numbered $1$ through
$d$, and there is a base point ($\ast$) marked on each of the input disks
as well as on the output disk. Finally, this information is considered
only up to planar isotopy. An {\em unoriented planar arc diagram} is the
same, except the orientation of the arcs is forgotten.
\end{definition}

\parpic(30mm,12mm)[r]{\raisebox{-14mm}{$
  \eps{12mm}{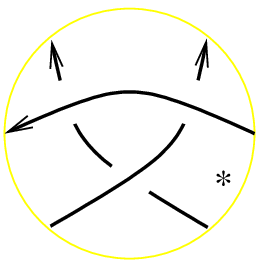}
  \in\calT^0_{\uparrow\downarrow\downarrow\downarrow\uparrow\uparrow}
$}}
\begin{definition}
Let $\calT^0(k)$ denote the collection of all $k$--ended unoriented tangle
diagrams (unoriented tangle diagrams in a disk, with $k$ ends on the
boundary of the disk) in a based disk (a disk with a base point marked on
its boundary). Likewise, if $s$ is a string of in ($\uparrow$) and out
($\downarrow$) symbols with a total length of $|s|$, let $\calT^0(s)$
denote the collection of all $|s|$--ended oriented tangle diagrams in a
based disk with incoming/outgoing strands as specified by $s$, starting
at the base point and going around counterclockwise. Let $\calT(k)$ and
$\calT(s)$ denote the respective quotient of $\calT^0(k)$ and $\calT^0(s)$
by the three Reidemeister moves (so those are spaces of tangles rather
than tangle diagrams).
\end{definition}

Clearly every $d$--input unoriented planar arc diagram $D$ defines
operations (denoted by the same symbol)
\[
  D\co \calT^0(k_1)\times\dots\times\calT^0(k_d)\to\calT^0(k)
  \quad\text{and}\quad
  D\co \calT(k_1)\times\dots\times\calT(k_d)\to\calT(k)
\]
by placing the $d$ input tangles or tangle diagrams into the $d$
holes of $D$ (here $k_i$ are the numbers of arcs in $D$ that end on the
$i$'th input disk and $k$ is the number of arcs that end on the output
disk). Likewise, if $D$ is oriented and $s_i$ and $s$
are the in/out strings read along the inputs and output of $D$ in the
natural manner, then $D$ defines operations
\[
  D\co \calT^0(s_1)\times\dots\times\calT^0(s_d)\to\calT^0(s)
  \quad\text{and}\quad
  D\co \calT(s_1)\times\dots\times\calT(s_d)\to\calT(s)
\]

These operations contain the identity operations on $\calT^{(0)}(k\text{ or }s)$ (take ``radial'' $D$ of the form $\eps{5mm}{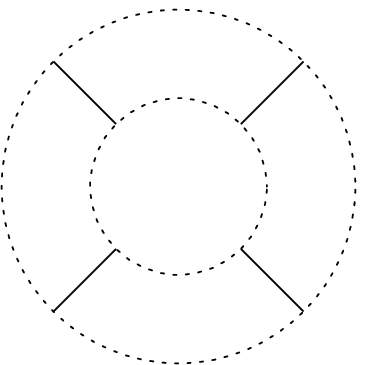}$) and are compatible
with each other (``associative'') in a natural way. In brief,
if $D_i$ is the result of placing $D'$ into the $i$th hole of
$D$ (provided the relevant $k$/$s$ match) then as operations,
$D_i=D\circ(I\times\dots\times D'\times\dots\times I)$.

In the spirit of Jones~\cite{Jones:PlanarAlgebrasI} we call a collection
of sets $\calP(k)$ (or $\calP(s)$) along with operations $D$ defined for each
unoriented planar arc diagram (oriented planar arc diagram) a {\em planar
algebra} (an {\em oriented planar algebra}), provided the radial $D$'s act
as identities and provided associativity conditions as above hold. Thus
as first examples of planar algebras (oriented or not) we can
take $\calT^{(0)}(k\text{ or }s)$.

\parpic[r]{$\eps{1in}{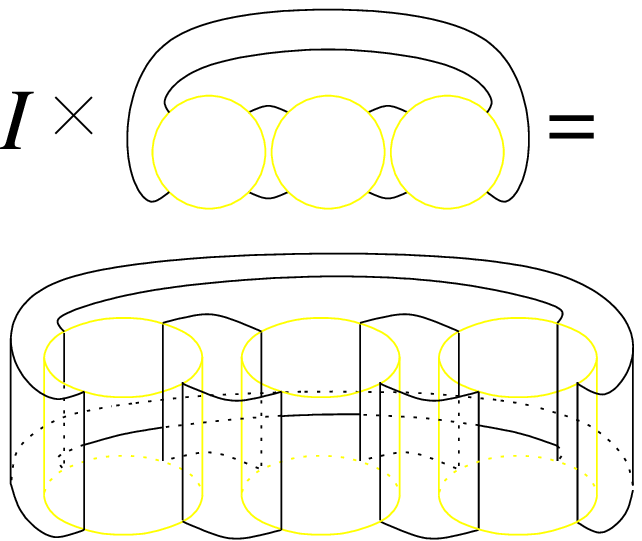}$}
Another example of a planar algebra (unoriented) is the full collection
$\Obj(\Cobl^3)$ of objects of the category $\Cobl^3$ --- this is in fact
the ``flat'' (no crossings) sub planar algebra of $(\calT(k))$. An
even more interesting example is the full collection $\Mor(\Cobl^3)$
of morphisms of $\Cobl^3$ --- indeed, if $D$ is a $d$--input unoriented
planar arc diagram then $D\times[0,1]$ is a vertical cylinder with $d$
vertical cylindrical holes and with vertical curtains connecting those. One
can place $d$ morphisms of $\Cobl^3$ (cobordisms) inside the cylindrical
holes and thus get an operation $D\co (\Mor(\Cobl^3))^d\to\Mor(\Cobl^3)$. Thus
$\Mor(\Cobl^3)$ is also a planar algebra.

A morphism $\Phi$ of planar algebras (oriented or not) $(\calP^a(k))$
and $(\calP^b(k))$ (or $(\calP^a(s))$ and $(\calP^b(s))$) is a collection of maps
(all denoted by the same symbol)
$\Phi\co \calP^a(k\text{ or }s)\to \calP^b(k\text{ or }s)$ satisfying
$\Phi\circ D=D\circ(\Phi\times\dots\times\Phi)$ for every $D$.

We note that every unoriented planar algebra can also be regarded as an
oriented one by setting $\calP(s):=\calP(|s|)$ for every $s$ and by otherwise
ignoring all orientations on planar arc diagrams $D$.

For any natural number $k$ let $\Kob(k):=\Kom(\Mat(\Cobl^3(B_k)))$ and
likewise let $\Kobh(k):=\Komh(\Mat(\Cobl^3(B_k)))$ where
$B_k$ is some placement of $k$ points along a based circle.

\begin{theorem} \label{thm:PlanarAlgebra}$\phantom{99}$
\begin{enumerate}
\item The collection $(\Kob(k))$ has a natural structure of a planar
  algebra.
\item The operations $D$ on $(\Kob(k))$ send homotopy equivalent complexes
  to homotopy equivalent complexes and hence the collection $(\Kobh(k))$
  also has a natural structure of a planar algebra..
\item The Khovanov bracket $\llbracket\cdot\rrbracket$ descends to
  an oriented planar algebra morphism
  $\llbracket\cdot\rrbracket\co (\calT(s))\to(\Kobh(s))$.
\end{enumerate}
\end{theorem}

We note that this theorem along with the results of the previous
section complete the proof of Theorem~\ref{thm:invariance}.

\vskip 2mm
\par\noindent{\bf Abbreviated Proof of Theorem~\ref{thm:PlanarAlgebra}}\qua
The key point is to think of the operations $D$ as (multiple) ``tensor
products'', thus defining these operations on $\Kob$ in analogy with the
standard way of taking the (multiple) tensor product of a number of
complexes.

Start by endowing $\Obj(\Mat(\Cobl^3))$ and $\Mor(\Mat(\Cobl^3))$ with a
planar algebra structure by extending the planar algebra structure of
$\Obj(\Cobl^3)$ and of $\Mor(\Cobl^3)$ in the obvious multilinear manner. Now
if $D$ is a $d$--input planar arc diagram with $k_i$ arcs ending on the
$i$'s input disk and $k$ arcs ending on the outer boundary, and if
$(\Omega_i, d_i)\in\Kob(k_i)$ are complexes, define the complex
$(\Omega,d)=D(\Omega_1,\dots,\Omega_d)$ by
\begin{equation} \label{eq:KobPA}
\begin{split}
    \Omega^r & :=
    \bigoplus_{r=r_1+\dots+r_d} D(\Omega_1^{r_1},\dots,\Omega_d^{r_d})
  \\
    \left.d\right|_{D(\Omega_1^{r_1},\dots,\Omega_d^{r_d})} & :=
    \sum_{i=1}^d (-1)^{\sum_{j<i}r_j}
      D(I_{\Omega_1^{r_1}},\dots,d_i,\dots,I_{\Omega_d^{r_d}}).
\end{split}
\end{equation}
With the definition of $D(\Omega_1,\dots,\Omega_d)$ so similar to the
standard definition of a tensor product of complexes, our reader
should have no difficulty verifying that the basic properties
of tensor products of complexes transfer to our context. Thus
a morphism $\Psi_i\co \Omega_{ia}\to\Omega_{ib}$ induces a morphism
$D(I,\dots,\Psi_i,\dots,I): D(\Omega_1,\dots,\Omega_{ia},\dots,\Omega_d)\to
D(\Omega_1,\dots,\Omega_{ib},\dots,\Omega_d)$ and homotopies at the
level of the tensor factors induce homotopies at the levels of tensor
products. This concludes our abbreviated proof of parts (1) and (2)
of Theorem~\ref{thm:PlanarAlgebra}.

\parpic[r]{$\eps{25mm}{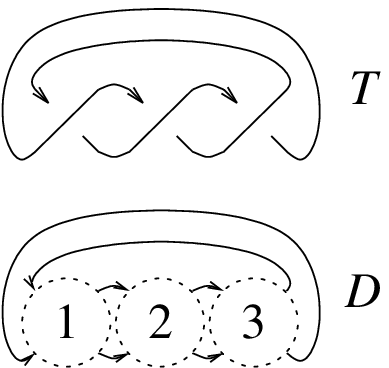}$}
Let $T$ be a tangle diagram with $d$ crossings, let $D$ be the $d$--input
planar arc diagram obtained from $T$ by deleting a disk neighborhood
of each crossing of $T$, let $X_i$ be the $d$ crossings of $T$, so that
each $X_i$ is either an $\overcrossing$ or an $\undercrossing$ (possibly
rotated). Let $\Omega_i$ be the complexes $\llbracket X_i\rrbracket$;
so that each $\Omega_i$ is either $\llbracket\overcrossing\rrbracket
= \left(\xymatrix{\underline\smoothing\ar[r]^\HSaddleSymbol &
\hsmoothing}\right)$ or $\llbracket\undercrossing\rrbracket = \left(
\xymatrix{\hsmoothing\ar[r]^\ISaddleSymbol & \underline\smoothing}\right)$
(possibly rotated, and we've underlined the $0$'th term in each complex). A
quick inspection of the definition of $\llbracket T\rrbracket$ (ie, of
\figref{fig:Main}) and of Equation~\eqref{eq:KobPA} shows that
\[ \llbracket D(X_1,\dots,X_d) \rrbracket
  = \llbracket T\rrbracket
  = D(\Omega_1,\dots,\Omega_d)
  = D(\llbracket X_1\rrbracket,\dots,\llbracket X_d\rrbracket).
\]
This proves part (3) of Theorem~\ref{thm:PlanarAlgebra} in the restricted
case where all inputs are single crossings. The general case follows from
this case and the associativity of the planar algebras involved. \qed

\section{Grading and a minor refinement} \label{sec:grading}

In this short section we introduce gradings into the picture, leading
to a refinement of Theorem~\ref{thm:invariance}. While there isn't any
real additional difficulty in the statement or proof of the refinement
(Theorem~\ref{thm:graded} below), the benefits are great --- the
gradings allow us to relate $\llbracket\cdot\rrbracket$ to the Jones
polynomial (Sections~\ref{sec:Homology} and~\ref{sec:Euler}) and allow
us to easily prove the invariance of the extension of
$\llbracket\cdot\rrbracket$ to 4--dimensional cobordisms
(Section~\ref{sec:EmbeddedCobordisms}).

\begin{definition} A {\em graded category} is a pre-additive category $\calC$
with the following two additional properties:
\begin{enumerate}
\item For any two objects $\calO_{1,2}$ in $\calC$, the morphisms
  $\Mor(\calO_1,\calO_2)$ form a graded Abelian group, the composition
  maps respect the gradings (ie, $\deg f\circ g=\deg f+\deg g$ whenever
  this makes sense) and all identity maps are of degree $0$.
\item There is a $\bbZ$--action $(m,\calO)\mapsto\calO\{m\}$, called
  ``grading shift by $m$'',  on the objects of $\calC$. As plain
  Abelian groups, morphisms are unchanged by this action,
  $\Mor(\calO_1\{m_1\},\calO_2\{m_2\}) = \Mor(\calO_1,\calO_2)$. But
  gradings do change under the action; so if
  $f\in\Mor(\calO_1,\calO_2)$ and $\deg f=d$, then as an element of
  $\Mor(\calO_1\{m_1\},\calO_2\{m_2\})$ the degree of $f$ is $d+m_2-m_1$.
\end{enumerate}
\end{definition}

We note that if an pre-additive category only has the first property above,
it can be `upgraded' to a category $\calC'$ that has the second property
as well. Simply let the objects of $\calC'$ be ``artificial'' $\calO\{m\}$
for every $m\in\bbZ$ and every $\calO\in\Obj(\calC)$ and it is clear
how to define a $\bbZ$--action on $\Obj(\calC')$ and how to define and grade
the morphisms of $\calC'$. In what follows, we will suppress the prime from
$\calC'$ and just call is $\calC$; that is, whenever the morphism groups
are graded, we will allow ourselves to grade-shift the objects of $\calC$.

We also note that if $\calC$ is a graded category then $\Mat(\calC)$ can
also be considered as a graded category (a matrix is considered homogeneous
of degree $d$ iff all its entries are of degree $d$). Complexes in
$\Kom(\calC)$ (or $\Kom(\Mat(\calC))$) become graded in a similar way.

\begin{definition} \label{def:degree}
Let $C\in\Mor(\Cob^3(B))$ be a cobordism in a cylinder, with $|B|$
vertical boundary components on the side of the cylinder. Define $\deg
C:=\chi(C)-\frac12|B|$, where $\chi(C)$ is the Euler characteristic of
$C$.
\end{definition}

\begin{exercise} \label{ex:degrees}
Verify that the degree of a cobordism is additive under vertical
compositions (compositions of morphisms in $\Cob^3(B)$) and under
horizontal compositions (using the planar algebra structure of
Section~\ref{sec:PlanarAlgebra}), and verify that the degree of a
saddle is $-1$ ($\deg\HSaddleSymbol=-1$) and that the degree of a
cap/cup is $+1$ ($\deg\fourwheel=\deg\fourinwheel=+1$). As every
cobordism is a vertical/horizontal composition of copies of
$\HSaddleSymbol$, $\fourwheel$ and $\fourinwheel$, this allows for a
quick computation of degrees.
\end{exercise}

Using the above definition and exercise we know that $\Cob^3$ is a graded
category, and as the $S$, $T$ and \FourTu{} relations are
degree-homogeneous, so is $\Cobl^3$. Hence so are the target categories
of $\llbracket\cdot\rrbracket$, the categories
$\Kobh=\Komh(\Mat(\Cobl^3))$ and $\Kobh=\Kob\!/(\text{homotopy})$.

\begin{definition} Let $T$ be a tangle diagram with $n_+$ positive
crossings and $n_-$ negative crossings. Let $\Kh(T)$ be the complex whose
chain spaces are $\Kh^r(T):=\llbracket T\rrbracket\{r+n_+-n_-\}$ and whose
differentials are the same as those of $\llbracket T\rrbracket$:
%\begin{figure}\anchor{fig:shifts}
\[
  \def\st{\displaystyle}
  \begin{array}{rccccccc}
    \st\llbracket T\rrbracket: & \quad &
      \st\llbracket T\rrbracket^{-n^-} & \st\longrightarrow &
      \st\cdots & \st\longrightarrow &
      \st\llbracket T\rrbracket^{n_+} \\
    &&&&&& \\
    \st\Kh(T): & &
    \st\llbracket T\rrbracket^{-n^-}\{n_+-2n_-\} & \st\longrightarrow &
      \st\cdots & \st\longrightarrow &
      \st\llbracket T\rrbracket^{n_+}\{2n_+-n_-\}
  \end{array}
\]
\end{definition}
\eject
\begin{theorem} \label{thm:graded}$\phantom{99}$
\begin{enumerate}
\item All differentials in $\Kh(T)$ are of degree $0$.
\item $\Kh(T)$ is an invariant of the tangle $T$ up to degree-$0$ homotopy
  equivalences. That is, if $T_1$ and $T_2$ are tangle diagrams which
  differ by some Reidemeister moves, then there is a homotopy equivalence
  $F\co \Kh(T_1)\to\Kh(T_2)$ with $\deg F=0$.
\item Like $\llbracket\cdot\rrbracket$, $\Kh$ descends to an oriented
  planar algebra morphism $(\calT(s))\to(\Kob(s))$, and all the planar
  algebra operations are of degree $0$.
\end{enumerate}
\end{theorem}

\begin{proof} The first assertion follows from $\deg\HSaddleSymbol=-1$ and
from the presence of $r$ in the degree shift $\{r+n_+-n_-\}$ defining
$\Kh$. The second assertion follows from a quick inspection of the
homotopy equivalences in the proofs of invariance under $R1$ and $R2$ in
section~\ref{subsec:Proof}, and the third assertion follows from the
corresponding one for $\llbracket\cdot\rrbracket$ and from the additivity
of $n_+$ and $n_-$ under the planar algebra operations.
\end{proof}

\section{Applying a TQFT and obtaining a homology theory}
\label{sec:Homology}

So $\Kh$ is an up-to-homotopy invariant of tangles, and it has
excellent composition properties. But its target space, $\Kob$, is
quite unmanageable --- given two formal complexes, how can one decide
if they are homotopy equivalent?

In this section we will see how to take the homology of
$\Kh(T)$. In this we lose some of the information
in $\Kh(T)$ and lose its excellent composition
properties. But we get a computable invariant, strong enough to be
interesting.

Let $\calA$ be some arbitrary Abelian category\footnote{You are welcome
to think of $\calA$ as being the category of vector spaces or of
$\bbZ$--modules.}.  Any functor $\calF\co \Cobl^3\to\calA$ extends right
away (by taking formal direct sums into honest direct sums) to a
functor $\calF\co \Mat(\Cobl^3)\to\calA$ and hence to a functor
$\calF\co \Kob\to\Kom(\calA)$. Thus for any tangle diagram $T$,
$\calF\Kh(T)$ is an ordinary complex, and applying $\calF$ to all
homotopies, we see that $\calF\Kh(T)$ is an up-to-homotopy invariant of
the tangle $T$. Thus the isomorphism class of the homology
$H(\calF\Kh(T))$ is an invariant of $T$.

If in addition $\calA$ is graded and the functor
$\calF$ is degree-respecting in the obvious sense, then the homology
$H(\calF\Kh(T))$ is a graded invariant of $T$. And if $\calF$ is only
partially defined, say on $\Cobl^3(\emptyset)$, we get a partially defined
homological invariant --- in the case of $\Cobl^3(\emptyset)$, for example,
its domain will be knots and links rather than arbitrary tangles.

We wish to postpone a fuller discussion of the possible choices for
such a functor $\calF$ to Section~\ref{sec:MoreOnCobl} and just give
the standard example here. Our example for $\calF$ will be a {\em TQFT}
--- a functor on $\Cob^3(\emptyset)$ valued in the category $\ZMod$ of
graded $\bbZ$--modules which maps disjoint unions of to tensor products.
It is enough to define $\calF$ on the generators of
$\Cob^3(\emptyset)$: the object $\bigcirc$ (a single circle) and the
morphisms $\fourwheel$, $\fourinwheel$,
\raisebox{2pt}{$\eps{5mm}{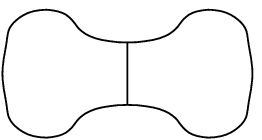}$} and
\raisebox{2pt}{$\eps{5mm}{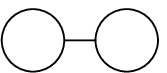}$} (the cap, cup, pair of
pants and upside down pair of pants).

\begin{definition} \label{def:Vdef}
Let $V$ be the graded $\bbZ$--module freely generated by
two elements $\{v_\pm\}$ with $\deg v_\pm=\pm1$. Let $\calF$ be the TQFT
defined by $\calF(\bigcirc)=V$ and by $\calF(\fourwheel)=\epsilon\co \bbZ\to V$,
$\calF(\fourinwheel)=\eta\co V\to\bbZ$,
$\calF(\!\raisebox{2pt}{$\eps{5mm}{POPSymbol}$}\!)=\Delta\co V\to V\otimes V$ and
$\calF(\!\raisebox{2pt}{$\eps{5mm}{InvertedPOPSymbol}$}\!)=m\co V\otimes V\to V$,
where these maps are defined by
\begin{alignat*}{1}
  \calF(\fourwheel)=\epsilon: & \begin{cases}
    1\mapsto v_+
  \end{cases} \\
  \calF(\fourinwheel)=\eta: & \begin{cases}
    v_+\mapsto 0 & v_-\mapsto 1
  \end{cases} \\
  \calF(\raisebox{2pt}{$\eps{5mm}{POPSymbol}$})=
  \Delta: & \begin{cases}
    v_+ \mapsto v_+\otimes v_- + v_-\otimes v_+ &\\
    v_- \mapsto v_-\otimes v_- &
  \end{cases} \\
  \qquad
  \calF(\raisebox{2pt}{$\eps{5mm}{InvertedPOPSymbol}$})=
  m: & \begin{cases}
    v_+\otimes v_-\mapsto v_- &
    v_+\otimes v_+\mapsto v_+ \\
    v_-\otimes v_+\mapsto v_- &
    v_-\otimes v_-\mapsto 0.
  \end{cases}
\end{alignat*}
\end{definition}

\begin{proposition} $\calF$ is well defined and degree-respecting. It
descends to a functor $\Cobl^3(\emptyset)\to\ZMod$.
\end{proposition}

\proof It is well known that $\calF$ is well defined --- ie,
that it respects the relations between our set of generators for
$\Cob^3$, or the relations defining a {\em Frobenius algebra}. See
eg,~\cite{Khovanov:Categorification}. It is easy to verify that $\calF$
is degree-respecting, so it only remains to show that $\calF$ satisfies
the $S$, $T$ and \FourTu{} relations:

\begin{itemize}
\item {$S$.\quad}A sphere is a cap followed by a cup, so we have to show
that $\eta\circ\epsilon=0$. This holds.

\item {$T$.\quad}A torus is a cap followed by a pair of pants followed
by an upside down pair of pants followed by a cup, so we have to compute
$\eta\circ m\circ\Delta\circ\epsilon$. That's not too hard:
\[
  1
  \overset{\epsilon}{\longmapsto} v_+
  \overset{\Delta}{\longmapsto} v_+\otimes v_- + v_-\otimes v_+
  \overset{m}{\longmapsto} v_-+v_-
  \overset{\eta}{\longmapsto} 1+1
  = 2.
\]

\item {$\FourTu$.\quad}We will show that the equality $L=R$ holds in
$V^{\otimes 4}$, where $L$ is given by
$\left(\calF(\eps{14mm}{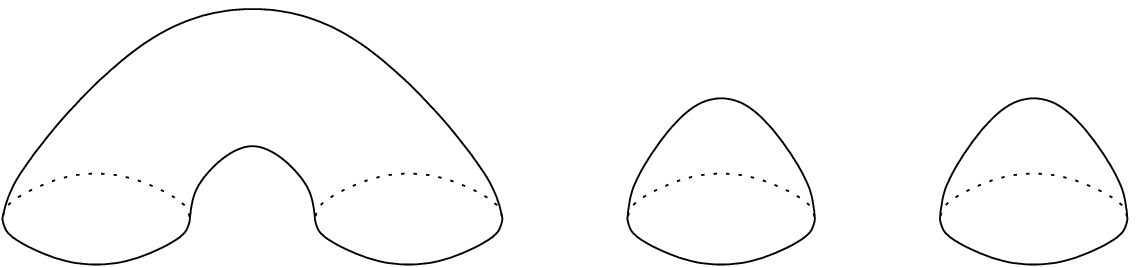})+\calF(\eps{14mm}{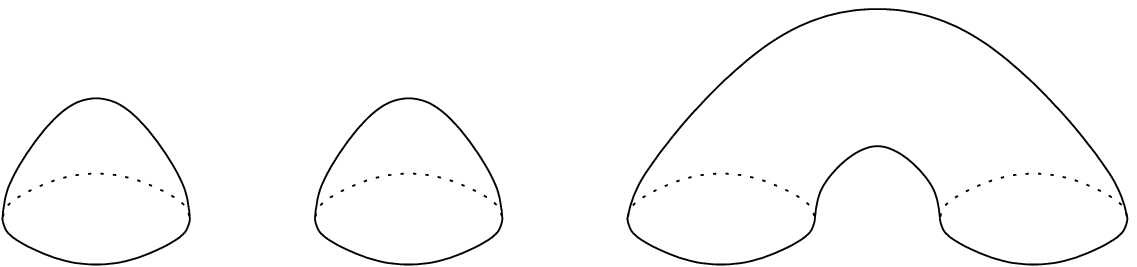})\right)(1)$ and
likewise $R$ is given by
$\left(\calF(\eps{14mm}{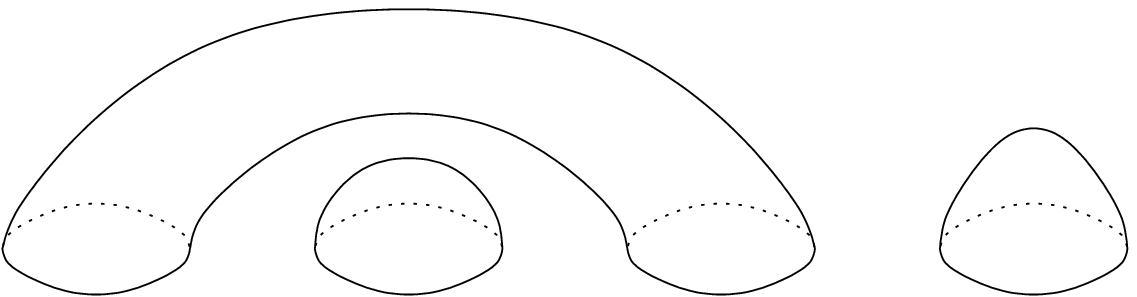})+\calF(\eps{14mm}{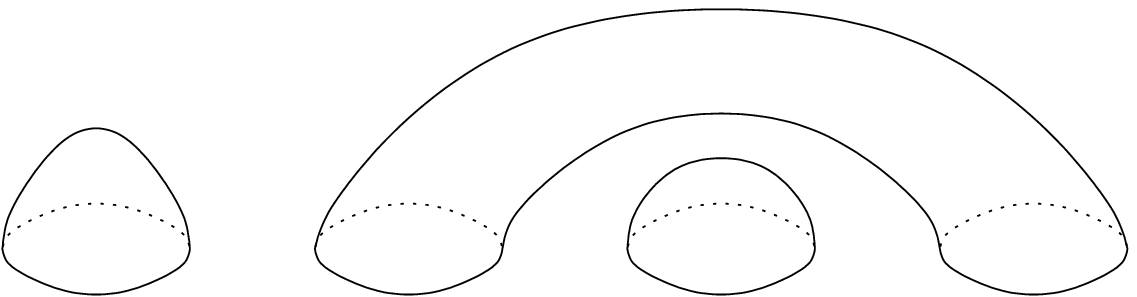})\right)(1)$.
Indeed, $\left(\calF(\eps{14mm}{C12})\right)(1)=\Delta\epsilon
1\otimes\epsilon 1\otimes\epsilon 1=v_-\otimes v_+\otimes v_+\otimes
v_+ + v_+\otimes v_-\otimes v_+\otimes v_+=:v_{-+++}+v_{+-++}$ and
similarly $\left(\calF(\eps{14mm}{C34})\right)(1)=v_{++-+}+v_{+++-}$
and so $L=v_{-+++}+v_{+-++}+v_{++-+}+v_{+++-}$. A similar computation
shows $R$ to be the same.\qed
\end{itemize}

Thus following the discussion at the beginning of this section, we know
that for any $r$ the homology $H^r(\calF\Kh(K))$ is an invariant of the
knot or link $K$ with values in graded $\bbZ$--modules.

A quick comparison of the definitions shows that $H^\star(\calF\Kh(K))$
is equal to Khovanov's categorification of the Jones polynomial and
hence that its graded Euler characteristic is the Jones polynomial $\hatJ$
(see~\cite{Khovanov:Categorification, Bar-Natan:Categorification}). In
my earlier paper~\cite{Bar-Natan:Categorification} I computed
$H^\star(\calF\Kh(K))\otimes\bbQ$ for all prime knots and links with up
to $11$ crossings and found that it is strictly a stronger knot and link
invariant than the Jones polynomial. (See some further computations and
conjectures at \cite{Khovanov:Patterns, Shumakovitch:Torsion}).

\section{Embedded cobordisms} \label{sec:EmbeddedCobordisms}

\subsection{Statement} \label{subsec:ECS}

Let $\Cob^4(\emptyset)$ be the category whose objects are oriented
based\footnote{``Based'' means that one of the crossings is starred.
The only purpose of the basing is to break symmetries and hence to make
the composition of morphisms unambiguous. For most purposes the basing
can be ignored.} knot or link diagrams in the plane, and whose
morphisms are 2--dimensional cobordisms between such knot/link diagrams,
generically embedded in $\bbR^3\times[0,1]$. Let $\Cobi^4(\emptyset)$
be the quotient of $\Cob^4(\emptyset)$ by isotopies (the $/i$ stands
for ``modulo isotopies'').

\parpic[r]{$\eps{2.2in}{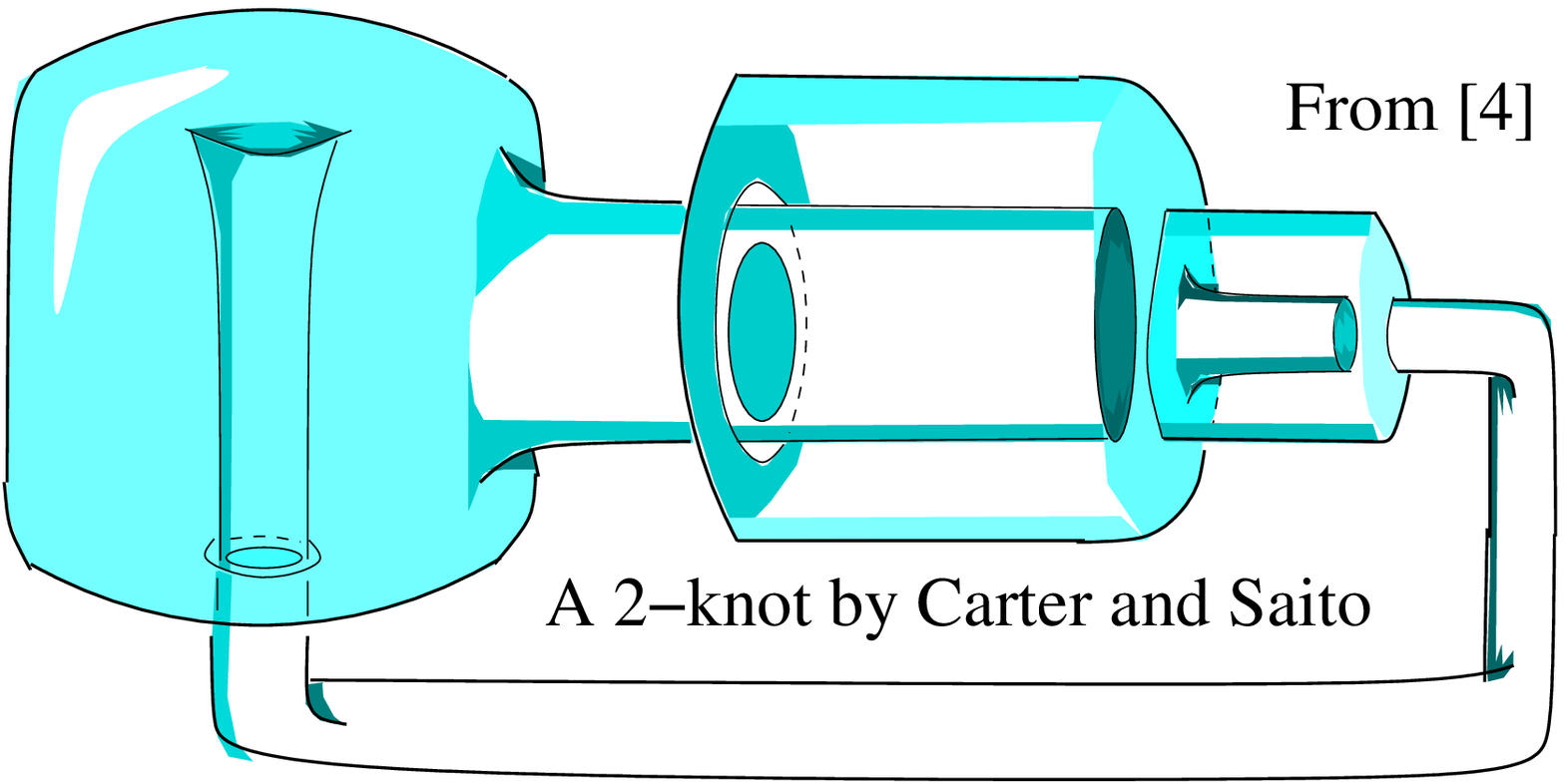}$}
Note that the endomorphisms in $\Cobi^4(\emptyset)$ of the empty link
diagram are simply 2--knots, 2--dimensional knots (or links) in
$\bbR^3\times(0,1)\equiv\bbR^4$. Hence much of what we will say below
specializes to 2--knots. Some wonderful drawings of 2--knots and other
cobordisms in 4--dimensional space are in the book by Carter and
Saito,~\cite{CarterSaito:KnottedSurfaces}.

\begin{figure}\anchor{fig:StringInteractions}
\[ \eps{5in}{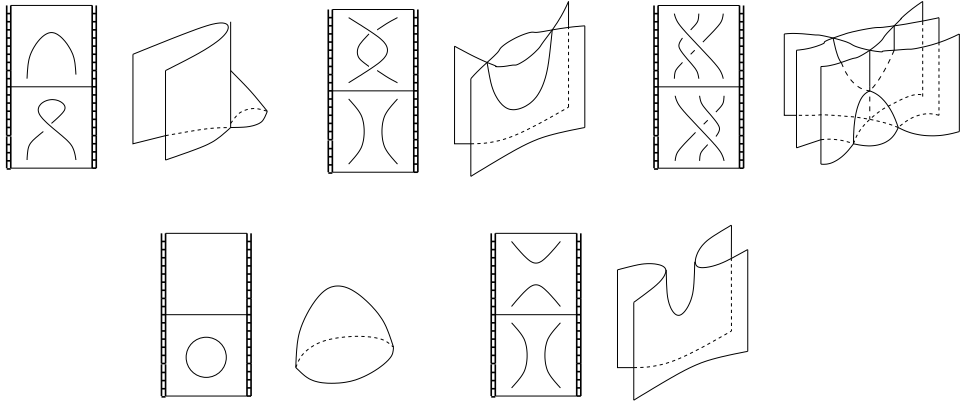} \]
\caption{
  Elementary string interactions as movie clips and 3D projections of
  their 4D realizations, taken from~\cite{CarterSaito:KnottedSurfaces}.
  (All clips are reversible.)
} \label{fig:StringInteractions}
\end{figure}

Thinking of the last coordinate in $\bbR^3\times[0,1]$ as time and
projecting $\bbR^3$ down to the plane, we can think of every cobordism
in $\Cob^4(\emptyset)$ as a movie whose individual frames are
knot/link diagrams (with at most finitely many singular exceptions).
And if we shoot at a sufficiently high frame rate, then between any
consecutive frames we will see (at most) one of the ``elementary
string interactions'' of \figref{fig:StringInteractions} --- a
Reidemeister move, a cap or a cup, or a saddle. Thus the category
$\Cob^4(\emptyset)$ is generated by the cobordisms corresponding to the
three Reidemeister moves and by the cobordisms $\fourwheel$, $\fourinwheel$
and $\HSaddleSymbol$ (now thought of as living in $4D$).

We now define a functor $\Kh\co \Cob^4(\emptyset)\to\Kob(\emptyset)$. On
objects, we've defined $\Kh$ already as the (formal) Khovanov homology of a
given knot/link diagram. On the generating morphisms of
$\Cob^4(\emptyset)$ we define $\Kh$ as
follows:
\begin{itemize}
\item Reidemeister moves go to the chain complex morphisms inducing the
  homotopy equivalences between the `before' and `after' complexes, as
  constructed within the proof of the invariance theorem
  (Theorem~\ref{thm:invariance}) in Section~\ref{subsec:Proof}.
\item The cobordism $\HSaddleSymbol\co \smoothing\to\hsmoothing$ induces a
  morphism $\llbracket\HSaddleSymbol\rrbracket\co 
  \llbracket\smoothing\rrbracket \to \llbracket\hsmoothing\rrbracket$
  just as within the proof of invariance under R3, and just as there it
  can be interpreted in a `skein theoretic' sense, where each symbol
  $\HSaddleSymbol$, $\smoothing$ or $\hsmoothing$ represents a small
  neighborhood within a larger context. But $\llbracket K\rrbracket$
  differs from $\Kh(K)$ only by degree shifts, so the cobordism
  $\HSaddleSymbol\co \smoothing\to\hsmoothing$ also induces a morphism
  $\Kh(\HSaddleSymbol)\co \Kh(\smoothing)\to\Kh(\hsmoothing)$, as required.
\item Likewise, the cobordisms $\fourwheel\co \emptyset\to\bigcirc$ and
  $\fourinwheel\co \bigcirc\to\emptyset$ induce morphisms of complexes
  $\Kh(\fourwheel)\co \Kh(\emptyset)\to\Kh(\bigcirc)$ and
  $\Kh(\fourinwheel)\co \Kh(\bigcirc)\to\Kh(\emptyset)$ (remember to interpret
  all this skein-theoretically --- so the $\emptyset$ symbols here don't
  mean ``the empty set'', but just ``the empty addition to some existing
  knot/link'').
\end{itemize}

\begin{theorem} \label{thm:KhIsFunctor}
Up to signs, $\Kh$ descends to a functor
$\Kh\co \Cobi^4(\emptyset)\to\Kobh(\emptyset)$. Precisely, let
$\PKob$ denote the projectivization of $\Kob$ --- same objects, but every
morphism is identified with its negative, and likewise let $\PKobh$ denote
the projectivization of $\Kobh$. Then $\Kh$ descends to a functor
$\Kh\co \Cobi^4(\emptyset)\to\PKobh(\emptyset)$.
\end{theorem}

The key to the proof of this theorem is to think locally. We need to
show that circular movie clips in the kernel of the projection
$\Cob^4\to\Cobi^4$ (such as the one in Equation~\eqref{eq:MMDemo}) map
to $\pm 1$ in $\Kob(\emptyset)$. As we shall see, the best way to do
so is to view such a clip {\em literally}, as cobordisms between
tangle diagrams, rather than symbolically, as skein-theoretic
fragments of ``bigger'' cobordisms between knot/link diagrams.

Cobordisms between tangle diagrams compose in many ways to produce bigger
cobordisms between tangle diagrams and ultimately to produce cobordisms between
knot/link diagrams or possibly even to produce 2--knots. Cobordisms between
tangle diagrams (presented, say, by movies) can be concatenated to give
longer movies provided the last frame of one movie is equal to the first
frame of the following movie. Thus cobordisms between tangle diagrams form
a category. Cobordisms between tangle diagrams can also be composed like
tangles, by placing them next to each other in the plane and connecting
ends using planar arc diagrams. Hence cobordisms between tangle diagrams
also form a planar algebra.

Thus our first task is to discuss those `things' (called ``canopolies''
below) which are both planar algebras and categories. Ultimately we
will prove that $\Kh$ is a morphism of canopolies between the canopoly
of tangle cobordisms and an appropriate canopoly of formal complexes.

\subsection{Canopolies and a better statement} \label{subsec:Canopolies}

\begin{definition} Let $\calP=(\calP(k))$ be a planar algebra. A
canopoly over $\calP$ is a collection of categories $\calC(k)$ indexed
by the non-negative integers so that $\Obj(\calC(k))=\calP(k)$, so that
the sets $\Mor(\calC(k))$ of all morphisms between all objects of
$\calC(k)$ also form a planar algebra, and so that the planar algebra
operations commute with the category operations (the compositions in
the various categories). A morphism between a canopoly $\calC^1$ over
$\calP^1$ and a canopoly $\calC^2$ over $\calP^2$ is a collection of
functors $\calC^1(k)\to \calC^2(k)$ which also respect all the planar
algebra operations. In a similar manner one may define `oriented'
canopolies $(\calC(s))$ over oriented planar algebras $(\calP(s))$ and
morphisms between such canopolies. Every unoriented canopoly can also
be regarded as an oriented one by setting $\calP(s):=\calP(|s|)$ and
$\calC(s):=\calC(|s|)$ and otherwise ignoring all orientations.
\end{definition}

\parpic[r]{$\includegraphics[width=18mm]{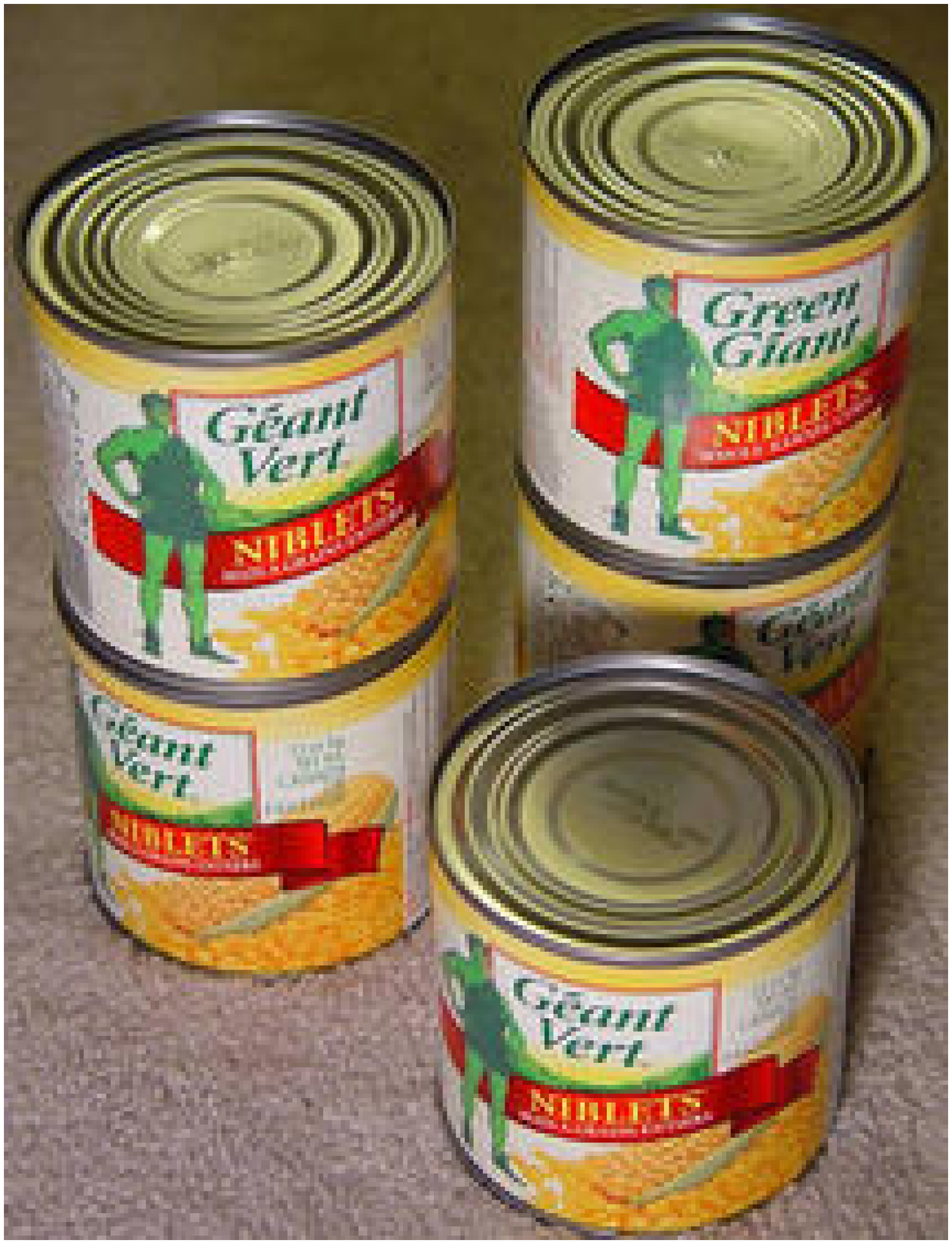}$}
A good way to visualize a canopoly is to think of $(\Mor(\calC(k)))$ as
a collection of `cans' with labels in $(\calP(k))$ on the tops and
bottoms and with $k$ vertical lines on the sides, along with
compositions rules that allow as to compose cans vertically when their
tops/bottoms match and horizontally as in a planar algebra, and so that
the vertical and horizontal compositions commute.

\vskip 1mm
\parpic[r]{$\eps{8mm}{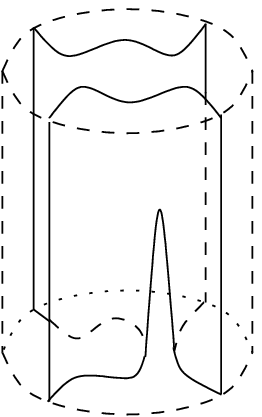}$}
\begin{example} $\Cob^3$ and $\Cobl^3$ are canopolies over the planar
algebra of crossingless tangles. A typical `can' is shown on the right.
\end{example}
\par\vskip 1mm

\parpic[r]{$\eps{1.8in}{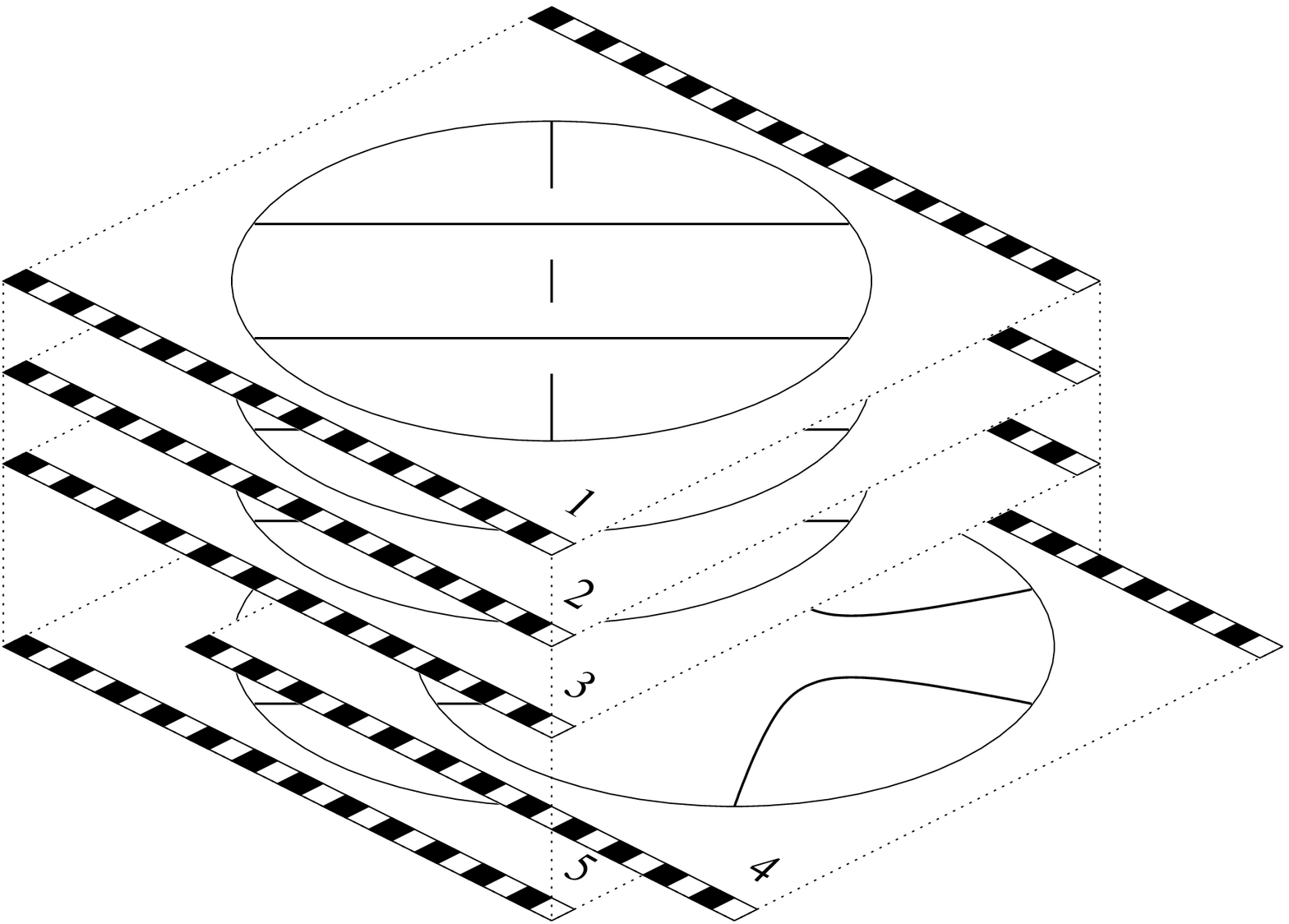}$}
\begin{example} For any finite set $B\subset S^1$ let $\Cob^4(B)$ be
the category whose objects are tangle diagrams in the unit disk $D$
with boundary $B$ and whose morphisms are generic 2--dimensional cobordisms
between such tangle diagrams embedded in
$D\times(-\epsilon,\epsilon)\times[0,1]$ with extra boundary (beyond
the top and the bottom) the vertical lines
$B\times(-\epsilon,\epsilon)\times[0,1]$. For any non-negative $k$, let
$\Cob^4(k)$ be $\Cob^4(B)$ with $B$ some $k$--element set in $S^1$. Then
$\Cob^4:=\bigcup_k\Cob^4(k)$ is a canopoly over the planar algebra of
tangle diagrams. With generic cobordisms visualized as movies, a can in
$\Cob^4$ becomes a vertical stack of frames, each one depicting an
intermediate tangle diagram. In addition, we mod $\Cob^4$ out by
isotopies and call the resulting canopoly $\Cobi^4$.
\end{example}

\begin{example} The collection $\Kob=\bigcup_k\Kob(k)$, previously
regarded only as a planar algebra, can also be viewed as a canopoly. In
this canopoly the `tops' and `bottoms' of cans are formal complexes and
the cans themselves are morphisms between complexes. Likewise
$\Kobh=\Kob/\!(\text{homotopy})$, $\PKob=\Kob/\pm1$ and
$\PKobh:=\Kobh/\pm1$ can be regarded as canopolies.
\end{example}

We note that precisely the same constructions as in
Section~\ref{subsec:ECS}, though replacing the empty boundary $\emptyset$ by
a general $k$ element boundary $B$, define a functor
$\Kh_0\co \Cob^4(B)\to\Kob(B)$ for any $B$. As these constructions are local,
it is clear that these functors assemble together to form a canopoly
morphism $\Kh_0\co \Cob^4\to\Kob$ from the canopoly of movie presentations of
four dimensional cobordisms between tangle diagrams to the canopoly of
formal complexes and morphisms between them.

We also note that the notion of a {\em graded canopoly} can be defined
along the lines of Section~\ref{sec:grading} --- grade the cans (but
not the planar algebras of the ``tops'' and ``bottoms'') and insist
that all the can composition operations be degree-additive. One easily
verifies that all the above mentioned canopolies are in fact graded,
with the gradings induced from the gradings of $\Cob^3$ and of $\Cob^4$
($\Cob^3$ was given a grading in Definition~\ref{def:degree} and
Exercise~\ref{ex:degrees}, and the same definition and exercise can be
applied without changing a word to $\Cob^4$). Clearly $\Kh_0$ is degree
preserving.

The following theorem obviously generalizes Theorem~\ref{thm:KhIsFunctor}
and is easier to prove:

\begin{theorem} \label{thm:Main} $\Kh_0$ descends to a degree preserving
canopoly morphism $\Kh: \Cobi^4\to\PKobh$ from the canopoly of four
dimensional cobordisms between tangle diagrams to the canopoly of
formal complexes with up to sign and up to homotopy morphisms between
them.
\end{theorem}

\subsection{Proof} \label{subsec:CobInvProof}
\begin{figure}\anchor{fig:MM1-5}
\[ \eps{4in}{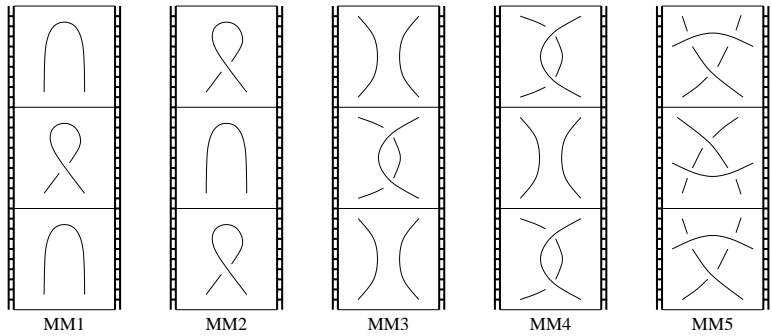} \]
\caption{
  Movie moves as in Carter and
  Saito~\cite{CarterSaito:KnottedSurfaces}. Type I: Reidemeister and
  inverses. These short clips are equivalent to ``do nothing'' identity
  clips.
} \label{fig:MM1-5}
\end{figure}

\begin{figure}\anchor{fig:MM6-10}
\[ \eps{4.8in}{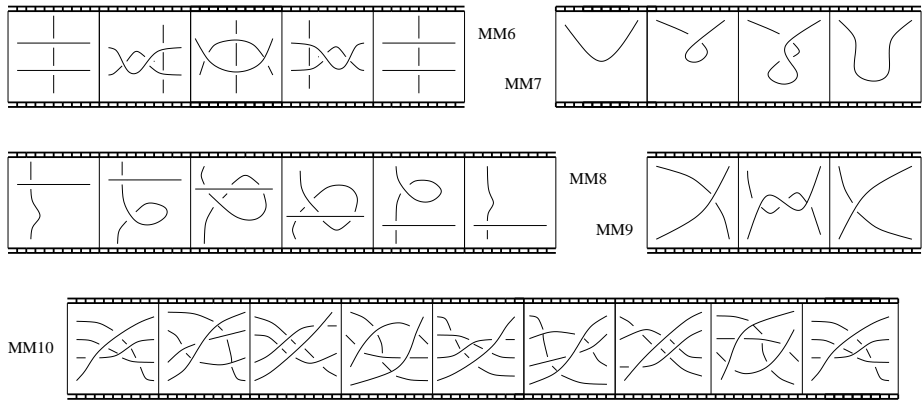} \]
\caption{
  Movie moves as in Carter and
  Saito~\cite{CarterSaito:KnottedSurfaces}. Type II: Reversible
  circular clips --- equivalent to identity clips.
} \label{fig:MM6-10}
\end{figure}

\begin{figure}\anchor{fig:MM11-15}
\[ \eps{\hsize}{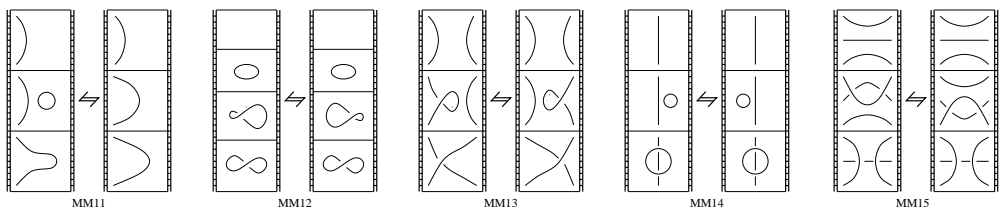} \]
\caption{
  Movie moves as in Carter and
  Saito~\cite{CarterSaito:KnottedSurfaces}. Type III: Non-reversible
  clips (can be read both from the top down and from the bottom up).
} \label{fig:MM11-15}
\end{figure}

We just need to show that $\Kh_0$ respects the relations in the kernel
of the projection $\Cob^4\to\Cobi^4$. These are the ``movie moves'' of
Carter and Saito~\cite{CarterSaito:KnottedSurfaces}, reproduced here in
Figures~\fref{fig:MM1-5}, \fref{fig:MM6-10} and~\fref{fig:MM11-15}. In
principle, this is a routine verification. All that one needs to do is
to write down explicitly the morphism of complexes corresponding to
each of the clips in those figures, and to verify that these morphisms
are homotopic to identity morphisms (in some cases) or to each other
(in other cases).

But this isn't as simple as it sounds, as many of the complexes
involved are quite complicated. The worst is of course $\MM_{10}$ of
\figref{fig:MM6-10} --- each frame in that clip involves a
6-crossing tangle, and hence a 6-dimensional cube of 64 smoothings, and
each of the 8 moves in $\MM_{10}$ is an $R3$ move, so the morphism
corresponding to it originates from the morphism displayed in
\figref{fig:R3Full}. Even if in principle routine, it obviously
isn't a simple task to show that the composition of 8 such beasts is
homotopic to the identity automorphism (of a 6-dimensional cube).

This is essentially the approach taken by Jacobsson
in~\cite{Jacobsson:Cobordisms}, where he was able to use clever tricks
and clever notation to reduce this complexity significantly, though
much complexity remains. At the end of the day the theorem is proven by
carrying out a number of long computations, but it remains a mystery
whether these computations {\em had} to work out, or is it just a
concurrence of lucky coincidences.

Our proof of Theorem~\ref{thm:Main} is completely different, though it
is very similar in spirit to Khovanov's
proof~\cite{Khovanov:Cobordisms}.  The key to our proof is the fact
that the complexes corresponding to many of the tangles appearing in
Figures~\fref{fig:MM1-5}, \fref{fig:MM6-10} and~\fref{fig:MM11-15} simply
have no automorphisms other than up-to-homotopy $\pm 1$ multiples of
the identity, and hence $\Kh_0$ has no choice but to send the clips in
Figures~\fref{fig:MM1-5} and~\fref{fig:MM6-10} to up-to-homotopy $\pm 1$
multiples of the identity.

We start with a formal definition of ``no automorphisms'' and then prove 4
short lemmas that together show that there are indeed many tangles whose
corresponding complexes have ``no automorphisms'':

\begin{definition} We say that a tangle diagram $T$ is $\Kh$--simple if
every degree $0$ automorphism of $\Kh(T)$ is homotopic to a $\pm 1$
multiple of the identity. (An automorphism, in this context, is a homotopy
equivalence of $\Kh(T)$ with itself).
\end{definition}

\vskip -5mm

\parpic(0.35in,0.35in)[r]{\raisebox{-12mm}{
  $\eps{0.35in}{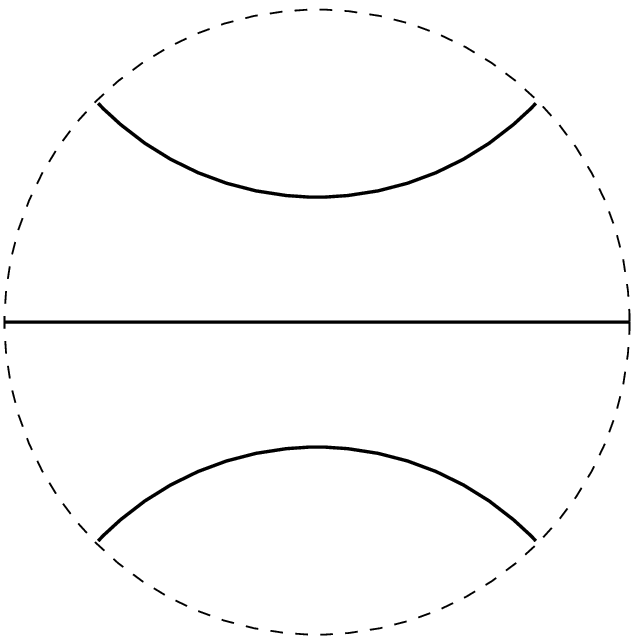}$
}}
%\begin{floatingfigure}[r]{0.35in}
%  \hspace{-5mm}\includegraphics[width=0.35in]{figs/Pairing.eps}
%\end{floatingfigure}
\begin{lemma} \label{lem:Pairings}
Pairings are $\Kh$--simple (a pairing is a tangle that has no
crossings and no closed components, so it is just a planar pairing of its
boundary points).
\end{lemma}

\begin{proof} If $T$ is a pairing then $\Kh(T)$ is the $0$--dimensional
cube of the $2^0$ smoothings of $T$ --- namely, it is merely the one
step complex consisting of $T$ alone at height $0$ and of no
differentials at all. A degree $0$ automorphism of this complex is a
formal $\bbZ$--linear combination of degree $0$ cobordisms with
top and bottom equal to $T$.

\parpic[r]{$\eps{0.85in}{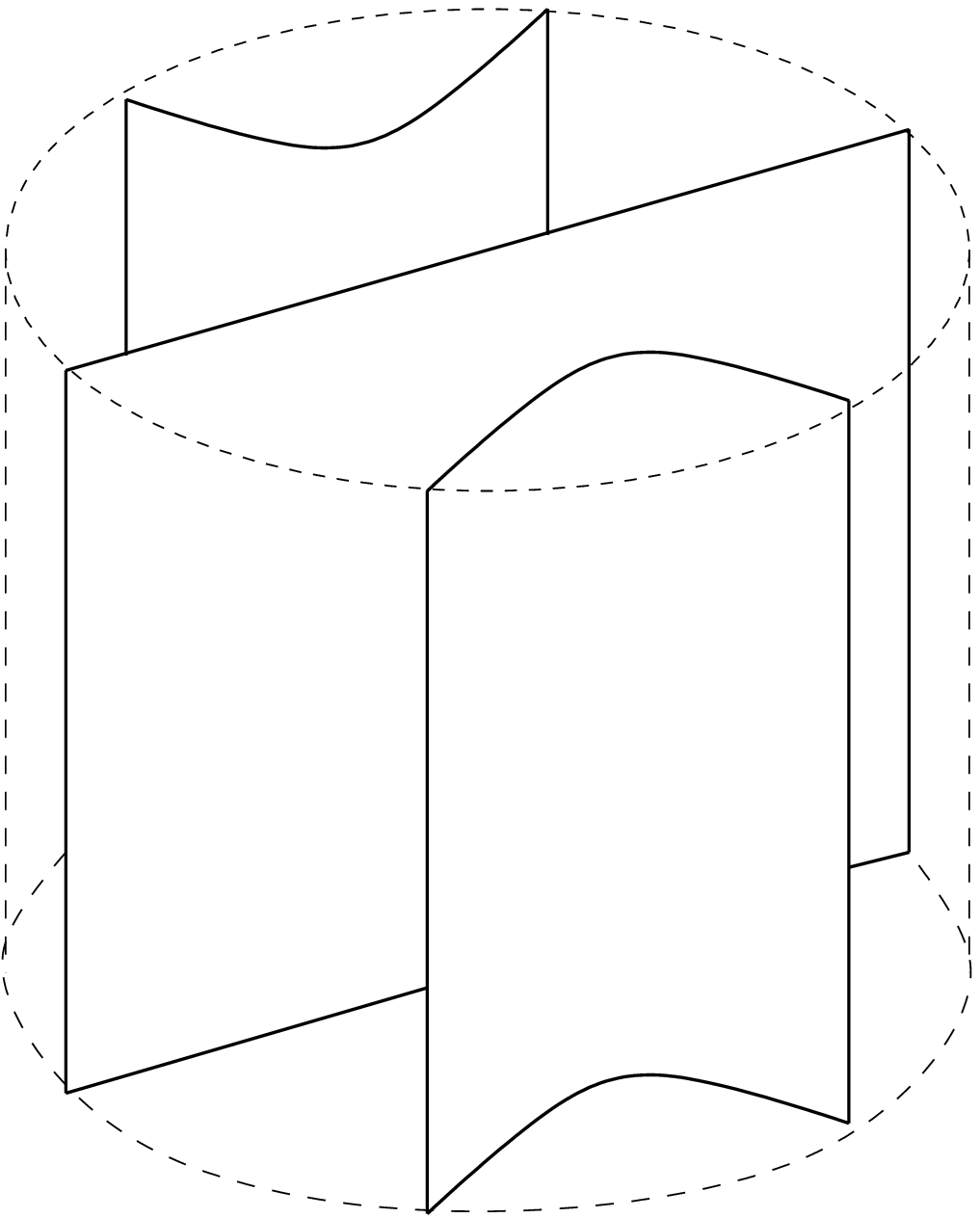}$}
Let us take one such cobordism and call it $C$. By the definition of
degrees in Section~\ref{sec:grading} it follows that $C$ must have
Euler characteristic equal to the number of its boundary components (which
is the same as the number of components of $T$ and half the number of
boundary points of $T$). If $C$ has no connected components with no
boundary, this forces $C$ to be a union of disks embedded vertically
(as ``curtains'') as on the right. Any tori (whose Euler characteristic
is $0$) in $C$ can be reduced using the $T$ relation and any higher
genus boundaryless components (with negative Euler characteristic) must
be balanced against spherical components (whose Euler characteristic is
positive). But the latter are $0$ by the $S$ relation.

Hence $C$ is the identity and so $\Kh(T)$ is a multiple of the identity.
But being invertible it must therefore be a $\pm 1$ multiple of the
identity.
\end{proof}

\begin{lemma} \label{lem:KhSimpleIsotopy}
If a tangle diagram $T$ is $\Kh$--simple and a tangle diagram
$T'$ represents an isotopic tangle, then $T'$ is also $\Kh$--simple.
\end{lemma}

\parpic[r]{$\displaystyle\xymatrix{
  \Kh(T') \ar@<2pt>[r]^F \ar[d]^{\alpha'} &
  \Kh(T) \ar@<2pt>[l]^G \ar[d]^\alpha \\
  \Kh(T') \ar@<2pt>[r]^F &
  \Kh(T) \ar@<2pt>[l]^G 
}$}
\begin{proof} By the invariance of $\Kh$ (Theorems~\ref{thm:invariance}
and~\ref{thm:PlanarAlgebra}), the complexes $\Kh(T')$ and $\Kh(T)$ are
homotopy equivalent. Choose a homotopy equivalence $F\co \Kh(T')\to\Kh(T)$
between the two (with up-to-homotopy inverse $G\co \Kh(T)\to\Kh(T')$), and
assume $\alpha'$ is a degree $0$ automorphism of $\Kh(T')$. As $T$ is
$\Kh$--simple, we know that $\alpha:=F\alpha'G$ is homotopic to $\pm I$
and so $\alpha'\sim GF\alpha'GF=G\alpha F\sim\pm GF\sim\pm I$, and so
$T'$ is also $\Kh$--simple.
\end{proof}

\parpic[r]{$\eps{0.45in}{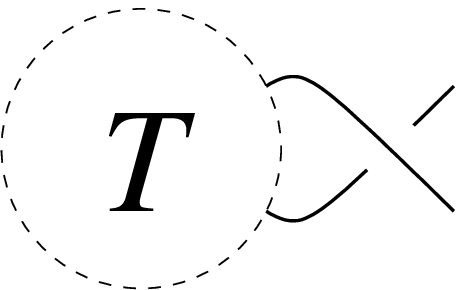}$}
Now let $T$ be a tangle and let $\TX$ be a tangle obtained from $T$ by
adding one extra crossing $X$ somewhere along the boundary of $T$, so
that exactly two (adjacent) legs of $X$ are connected to $T$ and two
remain free. This operation of adjoining a crossing is ``invertible'';
one can adjoin the inverse crossing $X^{-1}$ to get $\TX\!X^{-1}$ which is
isotopic to the original tangle $T$. This fact is utilized in the following
two lemmas to show that $T$ is $\Kh$--simple iff $\TX$ is $\Kh$--simple.

\parpic[r]{$\eps{0.8in}{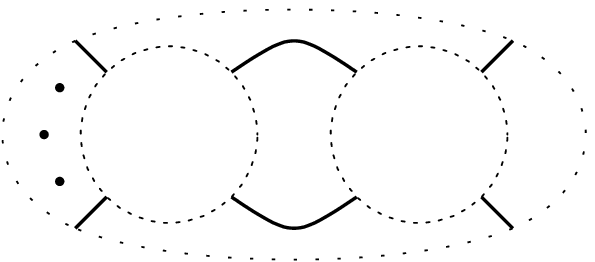}$}
But first a word about notation. In a canopoly there are many operations
and two of them will be used in the following proof. Within that proof
we will denote the horizontal composition of putting things side by
side by simple juxtaposition (more precisely, this is the planar
algebra operation corresponding to the planar arc diagram on the
right). Thus starting with a tangle $T$ and a crossing $X$ we get $\TX$
as in the previous paragraph.  The vertical composition of putting
things one atop the other will be denoted by $\circ$.

\begin{lemma} \label{lem:TXtoT}
If $\TX$ is $\Kh$--simple then so is $T$.
\end{lemma}

\begin{proof} Let $\alpha\co \Kh(T)\to\Kh(T)$ be a degree $0$
automorphism.  Using the planar algebra operations we adjoin a crossing
$X$ on the right to $T$ and to $\alpha$ to get an automorphism $\alpha
I_X\co \Kh(\TX)\to\Kh(\TX)$ (here $I_X$ denotes the identity automorphism
of $\Kh(X)$). As $\TX$ is $\Kh$--simple, $\alpha I_X\sim\pm I$.  We now
adjoin the inverse $X^{-1}$ of $X$ and find that $\alpha
I_{X\!X^{-1}}\co \Kh(\TX\!X^{-1})\to\Kh(\TX\!X^{-1})\sim\pm I$. But $X\!X^{-1}$ is
differs from $\hsmoothing$ by merely a Reidemeister move, so let
$F\co \Kh(\hsmoothing)\to\Kh(X\!X^{-1})$ be the homotopy equivalence between
$\Kh(\hsmoothing)$ and $\Kh(X\!X^{-1})$ (with up-to-homotopy inverse
$G\co \Kh(X\!X^{-1})\to\Kh(\hsmoothing)$). Now
$\alpha
  =\alpha I_{\Kh(\hsmoothing)}
  \sim\alpha(G\circ F)
  =(I_{\Kh(T)}G)\circ(\alpha I_{X\!X^{-1}})\circ(I_{\Kh(T)}F)
  \sim\pm(I_{\Kh(T)}G)\circ I\circ(I_{\Kh(T)}F)
  =\pm I_{\Kh(T)}(G\circ F)
  \sim \pm I,
$ and so $T$ is $\Kh$--simple.
\end{proof}

\begin{lemma} \label{lem:TtoTX}
If $T$ is $\Kh$--simple then so is $\TX$.
\end{lemma}

\begin{proof} Assume $T$ is $\Kh$--simple. By
Lemma~\ref{lem:KhSimpleIsotopy} so is $\TX\!X^{-1}$. Using
Lemma~\ref{lem:TXtoT} we can drop one crossing, the $X^{-1}$, and find that
$\TX$ is $\Kh$--simple.
\end{proof}

\vskip 2mm

We can finally get back to the proof of Theorem~\ref{thm:Main}. Recall
that we have to show that $\Kh_0$ respects each of the movie moves
$\MM_1$ through $\MM_{15}$ displayed in Figures~\fref{fig:MM1-5},
\fref{fig:MM6-10} and~\fref{fig:MM11-15}. These movie moves can be
divided into three types as follows.

{\bf Type I}\qua Performing a Reidemeister and then its inverse
(\figref{fig:MM1-5}) is the same as doing nothing. Applying $\Kh_0$
we clearly get morphisms that are homotopic to the identity --- this is
precisely the content of Theorem~\ref{thm:graded}.

\parpic[r]{$\eps{0.6in}{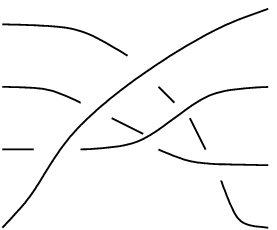}$}
{\bf Type II}\qua The reversible circular movie moves (``circular clips'')
of \figref{fig:MM6-10} are equivalent to the ``do nothing'' clips
that have the same initial and final scenes. This is the hardest and
easiest type of movie moves --- the hardest because it includes the
vicious $\MM_{10}$. The easiest because given our machinery the proof
reduces to just a few sentences. Indeed, $\Kh_0(\MM_{10})$ is an
automorphism of $\Kh(T)$ where $T$ is the tangle beginning and ending
the clip $\MM_{10}$. But using Lemma~\ref{lem:TtoTX} we can peel off
the crossings of $T$ one by one until we are left with a crossingless
tangle (in fact, a pairing), which by Lemma~\ref{lem:Pairings} is
$\Kh$--simple. So $T$ is also $\Kh$--simple and hence
$\Kh_0(\MM_{10})\sim\pm I$ as required. The same argument works for
$\MM_6$ through $\MM_9$. (And in fact, the same argument also works for
$\MM_1$ through $\MM_5$, though as seen above, these movie moves afford
an even easier argument).

{\bf Type III}\qua The pairs of equivalent clips appearing in
\figref{fig:MM11-15}. With some additional effort one can adapt the
proof for the type II movie moves to work here as well, but given the
low number of crossings involved, the brute force approach becomes
sufficiently gentle here. Indeed, we argue as follows.

\begin{myitemize}

\parpic[r]{\raisebox{-16mm}{$\eps{0.8in}{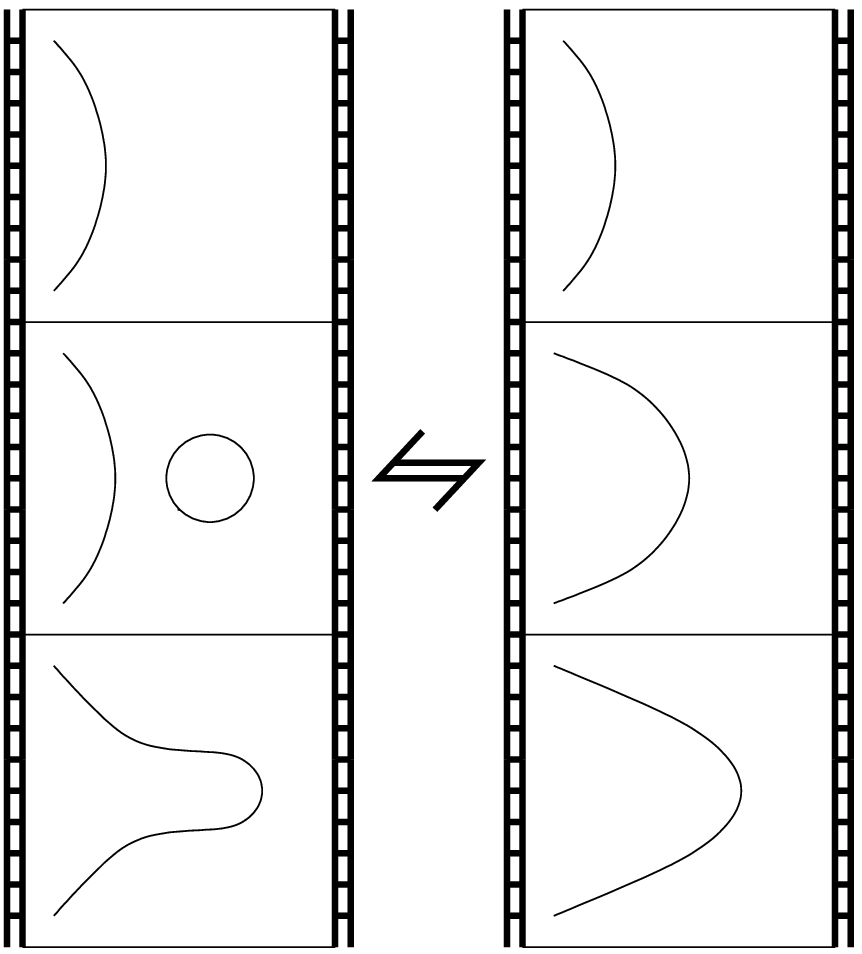}$}}
\item For $\MM_{11}$: Going down along the left side of
$\MM_{11}$ we get a morphism
$F\co \Kh\left(\eps{5mm}{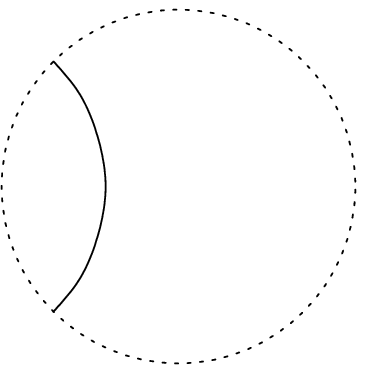}\right)\to\Kh\left(\eps{5mm}{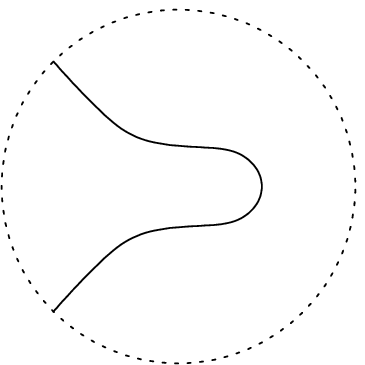}\right)$.
Both $\Kh\left(\eps{5mm}{MM11P1}\right)$ and
$\Kh\left(\eps{5mm}{MM11P2}\right)$ are one-step complexes,
$\underline{\eps{5mm}{MM11P1}}$ and $\underline{\eps{5mm}{MM11P2}}$
respectively, and $F$ is just the cobordism
$\eps{5mm}{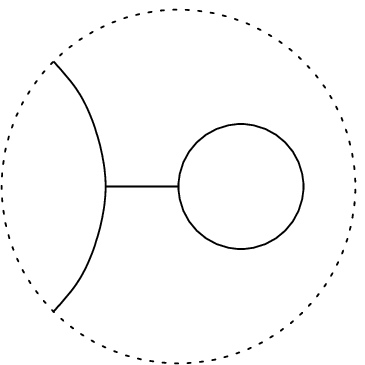}\circ\eps{5mm}{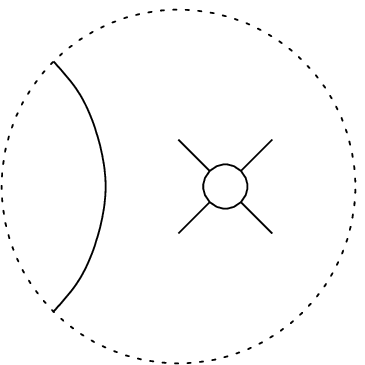}=\eps{7mm}{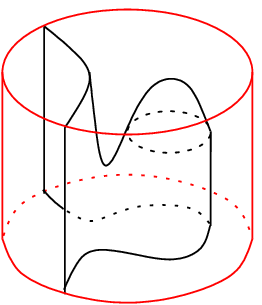}$ which is
isotopic to the identity cobordism
$\eps{5mm}{MM11P1}\to\eps{5mm}{MM11P2}$.  Going up along $\MM_{11}$,
you just have to turn all these figures upside down.

\vskip 2mm \parpic[r]{\raisebox{-16mm}{$\eps{0.8in}{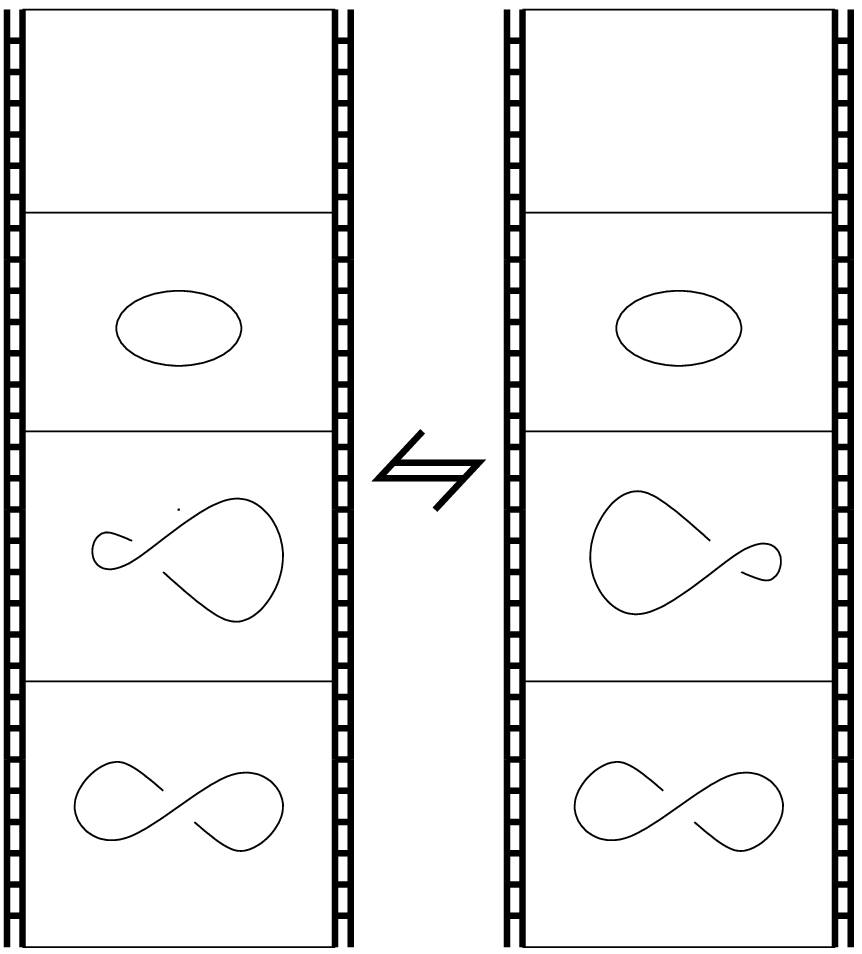}$}}
\item For $\MM_{12}$: At the top of $\MM_{12}$ we see the
empty tangle $\emptyset$ and $\Kh(\emptyset)=(\emptyset)$ is the
one-step complex whose only ``chain group'' is the empty smoothing
$\emptyset$. At the bottom,
$\Kh\left(\raisebox{1mm}{$\eps{5mm}{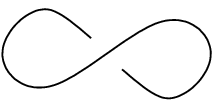}$}\right)$ is the two step
complex $\underline{\bigcirc\bigcirc}\to\bigcirc$. Going down the right
side of $\MM_{12}$ starting at $\Kh(\emptyset)$ we land in the height
$0$ part $\bigcirc\bigcirc$ of
$\Kh\left(\raisebox{1mm}{$\eps{5mm}{MM12P1}$}\right)$ and as can be
seen from the proof of invariance under $R1$, the resulting morphism
$\emptyset\to\bigcirc\bigcirc$ is the difference
$F_R=\eps{7mm}{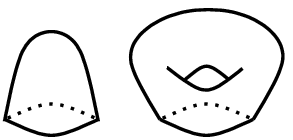}-\eps{7mm}{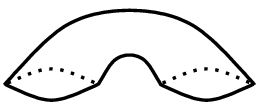}$. Likewise going along the left
side of $\MM_{12}$ we get the difference
$F_L=\eps{7mm}{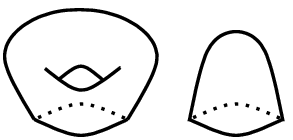}-\eps{7mm}{MM12P3}$. We leave it as an exercise to
the reader to verify that modulo the $\FourTu$ relation $F_L+F_R=0$ (hint:
Equation~\eqref{eq:CuttingNecks} below) and
hence the two side of $\MM_{12}$ are the same up to a sign. Going up
$\MM_{12}$ is even easier.

\vskip 2mm \parpic[r]{\raisebox{-16mm}{$\eps{0.8in}{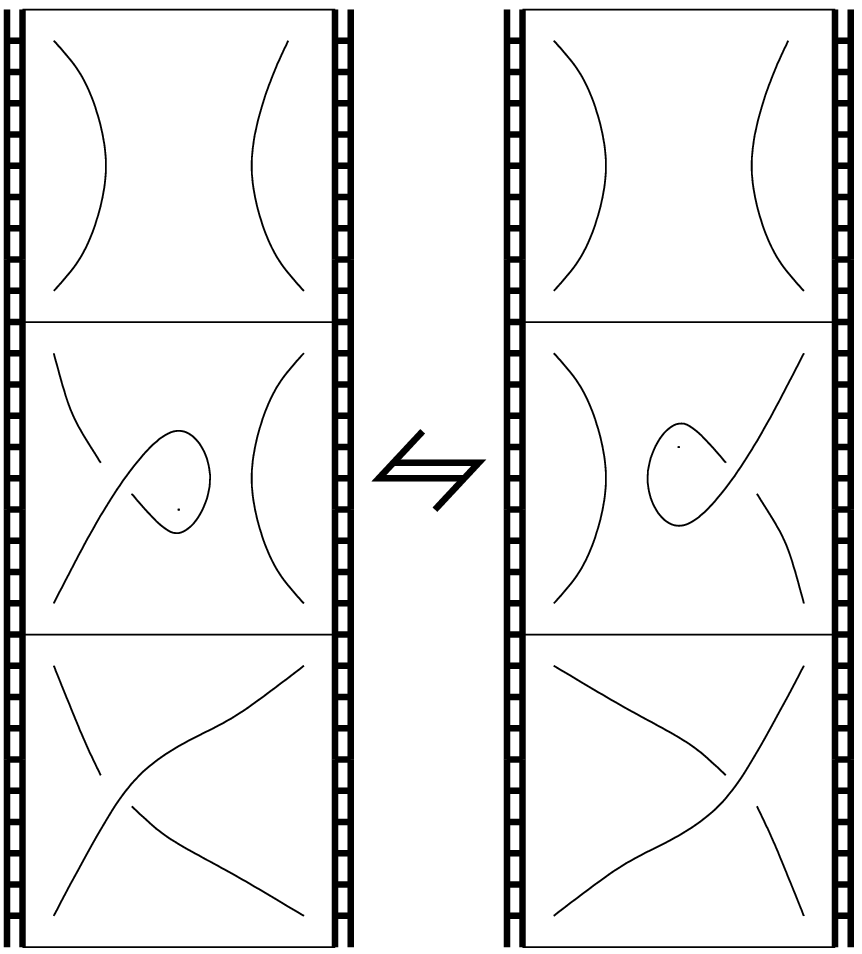}$}}
\item For $\MM_{13}$: The height 0 part of $\Kh(\slashoverback)$ is
$\smoothing$ so going down the two sides of $\MM_{13}$ we get two
morphisms $\smoothing\to\smoothing$, and both are obtained from the
morphism $F^0$ of \figref{fig:R1Invariance} by composing it with an
extra saddle. Tracing it through, we find that the left morphism is
$F_L=\eps{6mm}{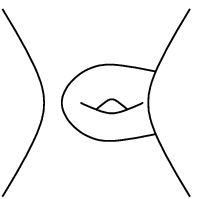}-\eps{6mm}{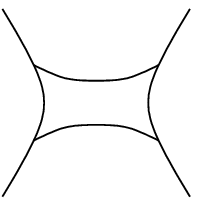}$ and the right morphism is
$F_R=\eps{6mm}{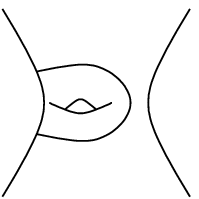}-\eps{6mm}{MM13P2}$ (here all cobordisms are
shown from above and $\eps{6mm}{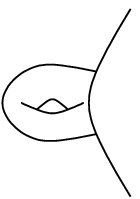}$ denotes a vertical curtain
with an extra handle attached and $\eps{6mm}{MM13P2}$ denotes two
vertical curtain connected by a horizontal tube). We leave it as an
exercise to the reader to verify that modulo the $\FourTu$ relation
$F_L+F_R=0$ (hint:~\eqref{eq:CuttingNecks}) and hence the two side of
$\MM_{13}$ are the same up to a sign. Going up $\MM_{13}$ is even
easier.

\vskip 2mm \parpic[r]{\raisebox{-16mm}{$\eps{0.8in}{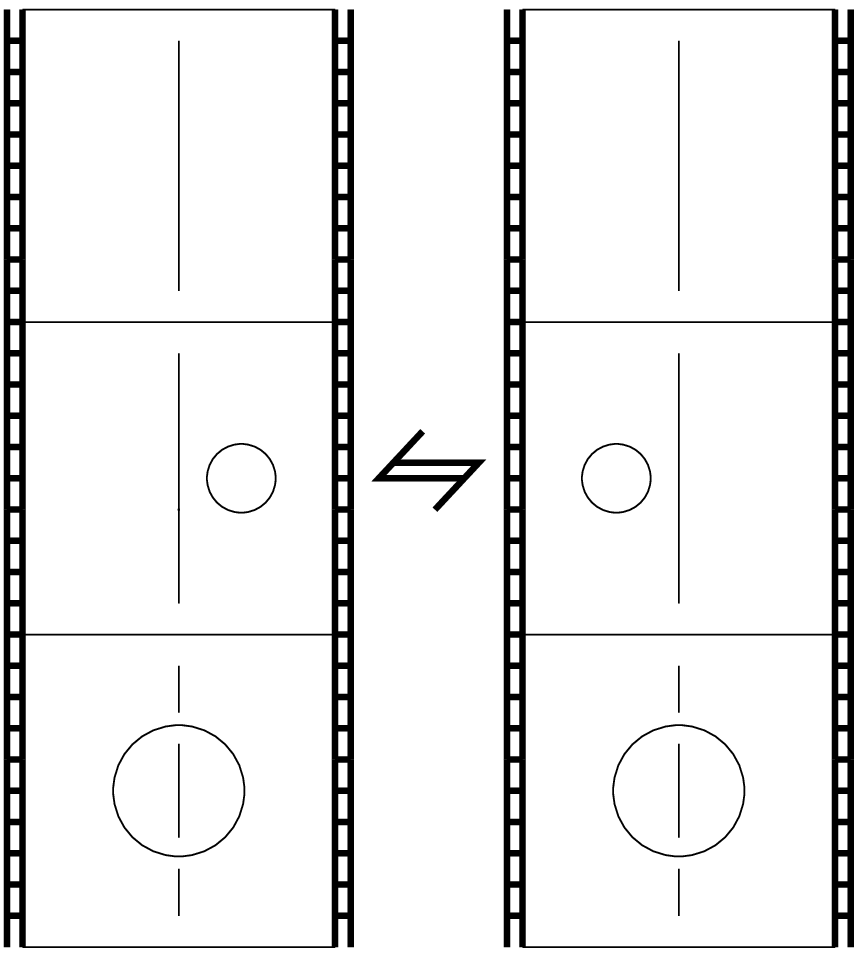}$}}
\item For $\MM_{14}$: The two height 0 parts of
$\Kh\left(\!\eps{3mm}{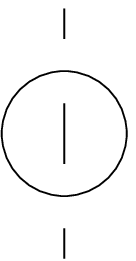}\!\right)$ are $|\circ$ and $\circ|$,
and using the map $F$ of \figref{fig:Reid2Proof} we see that the
four morphisms $|\to|\circ$ and $|\to\circ|$ obtained by tracing
$\MM_{14}$ from the top to the bottom either on the left or on the
right are all simple ``circle creation'' cobordisms (ie, caps) along
with a vertical curtain. In particular, the left side and the right
side of $\MM_{14}$ produce the same answer. A similar argument works
for the way up $\MM_{14}$.

\vskip 2mm \parpic[r]{\raisebox{-16mm}{$\eps{0.8in}{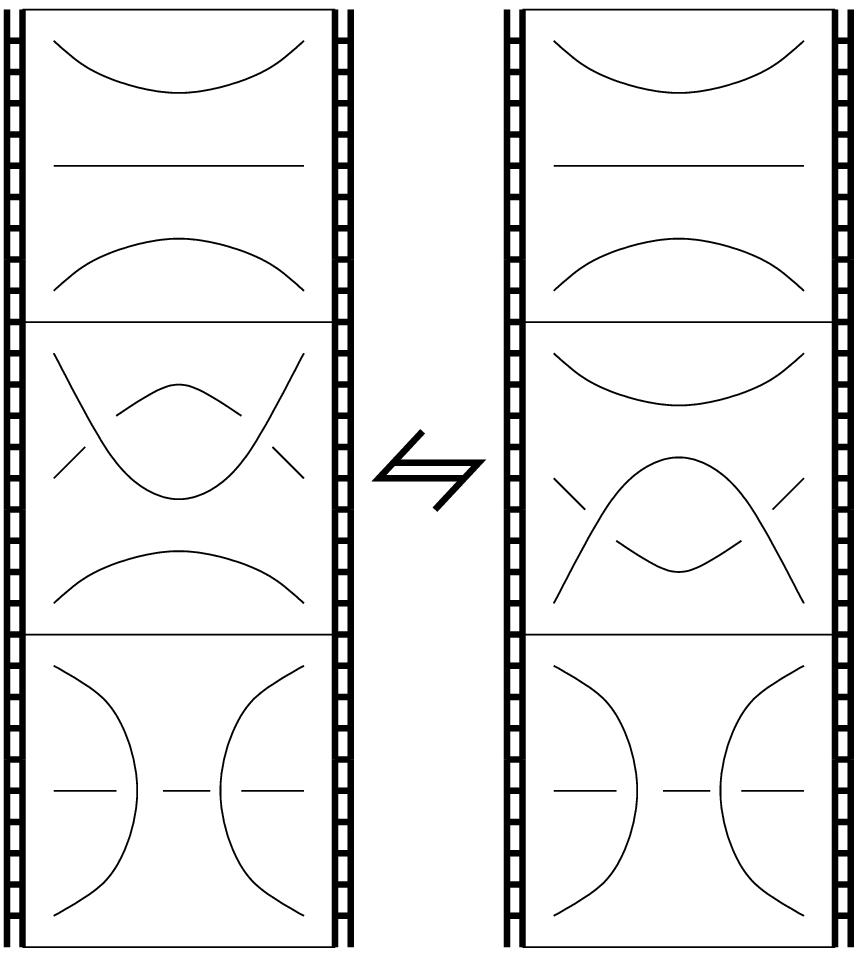}$}}
\item For $\MM_{15}$: Quite nicely, going down the two sides
of $\MM_{15}$ we get the two morphisms $\Psi_L$ and $\Psi_R$ of
\figref{fig:Reid3Proof}, and these two are the same. Going up
$\MM_{15}$ we get, in a similar manner, the two sides of the other
variant of the $R3$ move, as at the end of the invariance under $R3$
proof on just above Lemma~\ref{lem:ConeHomotopy}.

\end{myitemize}

\vskip 2mm

This concludes the proof of Theorems~\ref{thm:Main}
and~\ref{thm:KhIsFunctor}. \qed

\section{More on $\protect\Cobl^3$} \label{sec:MoreOnCobl}

In Section~\ref{sec:Homology} we've seen how a functor
$\Cobl^3\to\ZMod$ can take our theory (which now includes
Theorems~\ref{thm:KhIsFunctor} and~\ref{thm:Main} as well) into
something more computable. Thus we seek to construct many such
functors.  We start, right below, with a systematic construction of
such ``tautological'' functors.  Then in
Sections~\ref{subsec:OriginalKhovanov}--\ref{subsec:4TuZ2} we will
discuss several instances of the tautological construction, leading
back to the original Khovanov theory
(Section~\ref{subsec:OriginalKhovanov}), to the
Lee~\cite{Lee:AlternatingLinks} variant of the original Khovanov theory
(Section~\ref{subsec:Lee}) and to a new variant defined only over the two
element field $\bbF_2$ (Section~\ref{subsec:4TuZ2}).

\begin{definition} Let $B$ be a set of points in $S^1$ and let $\calO$
be an object of $\Cobl^3(B)$ (ie, a smoothing with boundary $B$;
often if $B=\emptyset$ we will choose $\calO$ to be the empty
smoothing). The tautological functor $\calF_\calO\co \Cob_l^3\to\ZMod$ is
defined on objects by $\calF_\calO(\calO'):=\Mor(\calO,\calO')$ and on
morphisms by composition on the left. That is, if $F\co \calO'\to\calO''$
is a morphism in $\Cobl^3(B)$ then
$\calF_\calO(F)\co \Mor(\calO,\calO')\to\Mor(\calO,\calO'')$ maps
$G\in\Mor(\calO,\calO')$ to $F\circ G\in\Mor(\calO,\calO'')$.
\end{definition}

At the moment we don't know the answer to the following question.

\begin{problem} Is the tautological construction faithful? Is there more
information in $\Kh$ beyond its composition with tautological functors?
Beyond the homology of its composition with tautological functors?
\end{problem}

The groups of morphisms of $\Cobl^3$ appear to be difficult to study.
Hence we will often simplify matters a bit by composing tautological
functors with some extra functors that forget some information. Examples
follow below.

\subsection{Cutting necks and the original Khovanov homology theory}
\label{subsec:OriginalKhovanov}

As a first example we take $B=\emptyset$ and $\calO=\emptyset$, we
forget all 2--torsion by tensoring with some ground ring in which
$2^{-1}$ exists (eg, $\Ztwo:=\bbZ[1/2]$) and we mod out by all surfaces
with genus greater than 1:
\[
  \calF_1(\calO'):=\Ztwo\otimes_\bbZ\Mor(\emptyset,\calO')/((g>1)=0).
\]

Taking $C$ to be a disjoint union of two twice-punctured disks in the
specification of the $\FourTu$ relation in Section~\ref{subsubsec:STFourTu}
we get the neck cutting relation

\begin{equation} \label{eq:CuttingNecks}
  2\begin{array}{c}\includegraphics[height=10mm]{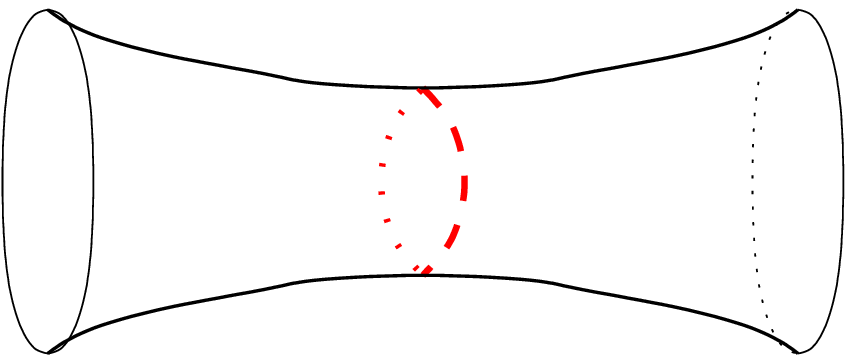}\end{array}
  =\begin{array}{c}\includegraphics[height=10mm]{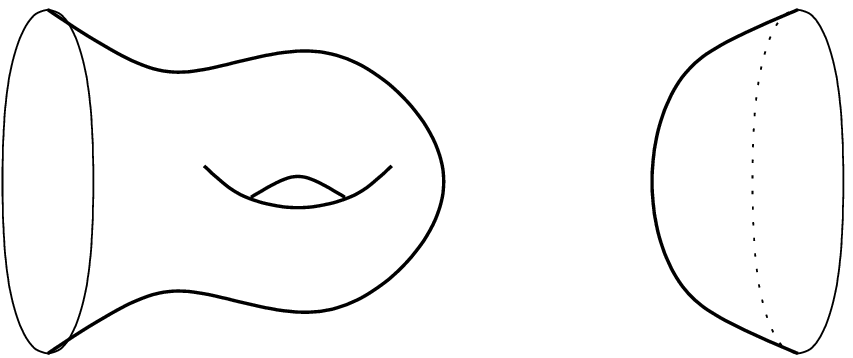}\end{array}
  +\begin{array}{c}\includegraphics[height=10mm]{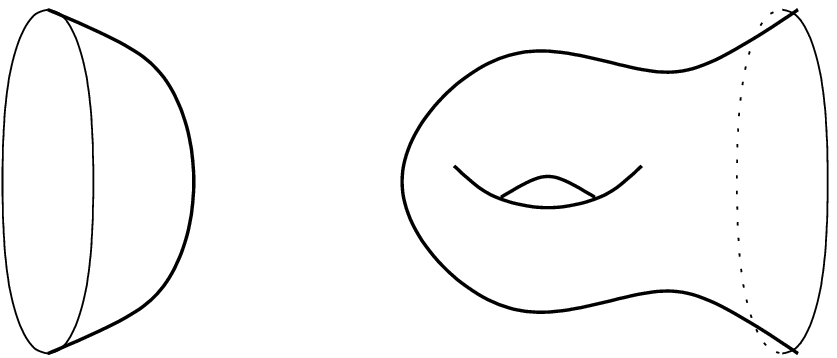}\end{array}
\end{equation}
\vskip 1mm

If $2$ is invertible, the neck cutting relation means that we can ``cut
open'' any tube inside a cobordism (replacing it by handles that are
localized to one side of the tube).  Repeatedly cutting tubes in this
manner we see that $\Ztwo\otimes_\bbZ\Mor(\emptyset,\calO')$ is
generated by cobordisms in which every connected component touchs at
most one boundary curve. Further reducing using the $S$, $T$ and
$((g>1)=0)$ relations we get to cobordisms in which every connected
component has precisely one boundary curve and is either of genus $0$
or of genus $1$. So if $\calO'$ is made of $k$ curves then
$\calF_1(\calO')=V^{\otimes k}$ where $V$ is the $\Ztwo$--module
generated by $v_+:=\eps{4mm}{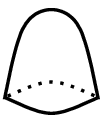}$ and by $v_-:=\frac12\eps{5mm}{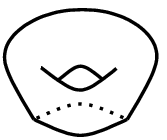}$.

\begin{exercise} Verify that with this basis for $V$ the (reduced)
tautological functor $\calF_1$ becomes the same as the functor $\calF$
of Definition~\ref{def:Vdef}, and hence once again we reproduce the
original Khovanov homology theory.
\end{exercise}

\vskip -4.5mm
\parpic(3.3in,0.9in)[r]{\raisebox{-18mm}{$
  \begin{array}{c|c|c}
    \rule[-6pt]{0pt}{18pt} \text{glue} &
      \left|\eps{4mm}{vp}\!\right\rangle &
      \left|\frac12\eps{5mm}{vm}\!\right\rangle \\
    \hline \rule[-9pt]{0pt}{24pt}
      \left\langle\frac12\eps{5mm}{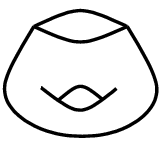}\!\right| &
      \frac12\left\langle\epsg{5mm}{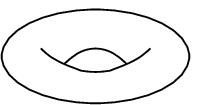}\!\right\rangle = 1 &
      \frac14\left\langle\eps{5mm}{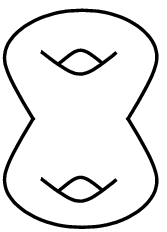}\!\right\rangle = 0 \\
    \hline \rule[-5pt]{0pt}{18pt}
      \left\langle\eps{4mm}{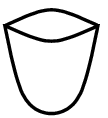}\!\right| &
      \left\langle\eps{5mm}{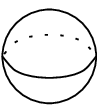}\!\right\rangle = 0 &
      \frac12\left\langle\epsg{5mm}{g1}\!\right\rangle = 1
  \end{array}
  \quad
  \begin{array}{c}
    \epsg{9mm}{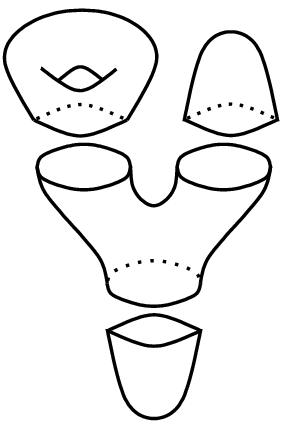} \\
    \text{\small cactus}
  \end{array}
$}}
\begin{hint} Using the bra-ket notation of quantum mechanics, the
bras dual to the kets $v_+=\left|\eps{4mm}{vp}\!\right\rangle$ and
$v_-=\left|\frac12\eps{5mm}{vm}\!\right\rangle$ are
$\left\langle\frac12\eps{5mm}{dvm}\!\right|$ and
$\left\langle\eps{4mm}{dvp}\!\right|$ respectively. Thus, for example,
the coefficient of $v_-$ within the product $v_-v_+$ is $\left\langle
\eps{4mm}{dvp}\! \left| \eps{8mm}{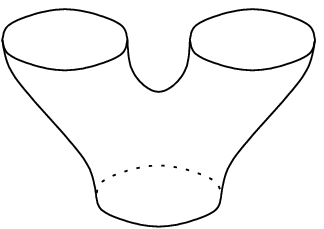}\! \right|
\frac12\eps{5mm}{vm}\!\eps{4mm}{vp}\! \right\rangle =
\frac12\left\langle\text{cactus}\right\rangle =
\frac12\left\langle\epsg{5mm}{g1}\!\right\rangle = 1$.
\end{hint}

\subsection{Genus $3$ and Lee's theory} \label{subsec:Lee}

At first glance it may appear that the relation $((g>1)=0)$ was
unnecessary in the above discussion --- every high genus surface
contains several ``necks'' and we can cut those
using~\eqref{eq:CuttingNecks} to get lower genus surfaces. This clearly
works if the genus $g$ is high enough to start with ($g\geq 4$ is
enough). Cutting the obvious neck in the genus $2$ surface
$\epsg{7mm}{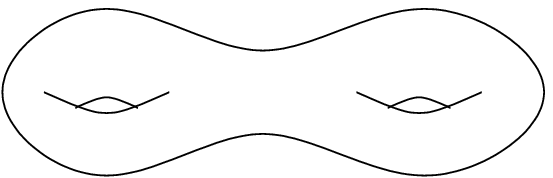}$ and reducing tori using the $T$ relation we find that
$\epsg{7mm}{g2}=0$ automatically. But the genus $3$ surface with no boundary
$\epsg{9mm}{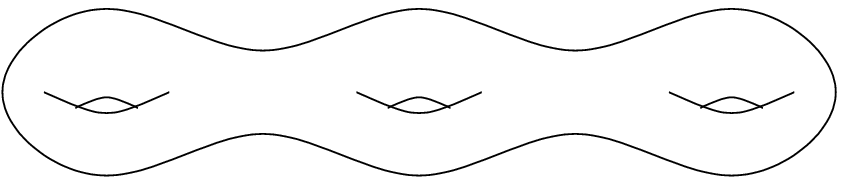}$ cannot be reduced any further.

Thus setting
\[
  \calF_2(\calO')
  := \Ztwo\otimes_\bbZ\Mor(\emptyset,\calO')/(\epsg{9mm}{g3}=8)
\]
we find that as a $\Ztwo$--module $\calF_2(\calO')$ is as in
the previous example and as in Definition~\ref{def:Vdef} (except that
the grading is lost), but $\Delta$ and $m$ are modified as follows:
\[
  \Delta_2: \begin{cases}
    v_+ \mapsto v_+\otimes v_- + v_-\otimes v_+ &\\
    v_- \mapsto v_-\otimes v_- + v_+\otimes v_+ &
  \end{cases}
  \qquad
  m_2: \begin{cases}
    v_+\otimes v_-\mapsto v_- &
    v_+\otimes v_+\mapsto v_+ \\
    v_-\otimes v_+\mapsto v_- &
    v_-\otimes v_-\mapsto v_+.
  \end{cases}
\]

This is precisely the Lee~\cite{Lee:AlternatingLinks} variant of the
original Khovanov theory (used by Rasmussen~\cite{Rasmussen:SliceGenus}
to give a purely combinatorial proof of the Milnor conjecture).

\subsection{Characteristic 2 and a new homology theory}
\label{subsec:4TuZ2}

\begin{figure}[t]\anchor{fig:H4Tu}
\[ 
  \begin{array}{c}\includegraphics[height=15mm]{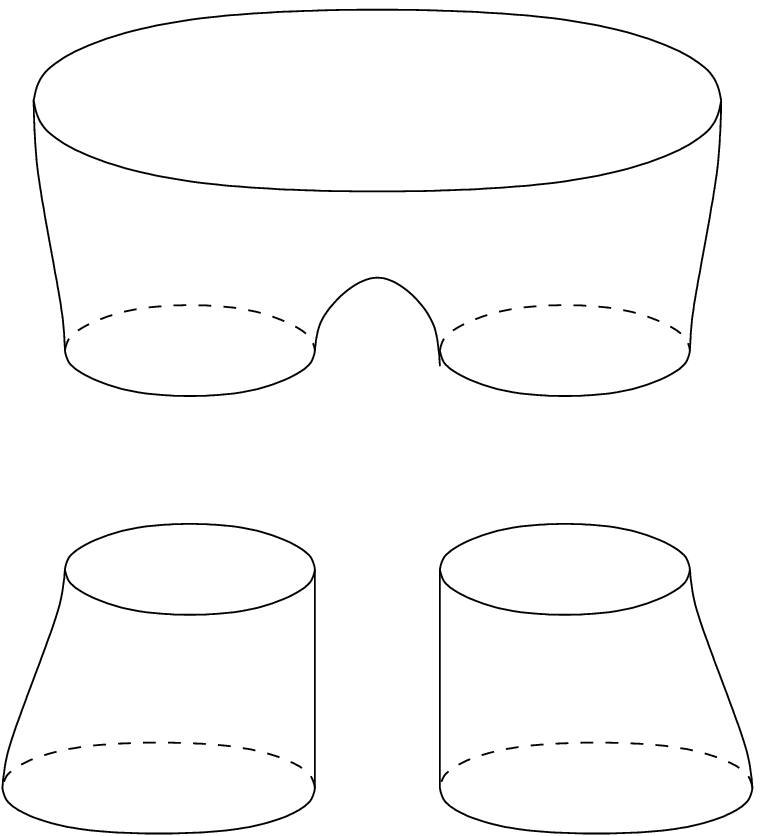}\end{array}
  \qquad\qquad
  \begin{array}{c}\includegraphics[height=15mm]{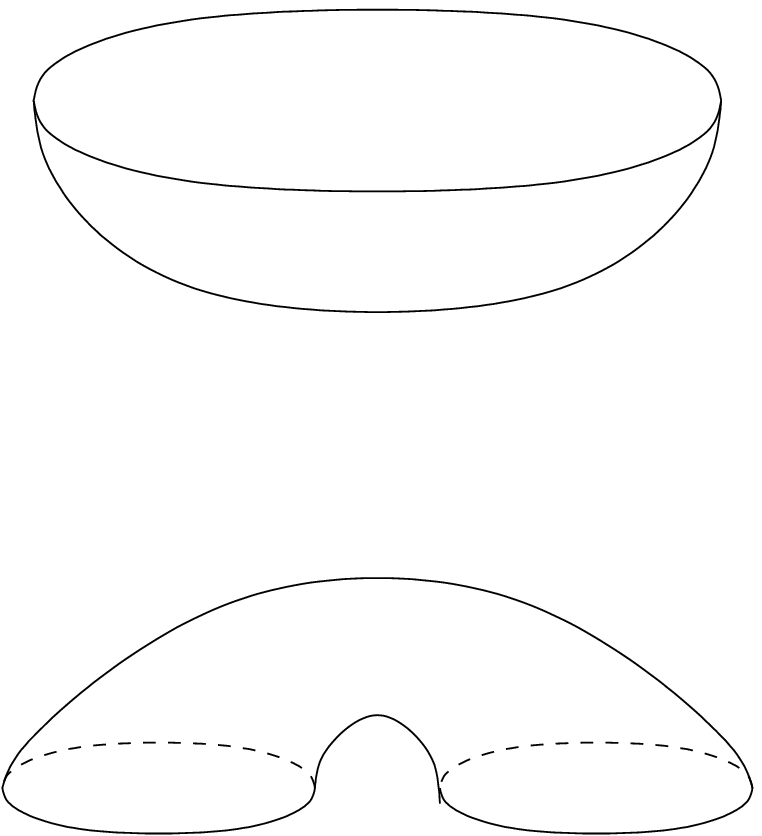}\end{array}
  =\begin{array}{c}\includegraphics[height=15mm]{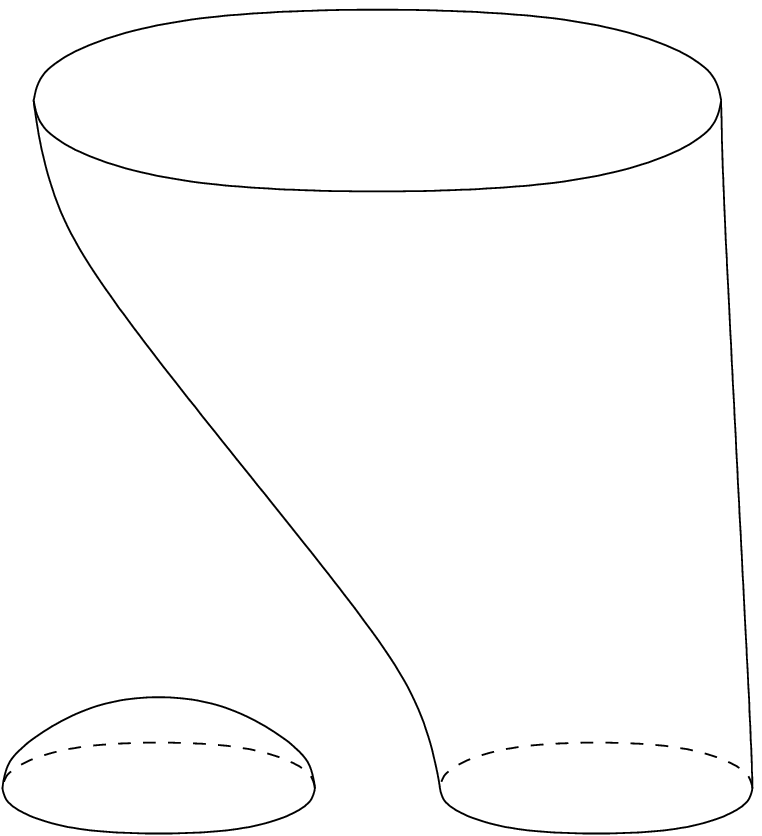}\end{array}
  +\begin{array}{c}\includegraphics[height=15mm]{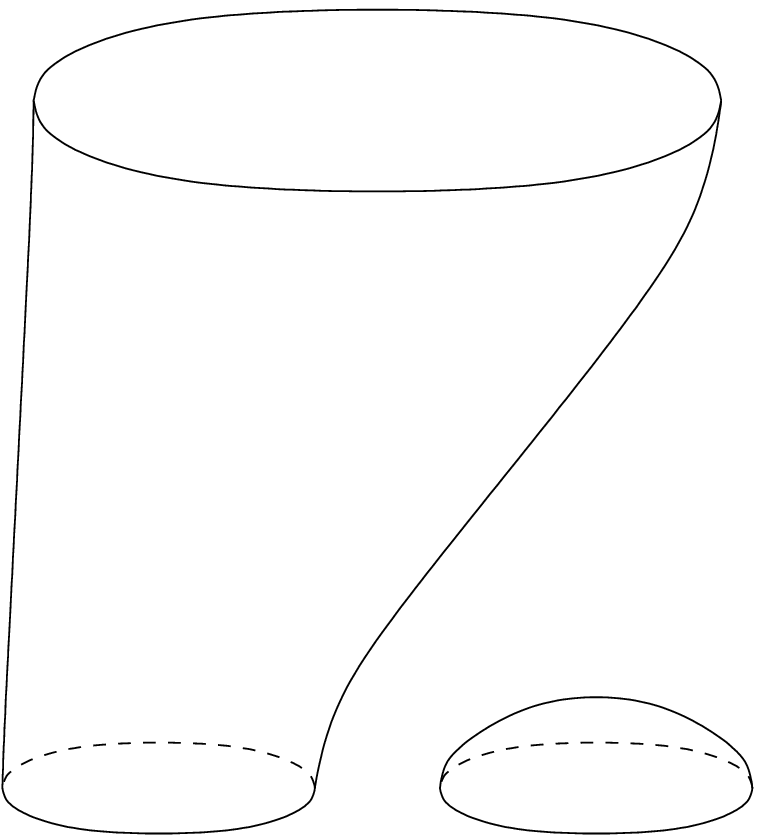}\end{array}
  +H\begin{array}{c}\includegraphics[height=15mm]{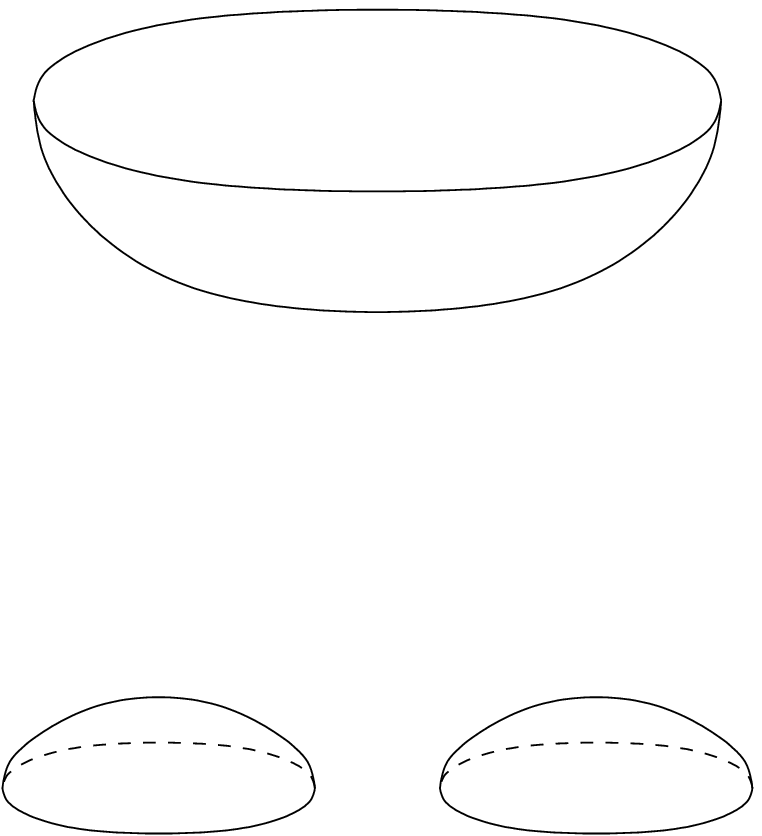}\end{array}
\]
\caption{
  A local picture and the corresponding $\FourTu$ relation (over 
  $\bbF_2[H]$, so signs can be disregarded and handles replaced by $H$'s).
} \label{fig:H4Tu}
\end{figure}
The other extreme is to tensor multiply with $\bbF_2$ (ie, to take
all linear combinations with coefficients in $\bbF_2$). In this case it
is convenient to take the ``target object'' $\calO$ to be a single
boundariless cycle $\bigcirc$ and set
\[
  \calF_3(\calO') := \bbF_2\otimes_\bbZ\Mor(\bigcirc,\calO').
\]
\parpic[r]{$
  \begin{array}{c}\includegraphics[height=10mm]{figs/CNL.eps}\end{array}
  =
  \begin{array}{c}\includegraphics[height=10mm]{figs/CNR.eps}\end{array}
$}
Over $\bbF_2$ the neck drops off the neck cutting
relation~\eqref{eq:CuttingNecks} and what remains is a relation (shown
on the right) saying that handles can be moved from one component to
another. We introduce a new variable $H$ with $\deg H=-2$ and with the
relation $H\epsg{4mm}{vp}=\epsg{5mm}{vm}$, so $H$ stands for ``a handle
inserted somewhere (anywhere) into a cobordism''. As modules over
$\bbF_2[H]$ the morphisms of $\bbF_2\otimes\Cobl^3$ are generated by
cobordisms with no handles, ie, by punctured spheres. (An exception
needs to be made for cobordisms with no boundary at all, in
$\Cobl^3(\emptyset)$. Our morphisms all have a source $\calO=\bigcirc$
and hence they always have at least one boundary component so this
exception is irrelevant in what follows.)

%\parpic[r]{$\eps{60mm}{4TuZ2Generator}$}
\vskip 2mm
\begin{floatingfigure}[r]{63mm}\anchor{fig:4TuZ2Generator}
\centering{$\eps{60mm}{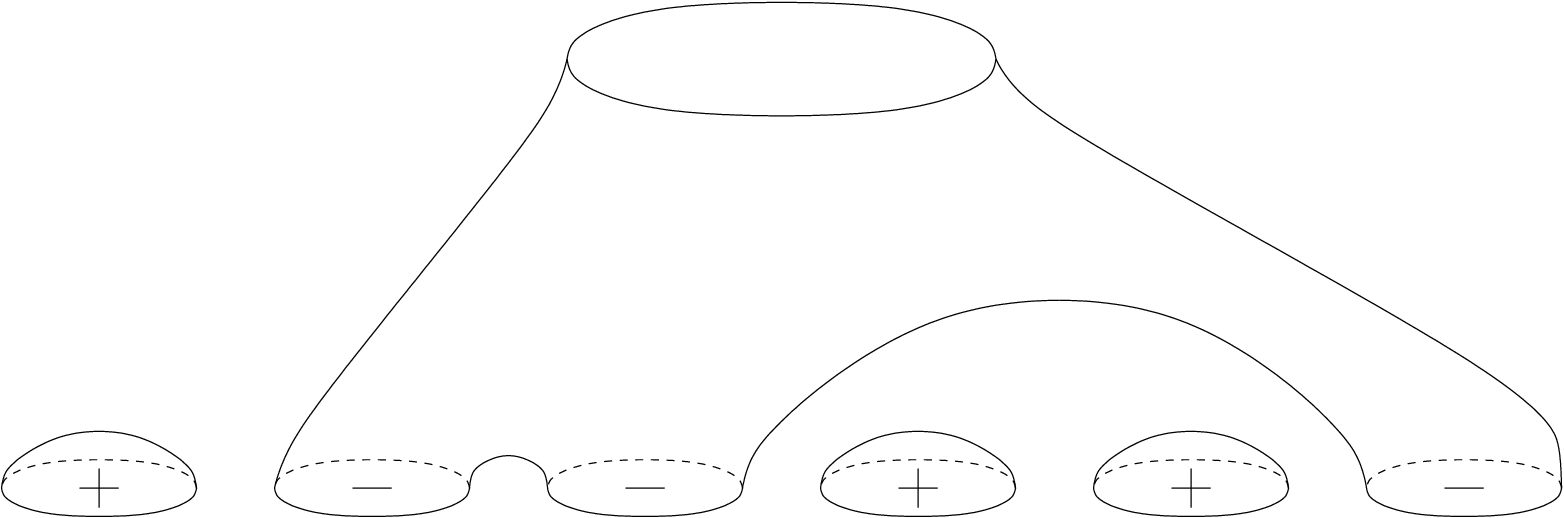}$}
\vskip -2mm
\caption{$v_{+--++-}\in\calF_3(\bigcirc^6)$} \label{fig:4TuZ2Generator}
\vskip -1mm
\end{floatingfigure}
Let $\bigcirc^k=\bigcirc\bigcirc\cdots\bigcirc$ ($k$ circles) be a
typical object of $\Cobl^3(\emptyset)$ and consider
$\calF_3(\bigcirc^k)= \bbF_2\otimes_\bbZ\Mor(\bigcirc,\bigcirc^k)$.
Over $\bbF_2[H]$ this module is generated by surfaces whose
components are punctured spheres with an overall number of $k+1$
punctures, $k$ of which corresponding to the $k$ circles in the
``target object'' $\bigcirc^k$ and a special puncture corresponding to
the ``source object'' $\bigcirc$. The $\FourTu$ relation of
\figref{fig:H4Tu} can be used to dissolve any component that has
more than one puncture and that does not contain the special puncture
into components that either contain just one puncture or contain the
special puncture. Hence $\calF_3(\bigcirc^k)$ is generated by
cobordisms such as the one in \figref{fig:4TuZ2Generator}, in which
there is a special spherical component whose boundary contains the
source object (and possibly several other circles from the target
object) and possibly several disjoint disks that cap the remaining
parts of the target object. As in \figref{fig:4TuZ2Generator},
these generators are in a natural bijective correspondence with the
elements of $V^{\otimes k}$ where $V$ (as before) is generated by two
elements $v_\pm$ with $\deg v_\pm=\pm 1$ (though this time $v_+$ and $v_-$
cannot be identified with the disk and half a punctured torus as in
Sections~\ref{subsec:OriginalKhovanov} and~\ref{subsec:Lee}).

\begin{exercise} Check that with this identification of
$\calF_3(\bigcirc^k)$ with $V^{\otimes k}$, the generating morphisms of
$\Cobl^3(\emptyset)$ map to
\begin{alignat*}{1}
  \Delta_3: & \begin{cases}
    v_+ \mapsto v_+\otimes v_- + v_-\otimes v_+ +Hv_+\otimes v_+ & \\
    v_- \mapsto v_-\otimes v_- &
  \end{cases} \\
  m_3: & \begin{cases}
    v_+\otimes v_-\mapsto v_- &
    v_+\otimes v_+\mapsto v_+ \\
    v_-\otimes v_+\mapsto v_- &
    v_-\otimes v_-\mapsto Hv_-.
  \end{cases}
\end{alignat*}
\end{exercise}

\begin{figure}[ht!]\anchor{fig:10n136Computed}
\begin{center} \scriptsize
  \def\merge{\multicolumn{2}{c|}{}}
  \def\putknot{
    \multicolumn{2}{c|}{\smash{\raisebox{1.5mm}{$\eps{1.1cm}{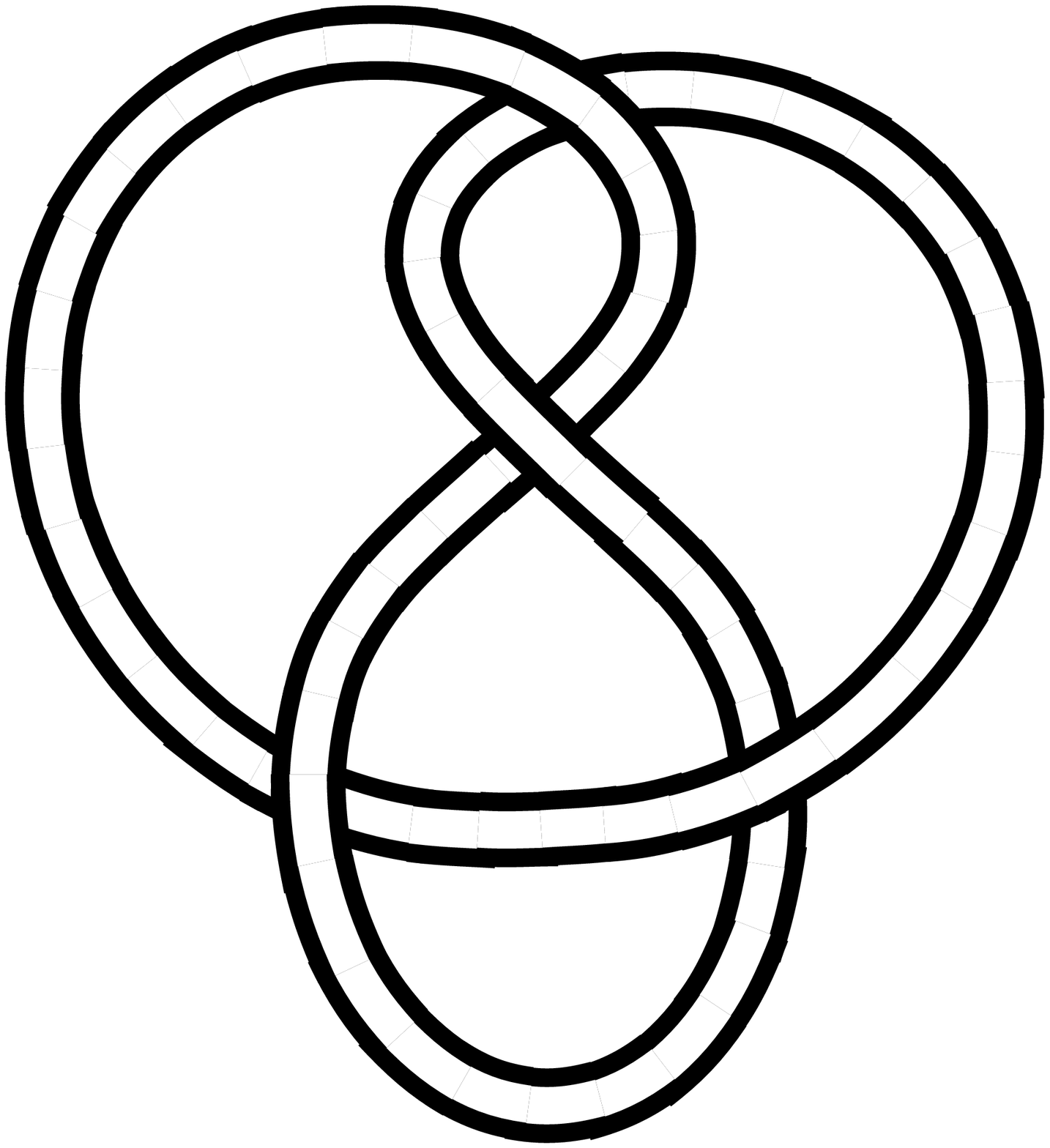}$}}}
  }
  \begin{tabular}{|c|c|c|c|c|c|} \hline
    \backslashbox{$j$}{$r$}
          & $-2$& $-1$&  0  &  1  &  2  \\ \hline
      5   &     &     &     &     &1,1,1\\ \hline
      3   &     &     &     &0,1,0&0,1,1\\ \hline
      1   &     &     &1,1,1&1,1,0&     \\ \hline
     $-1$ &     &1,1,1&1,1,2&  \merge   \\ \cline{1-4}
     $-3$ &0,1,0&0,1,1&0,0,2&  \merge   \\ \cline{1-4}
     $-5$ &1,1,0&     &0,0,2& \putknot  \\ \cline{1-4}
     $<-5$&     &     &0,0,2&  \merge   \\ \hline
  \end{tabular}
  \par \vskip 3mm
  \def\putknot{\multicolumn{2}{c|}{$\smash{\eps{1.35cm}{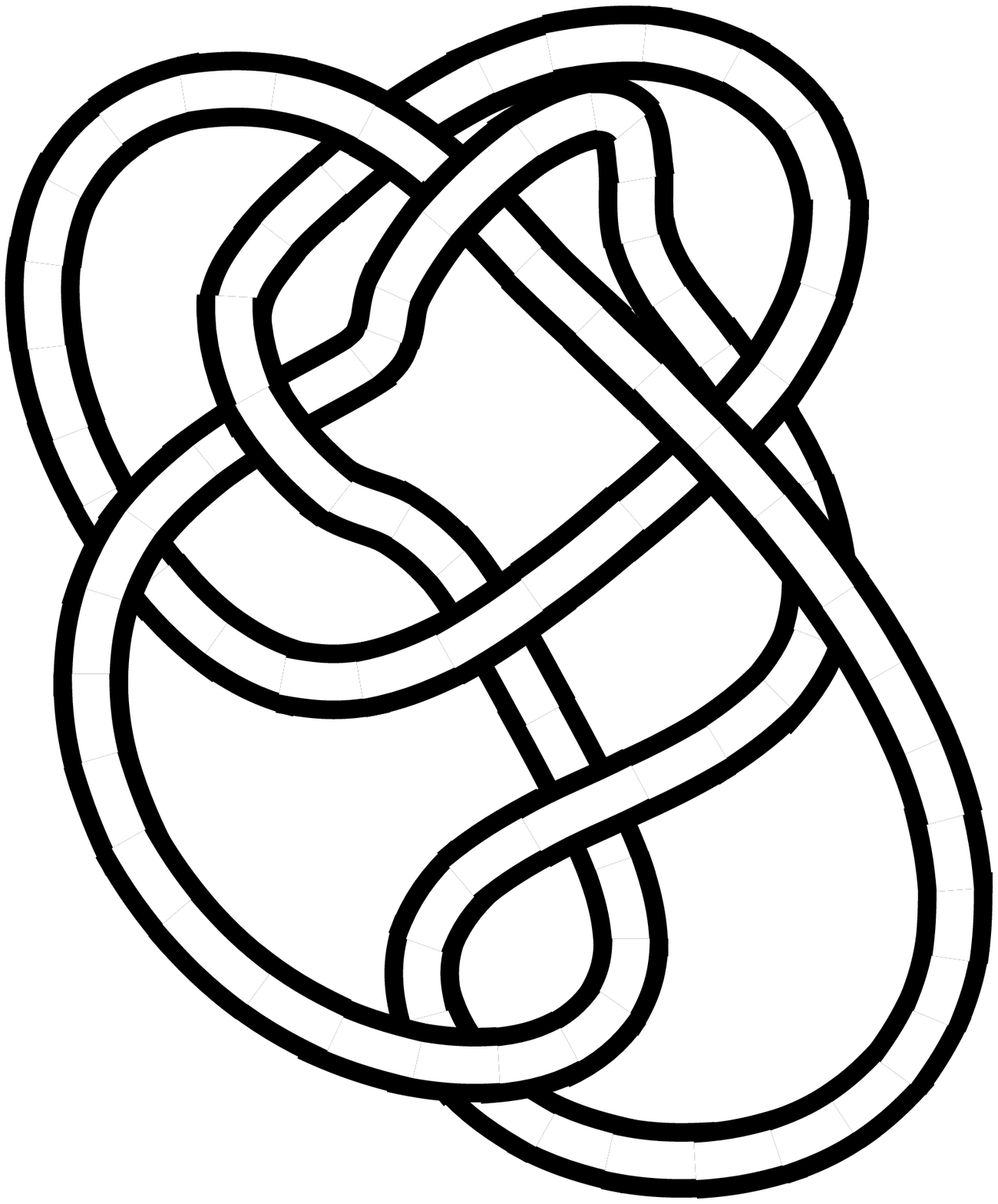}}$}}
  \begin{tabular}{|c|c|c|c|c|c|c|c|c|} \hline
    \backslashbox{$j$}{$r$}
          & $-4$& $-3$& $-2$& $-1$&  0  &  1  &  2  &  3  \\ \hline
      9   &     &     &     &     &     &     &     &1,1,1\\ \hline
      7   &     &     &     &     &     &     &1,2,1&0,1,1\\ \hline
      5   &     &     &     &     &     &1,2,1&1,2,1&     \\ \hline
      3   &     &     &     &     &2,3,2&1,2,1&     &     \\ \hline
      1   &     &     &     &1,3,1&2,4,3&     &     &     \\ \hline
     $-1$ &     &     &1,2,1&2,3,1&1,1,2&     &  \merge   \\ \cline{1-7}
     $-3$ &     &1,2,1&1,2,1&     &0,0,2&     &  \merge   \\ \cline{1-7}
     $-5$ &0,1,0&1,2,1&     &     &0,0,2&     & \putknot  \\ \cline{1-7}
     $-7$ &1,1,0&     &     &     &0,0,2&     &  \merge   \\ \cline{1-7}
     $<-7$&     &     &     &     &0,0,2&     &  \merge   \\ \hline
  \end{tabular}
\end{center}
\caption{
  The knots $4_1$ and  $10_{136}$ (as drawn
  by~\cite{Scharein:KnotPlot}) and their Betti numbers in a tabular
  form.  Each box contains the three numbers $b^\bbQ_{rj}$,
  $b^{\bbF_2}_{rj}$ and $b^3_{rj}$ in order. Empty boxes mean that all
  three Betti numbers are zero.
} \label{fig:10n136Computed}
\end{figure}

We know nothing in general about the homology of $\calF_3\Kh(L)$
associated with a knot/link $L$ and/or its relationship with the
original Khovanov homology $H(\calF\Kh(L))$. Let us describe here the
results of some sporadic computations that we have performed. We took
$H=1$ in the above formulas (losing the grading, of course) and
obtained a filtered theory, where $\calG_{\geq j}\calF^1_3\Kh(L)$
denotes the subcomplex of $\calF^1_3\Kh(L)$ made of chains of degrees
greater than or equal to $j$. We then computed the Betti numbers
$b^3_{rj}(L):=\dim_{\bbF_2}H^r(\calG_{\geq j}\calF^1_3\Kh(L))$ and
compared them with the Betti numbers of the original Khovanov homology
over $\bbQ$ and over $\bbF_2$ (ie, with $b^\bbQ_{rj}(L):=\dim_\bbQ
H(\bbQ\otimes\calF\Kh(L))$ and with $b^{\bbF_2}_{rj}(L):=\dim_{\bbF_2}
H(\bbF_2\otimes\calF\Kh(L))$), for several knots and links with up to
$10$ crossings. The results of these computations are best displayed as
two dimensional arrays of numbers, as in
\figref{fig:10n136Computed}.

\parpic[r]{\scriptsize\begin{tabular}{ccc}
  $\begin{array}{|c|c|}\hline&1\\\hline&\\\hline1&\\\hline\end{array}$ &
  $\begin{array}{|c|c|}\hline&1\\\hline1&1\\\hline1&\\\hline\end{array}$ &
  $\begin{array}{|c|c|}\hline&1\\\hline&1\\\hline&\\\hline\end{array}$ \\
  knight & tetris & shifted \\
  move   & piece  & pawn \\
\end{tabular}}
Readers familiar with~\cite{Khovanov:Patterns, Shumakovitch:Torsion,
Bar-Natan:Categorification} will notice that every ``knight move'' in
the array describing $b^\bbQ_{rj}$ gets replaced by a ``tetris piece''
in the array for $b^{\bbF_2}_{rj}$ and by a ``shifted pawn move'' in
the array for $b^3_{rj}$, while the special pawn in column $0$ of
$b^\bbQ_{rj}$ stays put for $b^{\bbF_2}_{rj}$ and gets replaced by a
column with entries $(1,2,2,\ldots)$ in the array for $b^3_{rj}$. This
observation repeats for all the knots that we have tested but we don't
know a general explanation (though see Conjecture~\ref{conj:Structure}
below).

\section{Trace groups, Euler characteristics and skein modules}
\label{sec:Euler}

We wish to push the relationship between $\Kh$ and the Jones
polynomial a little further\footnote{
This section was inspired by a talk A\,S~Sikora gave at the George
Washington University on May 2004 and
by~\cite{AsaedaPrzytyckiSikora:SkeinModule}. The results
in~\cite{AsaedaPrzytyckiSikora:SkeinModule} partially overlap ours.}.
Both $\Kh$ and $\hatJ$ can be computed
locally, first for small tangles and then, by composing small tangles,
the computations can be carried out for bigger tangles and whole
links.  But the relationship with the Jones polynomial discussed in
Section~\ref{sec:Homology} only works at the level of whole links. In
this section we will see how an appropriate notion of ``Euler
characteristic'' intertwines $\Kh$ and $\hatJ$ while respecting the
planar algebra structure (so progressive computations of $\Kh$ and
$\hatJ$ ``run in tandem'').

\subsection{Traces} \label{subsec:Traces}

We first have to generalize the notions of Euler characteristic, dimension
and trace to complexes over arbitrary pre-additive categories.

\begin{definition} A ``trace'' on a category $\calC$ with values in an
Abelian group $A$ is an additive map
$\tau\co \bigoplus_{\calO\in\Obj(\calC)}\Mor(\calO,\calO)\to A$ defined on
all the endomorphisms of objects in $\calC$ and with values in $A$,
which satisfies the ``trace relation'' $\tau(FG)=\tau(GF)$ whenever $F$
and $G$ are morphisms so that both $FG$ and $GF$ are endomorphisms
(ie, whenever $F\co \calO_1\to\calO_2$ and $G\co \calO_2\to\calO_1$ for
some $\calO_{1,2}\in\Obj(\calC)$).
\end{definition}

This definition is analogous to the standard definition of the trace of
a matrix, which is defined only for square matrices but satisfies the
trace relation $\tr(FG)=\tr(GF)$ even for non-square $F$ and $G$,
provided both $FG$ and $GF$ are square.

\begin{exercise} If $\tau$ is a trace on $\calC$ and if $(F_{ij})$ is an
endomorphism of some object in $\Mat(\calC)$, set
$\tau((F_{ij})):=\sum_i\tau(F_{ii})$ and show that the newly defined $\tau$
is a trace on $\Mat(\calC)$.
\end{exercise}

Assuming the presence of some fixed trace
$\tau$ on a pre-additive category $\calC$, we can now proceed to define
dimensions, Euler characteristics and Lefschetz numbers:

\begin{definition} The dimension (more precisely, the
``$\tau$--dimension'') of an object $\calO$ of $\calC$ or in
$\Mat(\calC)$ is the trace of the identity:
$\dim_\tau\calO:=\tau(I_\calO)$. If $(\Omega^r)$ is a complex in
$\Kom(\calC)$ or in $\Kom(\Mat(\calC))$ we define its Euler
characteristic (more precisely, its ``$\tau$--Euler characteristic'') to
be $$\chi_\tau((\Omega^r)):=\sum_r(-1)^r\dim_\tau\Omega^r.$$ Finally, if
$F=(F^r)$ is an endomorphism of $\Omega$ we define its Lefschetz number
(or ``$\tau$--Lefschetz number'') to be
$\tau(F):=\sum_r(-1)^r\tau(F^r)$. (All these quantities are members of
the Abelian group $A$ which may or may not be our underlying group of
scalars). Clearly, $\chi_\tau(\Omega)=\tau(I_\Omega)$.
\end{definition}

We claim that Lefschetz number and Euler characteristics are homotopy
invariant:

\begin{proposition}$\phantom{99}$
\begin{enumerate}
\item If $\Omega$ is a complex and $F$ and $G$ are homotopic 
endomorphisms $F,G\co \Omega\to\Omega$ then $\tau(F)=\tau(G)$.
\item If the complexes $\Omega_a$ and $\Omega_b$ are homotopy equivalent
then $\chi_\tau(\Omega_a)=\chi_\tau(\Omega_b)$.
\end{enumerate}
\end{proposition}

\begin{proof}$\phantom{99}$
\begin{enumerate}
\item Let $h$ be so that $F-G=hd+dh$. Then, using the trace relation for
the last equality,
\begin{eqnarray*}
  \tau(F)-\tau(G)
  & = \sum_r(-1)^r\tau(F^r-G^r)
  = \sum_r(-1)^r\tau(h^{r+1}d^r+d^{r-1}h^r) \\
  & = \sum_r(-1)^r\tau(h^{r+1}d^r-d^rh^{r+1}) = 0.
\end{eqnarray*}
\item If $F\co \Omega_a\to\Omega_b$ and $G\co \Omega_b\to\Omega_a$ induce the
homotopy equivalence then $GF\sim I_{\Omega_a}$ and $FG\sim I_{\Omega_b}$
and so
\[
  \chi_\tau(\Omega_a)
  = \tau(I_{\Omega_a})
  = \tau(GF)
  = \tau(FG)
  = \tau(I_{\Omega_b})
  = \chi_\tau(\Omega_b)
\]
\end{enumerate}
(using the trace relation once more, for the middle equality)
\end{proof}

\subsection{The trace group and the universal trace}
\label{subsec:TraceGroups} Given a pre-additive category $\calC$ it is
interesting to find all traces defined on it. Quite clearly, they all
factor through the ``universal trace'' defined below:

\begin{definition} \label{def:TraceGroup}
The ``trace group'' $\Xi(\calC)$ is\footnote{S\,D~Schack informed me that
this notion is due to F\,W~Lawvere and S\,H~Schanuel. It is also the same
as ``$0$th Hochschild--Mitchell homology of a category with coefficients in
itself'', but it is not the same as the ``Grothendieck group'' of an Abelian
category.}
\[ \Xi(\calC) :=
  \left.\bigoplus_{\calO\in\Obj(\calC)}\Mor(\calO,\calO)\right/
  \parbox{2.8in}{
    the trace relation: $FG=GF$ whenever $F\co \calO_1\to\calO_2$ and
    $G\co \calO_2\to\calO_1$.
  }
\]
The map which takes any endomorphism in $\calC$ to itself as a member of
$\Xi(\calC)$ is denoted $\tau_\star$ and called ``the universal trace of
$\calC$''. (It is, of course, a trace). We denote dimensions and Euler
characteristics defined using $\tau_\star$ (hence valued in $\Xi(\calC)$)
by $\dim_\star$ and $\chi_\star$.
\end{definition}

\begin{example} \label{exmp:MatrixTrace} (Told by D~Thurston; see
also~\cite{Stallings:CenterlessGroups}) Let
$\Mat(\bullet)$ be the category whose objects are non-negative integers
and whose morphisms are rectangular integer matrices of appropriate
dimensions:  $\Mor(n,m)=\{m\times n\text{ matrices}\}$. Then the trace
group $\Xi(\Mat(\bullet))$ is generated by square matrices. The trace
relation
\[
  \begin{pmatrix}0&f\\0&0\end{pmatrix}
  = \begin{pmatrix}0&f\\0&0\end{pmatrix}\begin{pmatrix}0&0\\0&1\end{pmatrix}
  = \begin{pmatrix}0&0\\0&1\end{pmatrix}\begin{pmatrix}0&f\\0&0\end{pmatrix}
  = \begin{pmatrix}0&0\\0&0\end{pmatrix}
\]
allows us to cancel off-diagonal entries (perhaps after extending the
relation by adding zero rows and columns). Diagonal matrices can be
written as sums of diagonal matrices that have just one non-zero entry on
the diagonal, and then empty row/column pairs can be removed using the
trace relation
\[
  \begin{pmatrix}g&0\\0&0\end{pmatrix}
  = \begin{pmatrix}g\\0\end{pmatrix}\begin{pmatrix}1&0\end{pmatrix}
  = \begin{pmatrix}1&0\end{pmatrix}\begin{pmatrix}g\\0\end{pmatrix}
  = \begin{pmatrix}g\end{pmatrix}.
\]
Thus $\Xi(\Mat(\bullet))$ is generated by the $1\times 1$ matrix
$\begin{pmatrix}1\end{pmatrix}$, and so the ordinary matrix trace is up
to scalars the unique trace on $\Mat(\bullet)$ and $\Xi(\Mat(\bullet))$ is
isomorphic to $\bbZ$.
\end{example}

\begin{example}
It is likewise easy to show that the trace group of the category of
finite dimensional vector spaces and linear maps over some ground field
$\bbF$ is $\bbF$ itself (generated by the identity on a one-dimensional
vector space) and that (up to scalars) the unique trace on finite
dimensional vector spaces is the ordinary trace.
\end{example}

\begin{exercise} Refine the argument in Example~\ref{exmp:MatrixTrace}
to show that $\Xi(\Mat(\calC))$ $=\Xi(\calC)$ for any pre-additive category
$\calC$.
\end{exercise}

\begin{example} Let $\Vect_0$ be the category of finite dimensional
graded vector spaces (over a field $\bbF$) with degree 0 morphisms
between them. Then, loosely speaking, the different homogeneous
components don't interact and so $\Xi(\Vect_0)$ has one generator in
each degree. The generator in degree $m$ can be taken to be the
identity on a one-dimensional graded vector space whose only non-zero
homogeneous component is in degree $m$. Denoting this generator by
$q^m$ we can identify $\Xi(\Vect_0)$ with the ring of Laurent
polynomials $\bbF[q,q^{-1}]$ (where $q$ is a formal variable which is
best thought of as carrying degree $1$). Thus in this case $\dim_\star$
and $\chi_\star$ are the $q$--dimension and $q$--Euler characteristic of
the theory of finite dimensional graded vector spaces (denoted $\qdim$
and $\chi_q$ in~\cite{Bar-Natan:Categorification}).
\end{example}

\subsection{The trace groups of $\Cobz^3$ and skein modules}
\label{subsec:SkeinModules} Theorems~\ref{thm:invariance},
\ref{thm:PlanarAlgebra} and~\ref{thm:graded} all provide us with
invariants of tangles with values in homotopy classes of complexes. In
each case we can find the trace group of the category underlying the
target complexes and then compute universal Euler characteristics, thus
getting potentially simpler ``non-homological'' invariants. Here we
only do it for the most refined of the three invariants, $\Kh$ of
Theorem~\ref{thm:graded}, though to keep things simpler, we also make
the two simplifications of Section~\ref{subsec:OriginalKhovanov} --- we
forget all 2--torsion by tensoring with $\Ztwo$ and we mod out by
all surfaces with genus greater than 1.  As we shall see, this recovers
the Jones polynomial for tangles.

Some preliminary definitions are required.  Let $B_k$ be a collection
of $k$ points in $S^1$. Recall from Section~\ref{sec:grading} that by
allowing artificial degree shifts $\Cobl^3(B_k)$ can be considered as a
graded category and from Theorem~\ref{thm:graded} that in as much as
$\Kh$ is concerned, we can restrict our attention to morphisms in
$\Cobl^3(B_k)$ that, degree shifts considered, are of degree $0$. Let
$\Cobz^3$ be the restricted category, tensored with $\Ztwo$ and
modulo $((g>1)=0)$:  the objects of $\Cobz^3$ are animals of the form
$S\{m\}$ where $S$ is a smoothing and $m$ is an integer indicating a
formal degree shift, and the morphisms from $S_1\{m_1\}$ to
$S_2\{m_2\}$ are $\Ztwo$--linear combinations of cobordisms whose
top is $S_1$ and bottom is $S_2$ and whose degree in the sense of
Definition~\ref{def:degree} is $m_1-m_2$, taken modulo the $S$, $T$,
$\FourTu$ and $((g>1)=0)$ relations. Let $\Xi_k$ be the trace group
$\Xi(\Cobz^3(k))$ of the category $\Cobz^3$. The collection $(\Xi_k)$
inherits a planar algebra structure from the planar algebra structure
of the morphisms of $\Cob^3$ and as the universal Euler characteristic
is homotopy invariant, the first part of the following theorem is
evident:

\begin{theorem} \label{thm:Euler}$\phantom{99}$
\begin{enumerate} 
\item $\chi_\star\circ\Kh$ is an invariant of tangles with values in
$(\Xi_k)$; in fact, $\chi_\star\circ\Kh$ is an oriented planar algebra
morphism $(\calT(s))\to(\Xi(s))$ (where $\Xi(s):=\Xi(|s|)$).
\item $\chi_\star\circ\Kh$ is the Jones polynomial for tangles (proof
follows the proof Proposition~\ref{prop:XiIsS}).
\end{enumerate}
\end{theorem}

Recall from~\cite{Przytycki:SkeinModules} that the Jones polynomial
$\hatJ$ for $k$--ended tangles\footnote{Our slightly non-standard
normalization was chosen to make the statement of
Theorem~\ref{thm:Euler} as simple as possible. For links,
$\hatJ(L)(q)=(q+q^{-1})J(L)(q^2)$, where $J$ is the standard Jones
polynomial, as normalized (say) in~\cite{Kauffman:Bracket}.} can be
defined via the ``skein relations'' $\hatJ\co  \overcrossing\mapsto
q\smoothing-q^2\hsmoothing$ and
$\hatJ\co \undercrossing\mapsto-q^{-2}\hsmoothing+q^{-1}\smoothing$. It
takes values in the ``skein module'' $\calS_k$ for $k$--ended tangles
--- the collection of $\bbZ[q,q^{-1}]$--linear combinations of $k$--ended
crossingless tangles modulo the extra relation $\bigcirc=q+q^{-1}$.

\begin{proposition} \label{prop:XiIsS} The trace group $\Xi_k$
is naturally isomorphic to the skein module $\calS_k$. If $S$ a
$k$--ended crossingless tangle (ie, a smoothing), the isomorphism
$\sigma\co \calS_k\to\Xi_k$ maps $q^mS$ to the identity automorphism of the
object $S\{m\}$ of $\Cobz^3(k)$, regarded as a member of the trace group
$\Xi_k=\Xi(\Cobz^3(k))$.
\end{proposition}

%\parpic[r]{$q^m\eps{0.45in}{sigma3}$}
\begin{proof} $\Xi_k$ is spanned by pairs $(S\{m\},C)$ where $S$ is a
$k$--ended smoothing, $m$ is an integer indicating a formal degree
shift, and $C$ is a degree $0$ cobordism $C\co S\to S$. Such pairs are
taken modulo the trace relation $FG=GF$ of Definition~\ref{def:TraceGroup}.
We can visualize such a pair as a cylinder with $C$ inside, with $S$ on
the top and at the bottom and with extra $q^m$ coefficient placed in front,
(below left).
{\small$$q^m\eps{0.45in}{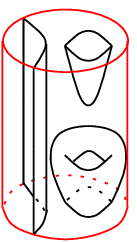}\qquad\qquad\qquad\qquad
  \def\s{\scriptstyle}
  q^{m_1}\eps{0.45in}{sigma3}
    \hspace{-5pt}{\begin{array}{c}\s G\\\ \\\s F\end{array}}
  =q^{m_1+1}\eps{0.45in}{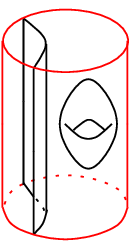}
    \hspace{-5pt}{\begin{array}{c}\ \\\s F\\\s G\end{array}}
$$}%
\iffalse
\parpic[r]{$
  \def\s{\scriptstyle}
  q^{m_1}\eps{0.45in}{sigma3}
    \hspace{-5pt}{\begin{array}{c}\s G\\\ \\\s F\end{array}}
  =q^{m_1+1}\eps{0.45in}{sigma5}
    \hspace{-5pt}{\begin{array}{c}\ \\\s F\\\s G\end{array}}
$}\fi
Let $F\co S_1\{m_1\}\to S_2\{m_2\}$ and $G\co S_2\{m_2\}\to
S_1\{m_1\}$ be arbitrary degree $0$ morphisms in $\Cobz^3(k)$. That
is, $F$ and $G$ are degree $m_1-m_2$ and $m_2-m_1$ %(respectively)
cobordisms within a cylinder (degrees measured as in
Definition~\ref{def:degree}), with tops $S_1$ and $S_2$ %(respectively)
and bottoms $S_2$ and $S_1$ (all respectively), see above right.  The
trace relation states that $C:=GF$, an endomorphism of $S_1\{m_1\}$,
is equal to $FG$, which is an endomorphism of $S_2\{m_2\}$. Setting
$m=m_2-m_1$ and visualizing as before, the trace relation becomes
``cut a degree $m$ piece $G$ off the top of $C$ and reglue it at the
bottom while multiplying by extra factor of $q^m$''. With this
interpretation of the trace relation it is easy to check that $\sigma$
respects the relation $\bigcirc=q+q^{-1}$ and is hence well defined.
Indeed, using the neck cutting relation~\ref{eq:CuttingNecks}, then
the trace relation and then the $T$ relation we get: {\small
\[
  \eps{0.25in}{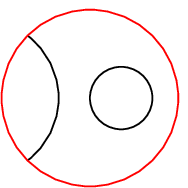}
  \overset{\sigma}{\longrightarrow} \eps{0.4in}{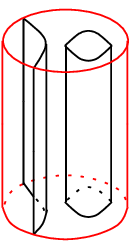}
  = \frac12\eps{0.4in}{sigma3} + \frac12\eps{0.4in}{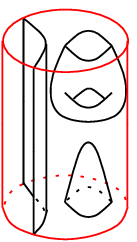}
  = \frac12(q+q^{-1})\eps{0.4in}{sigma5}
\]
\[
  = (q+q^{-1})\eps{0.4in}{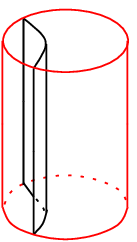}
  = (q+q^{-1})\sigma\left(\eps{0.25in}{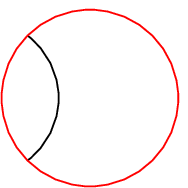}\right)
\]
}%
Now let $(S\{m\},C)$ be a general element in $\Xi_k$. Cutting necks in
$C$ using the neck cutting relation~\eqref{eq:CuttingNecks} as in
Section~\ref{subsec:OriginalKhovanov} and the relations $S$, $T$ and
$((g>1)=0)$, we can assume without loss of generality that every
connected component of $C$ has precisely one boundary component (which
is a cycle on the boundary of the cylinder). Furthermore, using the
trace relation as above, components attached to top boundary components
can be brought to the bottom (perhaps at the cost of some powers of
$q$), where they `cancel' the corresponding bottom boundary components
and create boundary-free components of $C$. And again, these can be
removed using the  $S$, $T$ and $((g>1)=0)$ relations. So without loss
of generality $C$ just has $k/2$ components with ``rectangular''
boundary, made of an arc on the top of the cylinder, an arc on the
bottom, and two arcs on the sides. As $\deg C=0$ it follows from
Definition~\ref{def:degree} that all those components of $C$ must be
disks, so $C$ is simply $S\times I$ for a cycle-free $k$--ended $S$. Hence
$\sigma$ is surjective.

We leave it to the reader to verify that the procedure described in the
previous paragraph does not depend on the choices within it (the only
apparent choice is the ordering of the necks for cutting) and hence it
defines a well-defined inverse for $\sigma$, concluding the proof of
Proposition~\ref{prop:XiIsS}.
\end{proof}

\vskip 2mm
{\bf Proof of part (2) of Theorem~\ref{thm:Euler}}\qua
\label{pf:Euler}
Given part (1) of the theorem, it is only necessary to verify part (2)
in the case of crossingless tangles (which is tautological) the
single-crossing tangles $\overcrossing$ and $\undercrossing$. The
latter two cases follow immediately from definition of the Jones
polynomial, $\overcrossing\to q\smoothing-q^2\hsmoothing$ and
$\undercrossing\to-q^{-2}\hsmoothing+q^{-1}\smoothing$, the
corresponding evaluations of $\Kh$, $\overcrossing\to
\left(\xymatrix{\underline{\smoothing\{1\}}\ar[r] &
\hsmoothing\{2\}}\right)$ and
$\undercrossing\to\left(\xymatrix{\hsmoothing\{-2\}\ar[r] &
\underline{\smoothing\{-1\}}}\right)$, the definition of the Euler
characteristic $\chi_\star$ and the identification of $\calS_k$ with
$\Xi_k$ via $\sigma$. \qed

\section{Odds and ends} \label{sec:OddsAndEnds}

\subsection{A structural conjecture}
We say that a complex in $\Kob(\emptyset)$ is ``basic'' if up to degree and
height shifts it is either one of the following two complexes:
\begin{itemize}
\item The one term complex
$\Omega_1\co\,\xymatrix{0\ar[r]&\bigcirc\ar[r]&0}$ whose only non-zero
term is a smoothing consisting of a single circle.
\item The two term complex
$\Omega_2\co\,\xymatrix{
  0\ar[r] &
  \bigcirc\ar[r]^d &
  \bigcirc\ar[r] &
  0
}$
whose two non-zero terms are both smoothings consisting of a single
circle and whose only non-zero differential $d=\eps{6mm}{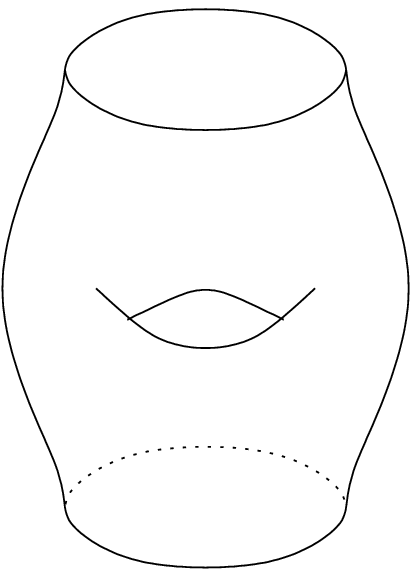}$
is the genus $1$ surface with
a circle boundary at the top and a circle boundary at the bottom.
\end{itemize}

We say that a link $L$ with $c$ components is Khovanov-basic if the
complex $\Kh(L)$ is homotopy equivalent to a direct sum of basic
complexes and the number of $\Omega_1$ terms appearing is $2^{c-1}$.

One can verify that Khovanov-basic links have Betti numbers consistent
with the knight-tetris-pawn observations of the previous section and
with Lee's Theorem~5.1 (see~\cite{Lee:AlternatingLinks}). But the Betti
numbers computed by Shumakovitch in~\cite[Appendices~A.4
and~A.5]{Shumakovitch:Torsion} show that some links are not
Khovanov-basic.

\begin{conjecture} \label{conj:Structure} Alternating links are
Khovanov-basic (and so are many other links, but we'd rather remain
uncommitted).
\end{conjecture}

\subsection{Dotted cobordisms} \label{subsec:DottedCob}
As an aside, in this section we briefly describe a weaker variant of
our theory which on links is equivalent to the original Khovanov
theory\footnote{And so it does not project to Lee's theory or to the
``new'' theory of Sections~\ref{subsec:Lee} and~\ref{subsec:4TuZ2}} but
is still rich enough for our tangles and cobordisms discussion to go
through mostly unchanged.

Extend the category $\Cob^3$ to a new category $\Cobd^3$ that has the same
objects as $\Cob^3$ and nearly the same morphisms --- the only difference
is that we now allow ``dots'' (of degree $-2$) that can be marked on
cobordisms and moved freely within each connected component of a given
cobordisms. We then form the quotient category $\Cobdl^3$ by reducing
$\Cobd^3$ modulo the local relations
\[
  \begin{array}{c}
    \includegraphics[height=1cm]{figs/S.eps}
  \end{array}\hspace{-2mm}=0,
  \qquad\qquad
  \begin{array}{c}
    \includegraphics[height=1cm]{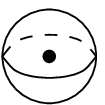}
  \end{array}\hspace{-2mm}=1,
  \qquad\qquad
  \begin{array}{c}\includegraphics[height=10mm]{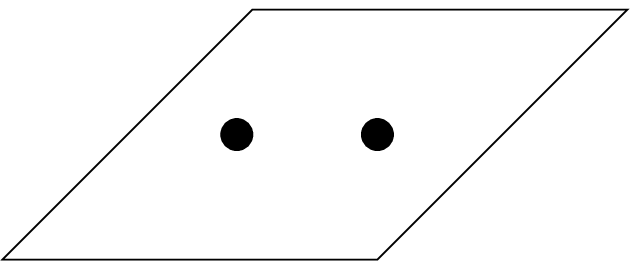}\end{array}
  \hspace{-4mm}=0,
\]
\[
  \text{and}\qquad
  \begin{array}{c}\includegraphics[height=10mm]{figs/CNN.eps}\end{array}
  =\begin{array}{c}\includegraphics[height=10mm]{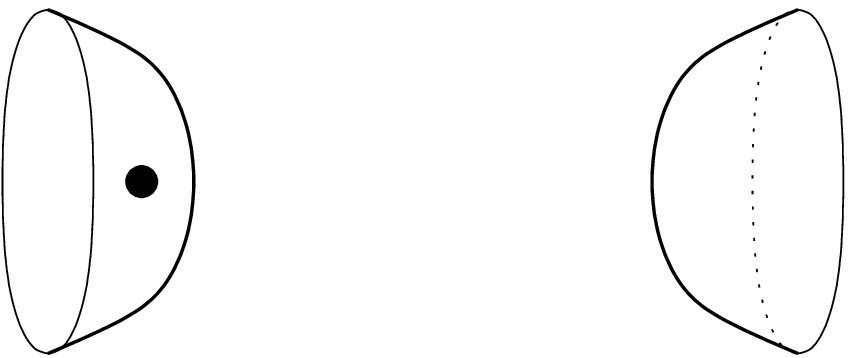}\end{array}
  +\begin{array}{c}\includegraphics[height=10mm]{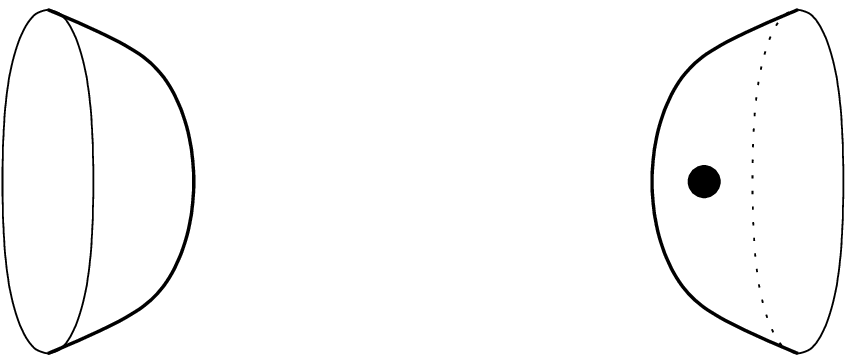}\end{array}.
\]

Pretty much everything done with $\Cobl^3$ can be repeated verbatim for
$\Cobdl^3$ (indeed, the $S$, $T$ and $\FourTu$ relations follow from the
above relations), and in particular, the appropriately modified
Theorems~\ref{thm:invariance} through~\ref{thm:Euler} hold. But now we have
a stronger neck cutting relation and tubes and handles can be removed
regardless of the ground ring. At the end we only need to consider
cobordisms in which every connected component is either a disk or a
singly-dotted disk, and these can be identified with $v_+$ and $v_-$ of
Section~\ref{sec:Homology}. Applying a tautological functor as in
Section~\ref{sec:MoreOnCobl} (with $\calO=\emptyset$) we get back to the
standard Khovanov homology without any restriction on the ground ring.

\subsection{Abstract cobordisms} \label{subsec:AbstractCobordism} As
defined, the cobordisms in the various variants of $\Cob^3$ considered
in this paper are all embedded into cylinders, and hence a full
classification of the morphisms in $\Cob^3$ could be as complicated as
knot theory itself. Two comments are in order.

\begin{itemize}

\item We could have just as well worked with un-embedded cobordisms ---
cobordisms whose boundaries are embedded in cylinders so as to have
gluing operations as required in this paper, but whose insides are
``abstract surfaces''. All the results of this paper continue to hold
in this setting as well.

\item In principle there are more embedded cobordisms than abstract
cobordisms, so keeping track of the embeddings potentially yields a
richer theory. Though if $2$ is invertible then the neck cutting
relation~\eqref{eq:CuttingNecks} shows that embedding information can
be forgotten anyway and so in this case the two theories are the same.

\end{itemize}

\subsection{Equivalent forms of the $\FourTu$ relation}
\label{subsec:4TuEquivs}
We have chosen to present the \FourTu{} relation the way we did because
in its present form it is handy to use within the proof
Theorem~\ref{thm:invariance}. But there are several equivalent or
nearly equivalent alternatives.

\begin{proposition} If the number $2$ is invertible, the neck cutting
relation of Equation~\eqref{eq:CuttingNecks} is equivalent to the
\FourTu{} relation.
\end{proposition}

\begin{proof} We've already seen in Equation~\ref{eq:CuttingNecks} that the
\FourTu{} relation implies the neck cutting relation regardless of the
ground ring. To go the other way, use $2^{-1}$ and the neck cutting
relation to cut the four tubes in the \FourTu{} relation and replace them
with 8 handles. These eight handles cancel out in pairs.
\end{proof}

\begin{proposition} The following two three-site relations are equivalent
to each other and to the \FourTu{} relation (over any ground ring):
\begin{eqnarray*}
  3S_1:\qquad& \eps{7cm}{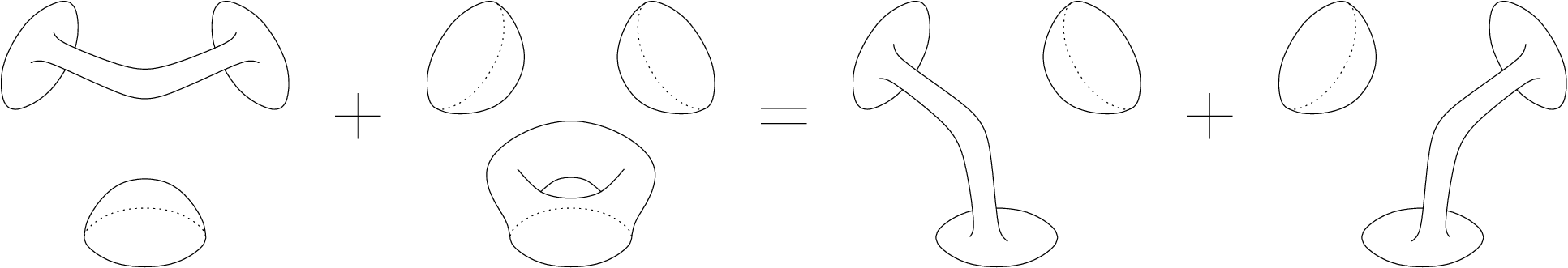} \\
  3S_2:\qquad& \displaystyle
  \sum_{\parbox{2cm}{\begin{center}\scriptsize\rm
    $0^\circ$, $120^\circ$, $240^\circ$ \newline
    rotations
  \end{center}}}
  \left(\eps{1.2cm}{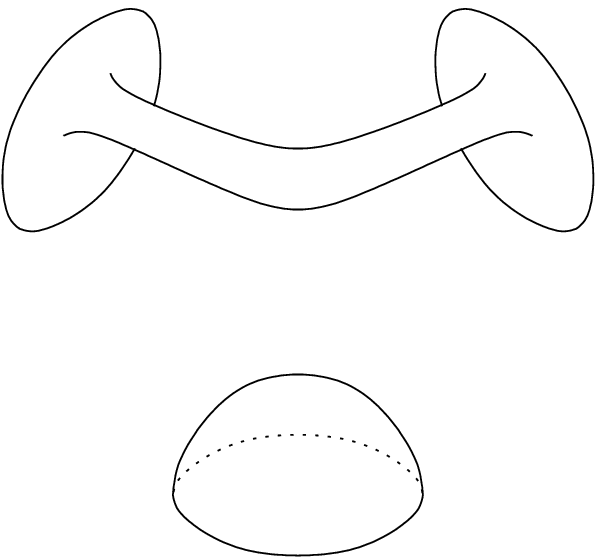}-\eps{1.2cm}{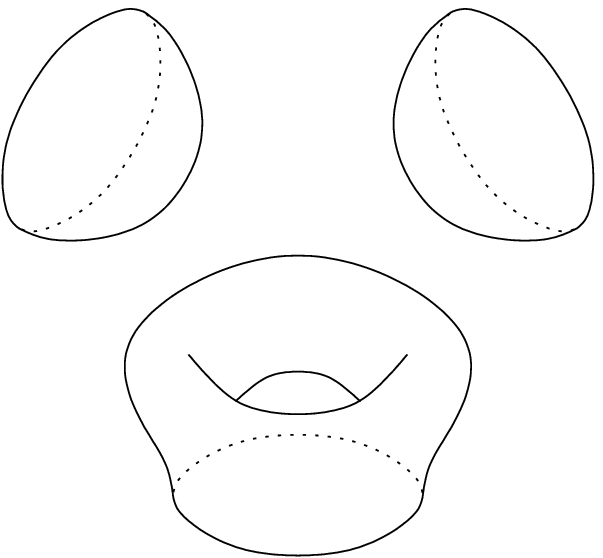}\right) = 0
\end{eqnarray*}
\end{proposition}

\begin{proof} $\FourTu\Rightarrow 3S_1$: The relation $3S_1$ is in fact
the relation used in the proof of invariance under $R1$; just as over
there, it follows from the \FourTu{} relation by specializing the general
\FourTu{} relation to the case when two of the disks $D_1$ through $D_4$
that appear in the definition of the \FourTu{} relation are on the same
connected component of the ``ambient'' cobordism.

$3S_1\Rightarrow 3S_2$: First, by putting the two upper disks in the
definition of $3S_1$ on the same connected component, we see that $3S_1$
implies the neck cutting relation. Now subtract the neck cutting relation
corresponding to the left-most tube in $3S_1$ from the entire $3S_1$
relation. The result is $3S_2$.

$3S_2\Rightarrow \FourTu$: As in the previous case, by putting the two
upper disks in the definition of $3S_2$ on the same connected
component, we see that $3S_2$ implies the neck cutting relation. Now
consider a cobordism $C$ with four special disks $D_1$ through $D_4$ as
in the definition of the \FourTu{} relation. The $3S_2$ relation
applied to $C$ at sites $1,2,3$ is $C_{12}+C_{23}+C_{13}-H_1-H_2-H_3$
where $C_{ij}$ are as in the definition of the \FourTu{} relation and
where $H_i$ is $C$ with a handle added at $D_i$. Likewise, the $3S_2$
relation applied to $C$ at sites $2,3,4$ is
$C_{23}+C_{34}+C_{24}-H_2-H_3-H_4$. The sum of these two $3S_2$
relations is $(C_{12}+C_{34}-C_{13}-C_{24}) + (2C_{13}-H_1-H_3) +
(2C_{24}-H_2-H_4) + (2C_{23}-H_2-H_3)$, and that's the sum of the
\FourTu{} relation and 3 neck cutting relations.
\end{proof}

Thus either one of $3S_1$ or $3S_2$ could have served as an alternative
foundation for our theory replacing \FourTu.

\subsection{Khovanov's $c$} \label{subsec:Khovanovc}
In his original paper~\cite{Khovanov:Categorification} on
categorification Khovanov introduced a more general knot homology
theory, defined over the polynomial ring $\bbZ[c]$ where $\deg c=2$.
The more general theory is defined using a functor, which we will call
$\calF_c$, which is similar to the functor $\calF$ of
Definition~\ref{def:Vdef}:
\begin{alignat*}{1}
  \epsilon_c: & \begin{cases} 1\mapsto v_+ \end{cases} \\
  \eta_c: & \begin{cases} v_+\mapsto -c & \\ v_-\mapsto 1 \end{cases} \\
  \Delta_c: & \begin{cases}
    v_+ \mapsto v_+\otimes v_- + v_-\otimes v_+ +c\,v_-\otimes v_-&\\
    v_- \mapsto v_-\otimes v_- &
  \end{cases} \\
  m_c: & \begin{cases}
    v_+\otimes v_-\mapsto v_- &
    v_+\otimes v_+\mapsto v_+ \\
    v_-\otimes v_+\mapsto v_- &
    v_-\otimes v_-\mapsto 0.
  \end{cases}
\end{alignat*}

This more general theory remains little studied, and unfortunately, it
doesn't fit inside our framework. The problem is that $\calF_c$ does not
satisfy the $S$ and the $\FourTu$ relations, and hence it does not descend
from $\Cob^3$ to $\Cobl^3$.

It is natural to seek for replacements for $S$ and $\FourTu$ which are
obeyed by $\calF_c$ and to try to repeat everything using those
replacements. A natural replacement for the $S$ relation may be the
relation $S_c:\ \eps{6mm}{S}=-c$, and the following relation may be a
suitable replacement for $\FourTu$ (in its $3S_2$ guise; see
Section~\ref{subsec:4TuEquivs}):
\[
  3S_c:\qquad \displaystyle
  \sum_{\parbox{2cm}{\begin{center}\scriptsize\rm
    $0^\circ$, $120^\circ$, $240^\circ$ \newline
    rotations
  \end{center}}}
  \left(\eps{1.2cm}{3S2}-\eps{1.2cm}{3S5}\right)
  = c\eps{1.2cm}{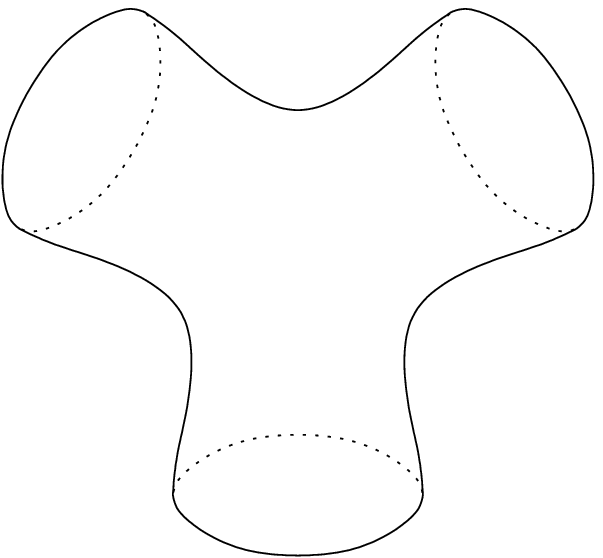}-c^2\eps{1.2cm}{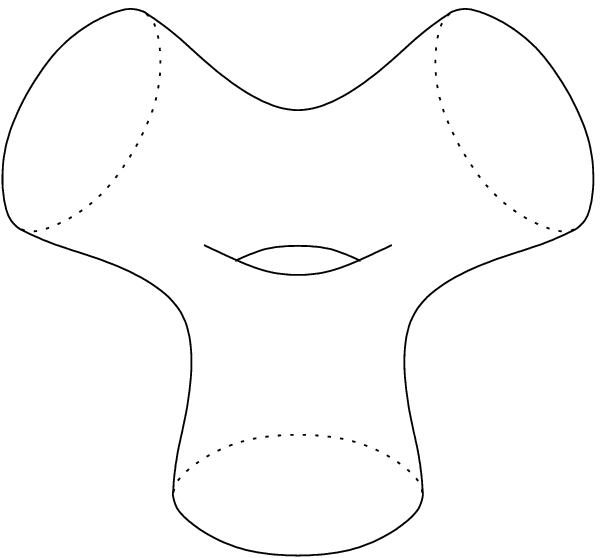}.
\]
Indeed, one may verify that $\calF_c$ satisfies $S_c$ and $3S_c$.

\begin{problem} Is there a parallel for our theory with $S_c$ and $3S_c$
replacing $S$ and $\FourTu$?
\end{problem}

We've been able to prove invariance under $R1$ using $S_c$ and $3S_c$, but
not invariance under $R2$ (one may hope that invariance under $R2$ will
imply invariance under $R3$ as in Section~\ref{subsec:Proof}, so only
invariance under $R2$ is really missing). Note, though, that the proof of
Lemma~\ref{lem:Pairings}, and hence of invariance under movie moves
(Theorems~\ref{thm:KhIsFunctor} and~\ref{thm:Main}), depend in a
fundamental way on the $S$ relation. Hence there is no reason to expect that
the theory defined by $\calF_c$ would have an invariant extension to
$\Cobi^4$. Indeed Jacobsson~\cite{Jacobsson:KhovanovConjecture} has shown
that such extension does not exist.

\subsection{Links on surfaces} 
\begin{figure}[ht!]\anchor{fig:LonS}
\[ \eps{2in}{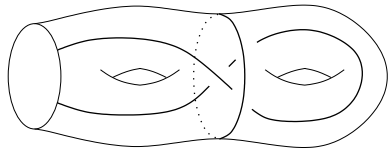} \]
\caption{
  A 3--crossing 2--component (one closed, one open) 2--ended ``tangle''
  drawn on a genus 2 surface with one boundary component $\Sigma_1^1$ (or
  embedded in $\Sigma_1^1\times I$).
}\label{fig:LonS} 
\end{figure}

Almost everything done in this paper is
local in nature, and so generalizes with no difficulty to links or
tangles drawn on surfaces (more precisely, embedded in thickened
surfaces; see \figref{fig:LonS}), in the spirit of
\cite{AsaedaPrzytyckiSikora:Surfaces}. The main challenge seems to be
to figure out the full collection of ``TQFTs'' that can be applied in
this case in order to get computable homology theories as in
Section~\ref{sec:Homology}.

\section{Glossary of notation}

We give a quick glossary of notation by section number:

\noindent
\begin{multicols}{2}
\begin{list}{}{
  \renewcommand{\makelabel}[1]{#1\hfil}
}

\item[$\llbracket\cdot\rrbracket$] the formal Khovanov bracket,
  \ref{sec:Narrative}.1.
\item[{$[\cdot]$}] height shift in complexes, \ref{subsec:Proof}.
\item[{$\{\cdot\}$}] degree shift, \ref{sec:grading}.
\item[$\fourwheel$] the cap $\fourwheel:\emptyset\to\bigcirc$,
  \ref{subsec:Proof}.
\item[$\fourinwheel$] the cup $\fourinwheel:\bigcirc\to\emptyset$,
  \ref{subsec:Proof}.
\item[$3S_i$] the three-site relations, \ref{subsec:4TuEquivs}.
\item[$\FourTu$] the four tube relation, \ref{subsubsec:STFourTu}.
\item[$\alpha$] an automorphism of a complex, \ref{subsec:CobInvProof}.
\item[$\Gamma$] the cone of a morphism, \ref{subsec:Proof}.
\item[$\partial T$] boundary of $T$, \ref{sec:Frame}.
\item[$\Delta$] a part of $\calF$, \ref{sec:Homology}.
\item[$\epsilon$] a part of $\calF$, \ref{sec:Homology}.
\item[$\eta$] a part of $\calF$, \ref{sec:Homology}.
\item[$\xi_i$] an edge of a cube, \ref{sec:Narrative}.6.
\item[$\Xi$] the trace group, \ref{subsec:TraceGroups}.
\item[$\Xi_k$] the trace group of $\Cobz^3(k)$, \ref{subsec:SkeinModules}.
\item[$\sigma$] the isomorphism $\calS_k\to\Xi_k$, \ref{subsec:SkeinModules}.
\item[$\tau$] a trace or a Lefschetz number \ref{subsec:Traces}.
\item[$\tau_\star$] the universal trace, \ref{subsec:TraceGroups}.
\item[$\Phi$] a morphism of planar algebras, \ref{subsubsec:HowTangles}.
\item[$\chi$] Euler characteristic, \ref{sec:grading}.
\item[$\chi_\tau$] $\tau$--Euler characteristic, \ref{subsec:Traces}.
\item[$\chi_\star$] the universal Euler characteristic,
  \ref{subsec:TraceGroups}.
\item[$\Psi$] a morphism of complexes, \ref{subsec:Proof}.
\item[$\Omega^r$] a chain group, \ref{sec:Frame}.
\item[$A$] an Abelian group, \ref{subsec:Traces}.
\item[$\calA$] an Abelian category, \ref{sec:Homology}.
\item[$\Alg$] the `algebraic picture', \ref{subsubsec:HowTangles}.
\item[$B$] a finite set of points on $S^1$, \ref{sec:Frame}.
\item[$c$] Khovanov's $c$, \ref{subsec:Khovanovc}.
\item[$C$] a cobordism, \ref{sec:intro}.
\item[$\calC$] a category, \ref{sec:Frame}.
\item[$\Cob^3$] either $\Cob^3(\emptyset)$ or $\Cob^3(B)$,
  \ref{sec:Frame}.
\item[$\Cob^3(\emptyset)$] 3D cobordisms, no vertical boundary,
  \ref{sec:Frame}. 
\item[$\Cob^3(B)$] 3D cobordisms, vertical boundary $B\times I$,
  \ref{sec:Frame}.
\item[$\Cobl^3$] the quotient of $\Cob^3$ by $S$, $T$ and $\FourTu$.
  Likewise for $\Cobl^3(\emptyset)$ and $\Cobl^3(B)$.
  \ref{subsubsec:STFourTu}.
\item[$\Cobz^3$] a variant of $\Cobl^3$, \ref{subsec:SkeinModules}.
\item[$\Cobd^3$] dotted cobordisms, \ref{subsec:DottedCob}.
\item[$\Cobdl^3$] dotted cobordisms modulo local relations,
  \ref{subsec:DottedCob}.
\item[$\Cob^4$] either $\Cob^4(\emptyset)$ or $\Cob^4(B)$,
  \ref{sec:EmbeddedCobordisms}.
\item[$\Cob^4(\emptyset)$] 4D cobordisms, no vertical boundary,
  \ref{sec:EmbeddedCobordisms}. 
\item[$\Cob^4(B)$] 4D cobordisms, vertical boundary $B\times I$,
  \ref{sec:EmbeddedCobordisms}.
\item[$\Cobi^4$] the quotient of $\Cob^4$ by isotopies.
  Likewise for $\Cobi^4(\emptyset)$ and $\Cobi^4(B)$.
  \ref{sec:EmbeddedCobordisms}.
\item[$d$, $d^r$] differentials, \ref{sec:Frame}.
\item[$D$] a planar arc diagram, \ref{sec:PlanarAlgebra}.
\item[$\dim_\tau$] the $\tau$--dimension, \ref{subsec:Traces}.
\item[$\dim_\star$] the universal dimension, \ref{subsec:TraceGroups}.
\item[$F$, $G$] morphisms (mainly between complexes), \ref{sec:Frame}. 
\item[$\calF$] a TQFT functor, \ref{subsubsec:HowTangles}, \ref{sec:Homology}.
\item[$\bbF_2$] the two element field, \ref{sec:MoreOnCobl}.
\item[$g$] genus, \ref{subsec:OriginalKhovanov}.
\item[$\calG_{\geq j}$] a filtration, \ref{subsec:4TuZ2}.
\item[$h$, $h^r$] homotopies, \ref{subsubsec:Homotopy}.
\item[$H$] a handle, \ref{subsec:4TuZ2}.
\item[$j$] a specific degree, \ref{subsec:4TuZ2}.
\item[$\hatJ$, $J$] the Jones polynomial, \ref{sec:Homology}, \ref{sec:Euler}.
\item[$K$] a knot, \ref{sec:Narrative}.1.
\item[$\Kh$] Khovanov homology, \ref{sec:intro}, \ref{sec:grading}.
\item[$\Kh_0$] Khovanov homology on movies, \ref{subsec:Canopolies}.
\item[$\Kob(\cdot)$] $\Kom(\Mat(\Cobl^3))$, complexes made of cobordisms,
  \ref{subsec:Statement}.
\item[$\Kobh(\cdot)$] $\Kob(\cdot)$ modulo homotopy, \ref{subsec:Statement}.
\item[$\PKob$] projectivized $\Kob$, \ref{sec:EmbeddedCobordisms}.
\item[$\PKobh$] projectivized $\Kobh$, \ref{sec:EmbeddedCobordisms}.
\item[$\Kom(\cdot)$] complexes over a category, \ref{sec:Frame}.
\item[$\Komh(\cdot)$] $\Kom(\cdot)$ modulo homotopy, \ref{sec:Frame}.
\item[$m$] a part of $\calF$, \ref{sec:Homology}.
\item[$\Mat(\cdot)$] matrices over a category, \ref{sec:Frame}.
\item[$\MM_i$] movie moves, \ref{subsec:CobInvProof}.
\item[$n$] number of crossings, \ref{sec:Narrative}.1.
\item[$n_\pm$] number of $\pm$ crossings, \ref{sec:Narrative}.2.
\item[$L$, $L_i$] links, \ref{sec:intro}.
\item[$\calO_i$] objects in a category, \ref{sec:grading}.
\item[$\calP(k)$] the sets making up an unoriented planar algebra,
  \ref{sec:PlanarAlgebra}.
\item[$\calP(s)$] the sets making up an oriented planar algebra,
  \ref{sec:PlanarAlgebra}.
\item[$q$] a formal variable, \ref{subsec:TraceGroups}.
\item[$r$] homological degree, \ref{sec:Frame}.
\item[$Ri$] Reidemeister moves, \ref{subsec:Statement}.
\item[$S$] a smoothing, \ref{subsec:SkeinModules}.
\item[$S$] the sphere relation, \ref{subsubsec:STFourTu}.
\item[$\calS_k$] a skein module, \ref{subsec:SkeinModules}.
\item[$T$] a tangle or a tangle diagram, \ref{subsubsec:WhyTangles}.
\item[$T$] the torus relation, \ref{subsubsec:STFourTu}.
\item[$\calT^0(k)$] unoriented $k$--ended tangle diagrams,
  \ref{sec:PlanarAlgebra}.
\item[$\calT^0(s)$] oriented $k$--ended tangle diagrams,
  \ref{sec:PlanarAlgebra}.
\item[$\calT(k)$] unoriented $k$--ended tangles, \ref{sec:PlanarAlgebra}.
\item[$\calT(s)$] oriented $k$--ended tangles, \ref{sec:PlanarAlgebra}.
\item[$\Top$] the `topological picture', \ref{subsubsec:HowTangles}.
\item[$v_\pm$] the generators of $V$, \ref{sec:Homology}.
\item[$V$] the $\bbZ$--module for a circle, \ref{sec:Homology}.
\item[$X$] a crossing, \ref{sec:PlanarAlgebra}, \ref{subsec:CobInvProof}.
\item[$\Ztwo$] $\bbZ$ localized at $2$, \ref{subsec:OriginalKhovanov}.
\item[$\ZMod$] $\bbZ$--modules, \ref{sec:Homology}.

\end{list}
\end{multicols}

\end{document}